\documentclass[oneside,english]{amsart}
\usepackage{tgtermes}
\usepackage[T1]{fontenc}
\usepackage[latin9]{inputenc}
\usepackage{color}
\usepackage{babel}
\usepackage{mathrsfs}
\usepackage{mathtools}
\usepackage{amstext}
\usepackage{amsthm}
\usepackage{amssymb}
\usepackage{graphicx}
\usepackage[unicode=true,pdfusetitle,
 bookmarks=true,bookmarksnumbered=false,bookmarksopen=true,bookmarksopenlevel=2,
 breaklinks=false,pdfborder={0 0 0},pdfborderstyle={},backref=false,colorlinks=true]
 {hyperref}

\makeatletter
\numberwithin{equation}{section}
\numberwithin{figure}{section}
\theoremstyle{plain}
\newtheorem{thm}{\protect\theoremname}
\theoremstyle{definition}
\newtheorem{defn}[thm]{\protect\definitionname}
\theoremstyle{remark}
\newtheorem{rem}[thm]{\protect\remarkname}
\theoremstyle{plain}
\newtheorem{prop}[thm]{\protect\propositionname}
\theoremstyle{plain}
\newtheorem{cor}[thm]{\protect\corollaryname}
\theoremstyle{definition}
\newtheorem{example}[thm]{\protect\examplename}
\theoremstyle{plain}
\newtheorem{lem}[thm]{\protect\lemmaname}

\makeatother

\providecommand{\corollaryname}{Corollary}
\providecommand{\definitionname}{Definition}
\providecommand{\examplename}{Example}
\providecommand{\lemmaname}{Lemma}
\providecommand{\propositionname}{Proposition}
\providecommand{\remarkname}{Remark}
\providecommand{\theoremname}{Theorem}

\begin{document}
\title{Mabuchi Solitons and Relative Ding Stability of Toric Fano Varieties}
\author{Yi Yao}
\begin{abstract}
As a generalization of K\"{a}hler-Einstein metrics for Fano manifolds
with nonvanishing Futaki invariant, Mabuchi solitons are critical
points of a Calabi-type energy functional. We study their existence
on toric Fano varieties and the underlying algebraic stability notion:
relative Ding stability. As a toy model for a YTD type correspondence,
a new feature is the emergence of a non-uniformly stable case. We
show a partial coercivity for the modified Ding functionals in this
case, and obtain singular Mabuchi solitons via a variational approach.
In the unstable case, we determine the maximal destabilizer which
is a simple convex function over the moment polytope, and establish
a Moment-Weight equality which connects the infimum of a Calabi-type
energy and the Berman-Ding invariant.
\end{abstract}

\email{yeeyoe@163.com }
\address{SCHOOL OF MATHEMATICS, HUNAN UNIVERSITY, CHANGSHA, 410012, CHINA}

\maketitle
\tableofcontents{}

\section{Introduction}

For a polarized K\"{a}hler manifold $(X,L)$, a central problem is
the existence of constant scalar curvature K\"{a}hler (cscK) metrics
in class $c_{1}(L)$. The Yau-Tian-Donaldson conjecture predicts the
existence is equivalent to K-polystability of $(X,L)$. In recent
years, there are great progresses toward this conjecture, such as
the confirmation for the Fano case, where $L=c_{1}(X)$ and cscK metrics
are K\"{a}hler-Einstein (KE) metrics, see \cite{CDS,Tian CPAM,BBJ},
etc.

Recall that cscK metrics can be realized as zeros of the moment map
in an infinitely dimensional GIT (geometric invariant theory) model
discovered by Fujiki \cite{Fujiki} and Donaldson \cite{Dona Fields}.
The associated Kempf-Ness function is the Mabuchi functional (or K-energy)
$M(\phi)$. K-stability was introduced by Tian (extended by Donaldson)
as an analogue of Hilbert-Mumford's criterion in GIT, via the limit
slope of $M(\phi)$ along Bergman geodesics. 

As a generalization of cscK metrics, the extremal metrics in $c_{1}(L)$
are defined as the critical points of the $L^{2}$-norm of the moment
map 
\[
Cal(\omega)=\frac{1}{V}\int_{X}\left(S(\omega)-\bar{S}\right)^{2}\omega^{n},\ \omega\in c_{1}(L),\ V=[\omega]^{n},
\]
called the Calabi energy, where $S(\omega)$ is the scalar curvature.
Their existence is expected to be equivalent to relative K-stability,
see \cite{Gabor thesis}. In \cite{Dona Ding}, Donaldson found a
new GIT model specialized for Fano manifolds, where the moment map
is
\[
c_{1}(X)\ni\omega\longmapsto1-e^{h_{\omega}},
\]
where $h_{\omega}$ is the Ricci potential satisfying $Ric(\omega)-\omega=i\partial\bar{\partial}h_{\omega}$
and $\int_{X}e^{h_{\omega}}\omega^{n}=V$. The associated Kempf-Ness
function is the Ding functional $D(\phi)$. Considering its limit
slopes along geodesic rays, it gives rise to the notion of Ding stability,
see \cite{Berman}. 

The starting point of this paper is to consider the analogue of Calabi's
energy in the new model, 
\[
Din(\omega)\coloneqq\frac{1}{V}\int_{X}\left(1-e^{h_{\omega}}\right)^{2}\omega^{n},\ \omega\in c_{1}(X),
\]
we call it the \textsl{Ding energy}. Our first observation (Theorem
\ref{thm: critical pt of DinE}) is that the critical points of $Din(\omega)$
are exactly the generalized KE metrics introduced by Mabuchi \cite{Mabu soliton}
for Fano manifolds with nonvanishing Futaki invariant. They are defined
to be the metrics $\omega\in c_{1}(X)$ such that 
\[
Z=\textrm{grad}_{\omega}^{(1,0)}(1-e^{h_{\omega}})\coloneqq\sum_{i,j}g^{i\bar{j}}\frac{\partial(1-e^{h_{\omega}})}{\partial\bar{z}^{j}}\frac{\partial}{\partial z^{i}}
\]
is a holomorphic vector field. We call them \textsl{Mabuchi solitons}
(comparing to K\"{a}hler-Ricci solitons where $\textrm{grad}_{\omega}^{(1,0)}h_{\omega}$
is holomorphic). Mabuchi \cite{Mabu soliton,Mabu 2 Vector field,Mabu 3 Calabi-Matsu,Mabu 4 Multiplier Herm}
shows they share many properties with extremal metrics. By \cite{Mabu soliton},
if $Z$ is holomorphic, it must coincide with the extremal vector
field defined in \cite{Futaki Mabuchi}, which can be determined by
the Futaki invariant after choosing a maximal compact subgroup of
$\textrm{Aut}^{0}(X)$, see Remark \ref{rem: extrm vec}. From the
viewpoint of PDE, Mabuchi soliton $\omega_{u}=\omega_{r}+i\partial\bar{\partial}u$
satisfying a Monge-Amp\`{e}re type equation 
\begin{equation}
\left(1-\theta_{Z}(\omega_{r})-Z(u)\right)\omega_{u}^{n}=e^{-u+h_{r}}\omega_{r}^{n},\label{Intro,MA equat}
\end{equation}
where $\theta_{Z}(\omega_{r})$ is the Hamiltonian function of $Z$
with respect to a reference metric $\omega_{r}$. Clearly, if there
is a smooth solution then $\sup_{X}\left(\theta_{Z}(\omega_{r})+Z(u)\right)<1$.
By \cite{Futaki Mabuchi}, this supremum is an invariant only depending
on $X$. 

As with extremal metrics, the existence of Mabuchi solitons should
correspond to some stability notion. The new GIT model hints us it
should be the relative version of Ding stability. In this paper, as
in \cite{Dona toric}, we introduce \textsl{relative Ding stability}
for toric Fano varieties and study its relation to the existence of
Mabuchi solitons. The full name should be relative Ding polystability.
In this paper, for simplicity, we call it relative Ding stability
or relative D-stability. In the toric setting, by assuming torus symmetry,
complex objects in K\"{a}hler geometry can be reduced to real objects
in convex analysis on $\mathbb{R}^{n}$. Then with the tools from
real Monge-Amp\`{e}re equations, problems usually become more tractable. 

The subsequent discussions are based on the following setting: $X_{P}$
is a toric Fano variety associated to polytope $P$, which is dual
to a Fano polytope (see Section \ref{subsec:Toric-Fano-Varieties}).
Actually, many results are valid for a general convex body $P$ containing
$0$ in its interior. 

\subsection*{Limit slopes of the Ding functional}

In \cite{Berman}, Berman obtained a lct (log canonical threshold)
formula for the limit slope of $D(\phi)$ along the geodesic ray induced
from a test-configurations $(\mathcal{X},\mathcal{L})$, the result
is called the Berman-Ding invariant denoted by $D^{\textrm{NA}}(\mathcal{X},\mathcal{L})$,
see Section \ref{subsec: Berman lct formula}. Ding stability is defined
in terms of this invariant. As with Ding stability, we will define
relative Ding stability in terms of the limit slopes of the modified
Ding functional $D_{\ell}(\phi)$ along toric geodesic rays. So first
we need to derive an explicit formula for the limit slopes in the
toric setting. 

Let $T_{N}$ ($S_{N}$) be the complex (compact) torus acting on $X$.
Let $\textrm{PSH}(-K_{X})^{S}$ be the space of $S_{N}$-invariant
psh (pluri-subharmonic) metrics on $-K_{X}$. By \cite{Guedj toric potential},
it can be identified with 
\[
\textrm{PSH}(P)=\{\phi:\mathbb{R}^{n}\rightarrow\mathbb{R}\ \textrm{convex}\mid\phi\leq\psi_{P}+C\ \textrm{on}\ \mathbb{R}^{n}\ \mathrm{for\ some}\ C\},
\]
where $\psi_{P}(x)=\sup_{p\in P}\left\langle x,p\right\rangle $ is
the support function of $P$. Its subspace $\textrm{PSH}_{b}(-K_{X})^{S}$
of bounded psh metrics can be identified with $\textrm{PSH}_{b}(P)\subset\textrm{PSH}(P)$
consisting of $\phi$ such that $\left|\phi-\psi_{P}\right|$ is bounded
over $\mathbb{R}^{n}$. The space of finite energy psh metrics $\mathcal{E}^{1}(-K_{X})^{S}$
can be identified with
\[
\mathcal{E}^{1}(P)=\{\phi\in\textrm{PSH}(P)\mid\int_{P}\phi^{*}dy<\infty\},
\]
where $\phi^{*}:P\rightarrow\mathbb{R}\cup\{\infty\}$ is the Legendre
dual. Restricting on $\mathcal{E}^{1}(P)$, the Ding functional is
given by (\ref{eq: toric Ding}). For various spaces of psh metrics
and the identifications, see Section \ref{sec: Dictionary}. 
\begin{thm}[Theorem \ref{thm: lim slope of toric geo ray}]
 Let $X$ be the toric Fano variety given by a polytope $P$. For
$\phi_{0}\in\textrm{PSH}_{b}(P)$ and convex $f\in L^{1}(P)$, let
$0\leq t\mapsto\phi_{t}=\left(\phi_{0}^{*}+tf\right)^{*}\in\mathcal{E}^{1}(P)$
be the associated toric geodesic ray with finite energy. Then the
limit slope of the Ding functional satisfies 
\begin{equation}
\lim_{t\rightarrow\infty}\frac{1}{t}D(\phi_{t})=D^{\mathrm{NA}}(f)\coloneqq-f(0)+\frac{1}{V_{P}}\int_{P}fdy.\label{Intro, limit slop of D}
\end{equation}
\end{thm}

Hence if $f$ is a rational PL (piecewise linear) convex function
on $P$, and $(\mathcal{X},\mathcal{L})$ is the toric test-configuration
for $(X,-K_{X})$ associated to datum $(P,f,L)$, where $L$ is a
large integer, then the Berman-Ding invariant $D^{\textrm{NA}}(\mathcal{X},\mathcal{L})$
should be equal to $D^{\textrm{NA}}(f)$. We verify this by Berman's
lct formula, see Theorem \ref{thm: NA Ding of toric TC}. It is interesting
to compare (\ref{Intro, limit slop of D}) with the formula of the
non-Archimedean Mabuchi functional $M^{\textrm{NA}}$ obtained by
Donaldson in \cite{Dona toric} (LHS below). We have 
\begin{equation}
\frac{1}{V_{P}}(\int_{\partial P}fd\sigma-n\int_{P}fdy)\geq-f(0)+\frac{1}{V_{P}}\int_{P}fdy,\label{eq: Intro, D vs M}
\end{equation}
which is a special case of the comparison between these two invariants
in the general setting. By the seminal work \cite{Wang-Zhu}, a toric
Fano manifold $X_{P}$ admits KE metrics if and only if its Futaki
invariant vanishes, and then is equivalent to that the origin is the
barycenter of $P$. By Jensen's inequality, (\ref{eq: Intro, D vs M})
clearly shows that if the origin is the barycenter of $P$ then $P$
is K-polystable. See Section \ref{subsec: D vs K} for more discussions. 

\subsection*{Relative Ding stability}

There is a unique affine function $\ell_{P}$ (we call the \textsl{Ricci
affine function}) such that 
\[
-h(0)+\int_{P}h\cdot\ell_{P}dy=0,
\]
for all affine functions $h$. Sometimes we denote it by $\ell$.
Then the Hamiltonian function of the extremal vector field $Z$ can
be expressed by $\ell$. The modified Ding functional $D_{\ell}$
takes Mabuchi solitons as its critical points. The expression is 
\[
D_{\ell}(\phi)=-\log\int_{\mathbb{R}^{n}}e^{-\phi}+\int_{P}\phi^{*}\ell dy,\ \phi\in\mathcal{E}^{1}(P).
\]
It is convex along toric geodesics. By (\ref{Intro, limit slop of D}),
its limit slopes along toric geodesic rays are 
\[
D_{\ell}^{\textrm{NA}}(f)\coloneqq-f(0)+\int_{P}f\ \ell dy,
\]
we call the relative Berman-Ding invariant. Then we define the relative
Ding stability of $X_{P}$ in terms of the sign of $D_{\ell}^{\textrm{NA}}$
on convex functions, see Definition \ref{def: rD-stable and unstable}
and \ref{def: Uniform rD-stable}. Finally, it turns out this stability
can be completely detected by the sign of $\inf_{P}\ell$, see Corollary
\ref{cor: summary of rD-stability}. In the following statements,
``D-stable'' is short for Ding stable, 
\begin{enumerate}
\item $X_{P}$ is relatively D-stable if and only if $\inf_{P}\ell_{P}\geq0$. 
\item $X_{P}$ is (non-)uniformly relative D-stable if and only if $\inf_{P}\ell_{P}>0$
$(=0)$. 
\item $X_{P}$ is relatively D-unstable if and only if $\inf_{P}\ell_{P}<0$. 
\end{enumerate}
Each of these cases does happen, see Example \ref{exa: 2 polygons}.
We find that $\mathbb{P}(1,1,2)$ is non-uniformly relative D-stable,
but the smooth such examples have not been found yet. In addition,
the relative D-semistability turns out to be equivalent to relative
D-stability. 

The situation is much more complicated for relative K-stability. Until
now, we only know some sufficient conditions for the relative K-stability
of toric varieties, see \cite{Zhou-Zhu Adv,Wang-Zhou,Yotsutani ZhouB},
etc. And to the author's knowledge, there is no example yet which
is relative K-stable but not uniformly. There are some implication
relations between relative K-stability and relative D-stability, see
Section \ref{subsec: rD vs rK}.

\subsection*{The uniformly stable case}

Corollary \ref{cor: unif D-stable} shows that $\inf_{P}\ell>0$ is
equivalent to the existence of a $\delta>0$ such that $D_{\ell}^{\textrm{NA}}(f)\geq\delta J^{T}(f)$
for all rational PL convex $f$, where 
\[
J^{T}(f)=\inf\left\{ \frac{1}{V_{P}}\int_{P}(f-l)dy-\inf_{P}(f-l)\middle|\mathrm{affine}\ l\right\} 
\]
is the reduced non-Archimedean J-functional introduced by Hisamoto
in \cite{Hisa toric}. It measures how far $f$ is from being affine. 

In this case, we show Mabuchi solitons exist. It needs to solve a
Monge-Amp\`{e} type equation on $\mathbb{R}^{n}$, 
\begin{equation}
\textrm{MA}(\phi)=\frac{e^{-\phi}}{\ell(d\phi)},\ C^{-1}\leq\ell\leq C\ \mathrm{on}\ P,\label{Intro, MS equation}
\end{equation}
where $\phi$ is required to be smooth and strictly convex on $\mathbb{R}^{n}$,
and belonging to $\textrm{PSH}_{b}(P)$. This problem has been solved
by Berndtsson-Berman \cite{Berman Berndtsson} via a variational approach.
When $X$ is smooth, it also can be solved via the continuation method
in \cite{Wang-Zhu}. We will use these results in Section \ref{sec: Exist result for unifm}. 

\subsection*{The non-uniformly stable case }

When $\inf_{P}\ell=0$, there is no $\delta>0$ such that $D_{\ell}^{\textrm{NA}}(f)\geq\delta J^{T}(f)$
for all rational PL convex $f$. We only have $D_{\ell}^{\textrm{NA}}(f)\geq0$
and the equality holds only for affine $f$. In this case, equation
(\ref{Intro, MS equation}) has not been considered in the previous
studies. Note that $\inf_{P}\ell=0$ amounts to $\sup_{X}\theta_{Z}=1$
for equation (\ref{Intro,MA equat}), so there is no smooth solution.
However, in Example \ref{exa: 8c conic MS}, where $X=\mathbb{P}(1,1,2)$
satisfying $\inf_{P}\ell=0$, we find an explicit solution for (\ref{Intro, MS equation})
which gives a metric with conic singularities along the divisor $\triangle_{\ell}$
corresponding to the facet $P\cap\{\ell=0\}$, where $\triangle_{\ell}$
is contained in the regular locus of $X$. Comparing with the smooth
solutions, a feature of it is that the term $\omega_{u}^{n}$ in (\ref{Intro,MA equat})
has a pole along $\triangle_{\ell}$. 

We show the existence of such solutions in general cases, following
the variational approach in \cite{Berman Berndtsson}. Although by
Corollary \ref{cor: unif D-stable}, $D_{\ell}$ is not coercive when
$\inf_{P}\ell=0$, but we can show $D_{\ell}$ has a partial coercivity
(see below). Combing with some tools, e.g. weighted Poincar\'{e}
inequality, the partial coercivity is sufficient for us to repeat
the arguments in \cite{Berman Berndtsson}. 
\begin{thm}[Theorem \ref{thm: Partial coercive}, \ref{thm: Exist nonunifm}]
 \label{Intro thm 2} Let $P\subset\mathbb{R}^{n}$ be a polytope
containing $0$ in its interior. We assume $\ell_{P}\geq0$ on $P$. 

(1) (Partial coercivity) For any $\epsilon\in(0,1)$, there is a constant
$C_{\epsilon}>0$ depending on $P$ and $\epsilon$ such that $D_{\ell}(\phi)\geq\epsilon J_{\ell}^{T}(\phi)-C_{\epsilon}$
for all $\phi\in\mathcal{E}^{1}(P)$. Where $J_{\ell}^{T}$ is the
modified reduced J-functional, see Definition \ref{def:def of reduced J}. 

(2) (Existence) There exists $\phi\in\textrm{PSH}_{b}(P)$ satisfying
$\mathrm{MA}_{\ell}(\phi)=e^{-\phi}$ on $\mathbb{R}^{n}$ in the
sense of Alexandrov, namely 
\[
\int_{\partial\phi(E)}\ell=\int_{E}e^{-\phi}
\]
for all Borel set $E\subset\mathbb{R}^{n}$. Its Legendre dual $\phi^{*}$
is H\"{o}lder continuous on $P$ for any exponent $\gamma\in(0,1)$.
When $\dim P\cap\{\ell=0\}\leq\frac{n}{2}$ (including $\inf_{P}\ell>0$),
we know $\phi$ is smooth and strictly convex over $\mathbb{R}^{n}$,
and $d\phi$ induces a diffeomorphism between $\mathbb{R}^{n}$ and
$\mathrm{Int}(P)$. 
\end{thm}

The key step in the proof of part (1) is to construct a subsolution
$\phi\in\textrm{PSH}_{b}(P)$ such that $C\cdot\ell(d\phi)\mathrm{MA}(\phi)\geq e^{-\phi}$
on $\mathbb{R}^{n}$ for some constant $C$. In Lemma \ref{lem: subsolut},
for polytopes, we find a Fubini-Study type potential satisfying this
condition. We expect such subsolutions to always exist for any convex
body. If it does, the above theorem could be extended to general convex
bodies. 

The solution $\phi$ obtained in part (2) is expected to be always
smooth and strictly convex on $\mathbb{R}^{n}$, we only verify this
when the zero locus of $\ell$ on $P$ is not too ``big''. The difficulty
comes from $\inf_{P}\ell=0$. By the equation $\textrm{MA}(\phi)=\frac{e^{-\phi}}{\ell(d\phi)}$,
we can not directly bound the RHS from above over compact sets. The
available results for the interior regularity for Monge-Amp\`{e}re
equations usually require the upper bound of the RHS. The above partial
result follows from a contradiction argument and a result of Mooney
\cite{Mooney} (a refinement of a result due to Caffarelli \cite{Caffarelli}).
For the boundary asymptotics of $\phi$, we expect it to satisfy a
modified Guillemin's condition, see Remark \ref{rem: Modified Guillemin}. 

This problem also relates to the optimal transport, since $d\phi$
gives an optimal transport plan from the measure $e^{-\phi}dx$ to
$\ell dy$. But the existing results (e.g. \cite{Philip}) for the
regularity of the optimal transport maps always assume the target
density function is bounded below by a positive number, so $\ell dy$
does not satisfy this when $\inf_{P}\ell=0$. 

\subsection*{The unstable case and the moment-weight equality}

Recall the setting of GIT, let $M\subset\mathbb{P}^{N}$ be a projective
manifold endowing with the Fubini-Study metric. Suppose $M$ is preserved
by the linear action of a compact Lie subgroup $G\subset U(N+1)$.
Its complexification $G^{c}$ also acts on $M$. We have a canonical
moment map $\mu:M\rightarrow\mathfrak{g}^{*}$ where $\mathfrak{g}=\mathrm{Lie}(G)$.
Then an orbit $G^{c}x\subset M$ is polystable if and only if it contains
a zero of $\mu$. For the unstable orbits, we have the \textit{moment-weight
inequality}, see \cite{M-W inequ}. For any point $x\in M$, we have
\begin{equation}
\inf_{g\in G^{c}}\left|\mu(gx)\right|\geq\sup_{\xi\in\mathfrak{g}\backslash\{0\}}\frac{-w(x,\xi)}{\left|\xi\right|},\label{eq: M-W ineq GIT}
\end{equation}
where $w(x,\xi)=\lim_{t\rightarrow\infty}\left\langle \mu\left(\exp(it\xi)x\right),\xi\right\rangle $
is the Hilbert-Mumford weight. Both sides of this inequality measure
how far the orbit $G^{c}x$ is from being stable. 

In the GIT model for cscK metrics \cite{Fujiki,Dona Fields}, the
analogue of inequality (\ref{eq: M-W ineq GIT}) is obtained by Donaldson
in \cite{Dona Lower bound}. For any polarized K\"{a}hler manifold
$(X,L)$, we have 
\begin{equation}
\inf_{\omega\in c_{1}(L)}\left\Vert S(\omega)-\bar{S}\right\Vert _{L^{2}(\omega^{n})}\geq\sup_{(\mathcal{X},\mathcal{L})}\frac{-\mathrm{DF}(\mathcal{X},\mathcal{L})}{\left\Vert \mathcal{X},\mathcal{L}\right\Vert _{2}},\label{eq: M-W general L}
\end{equation}
where the supremum is taken over all test-configurations for $(X,L)$
and $\mathrm{DF}(\mathcal{X},\mathcal{L})$ is the Donaldson-Futaki
invariant. For toric manifolds, Sz\'{e}kelyhidi \cite{Gabor optimal}
shows that it is an equality if the Calabi flow exists for long time. 

In the new GIT model for KE metrics on Fano manifolds \cite{Dona Ding},
the analogue of inequality (\ref{eq: M-W ineq GIT}) is 
\begin{equation}
\inf_{\omega\in c_{1}(X)}\left\Vert 1-e^{h_{\omega}}\right\Vert _{L^{2}(\omega^{n})}\geq\sup_{(\mathcal{X},\mathcal{L})}\frac{-D^{\textrm{NA}}(\mathcal{X},\mathcal{L})}{\left\Vert \mathcal{X},\mathcal{L}\right\Vert _{2}},\label{eq: M-W Fano}
\end{equation}
where $(\mathcal{X},\mathcal{L})$ are test-configurations for $(X,-K_{X})$
and the LHS is $Din(\omega)^{1/2}$. Since $D^{\textrm{NA}}$ is defined
as the limit slope of the Ding functional, (\ref{eq: M-W Fano}) directly
follows from the convexity of the Ding functional, see \cite{Hisa Crelle}
and Theorem 4.3 in \cite{Berman}. 

We show that (\ref{eq: M-W Fano}) takes equality for toric Fano manifolds.
Before that, we study the optimization problem proposed by the RHS
of (\ref{eq: M-W Fano}). When $(\mathcal{X},\mathcal{L})$ is a toric
test-configuration given by a rational PL convex $f$, the quotient
in (\ref{eq: M-W Fano}) is reduced to 
\[
\frac{D^{\textrm{NA}}(f)}{\left\Vert f-\bar{f}\right\Vert _{2,P}},\ \bar{f}=\frac{1}{V_{P}}\int_{P}fdy.
\]
We completely determine the minimizer, which turns out to be a simple
convex function and easy to compute. 
\begin{thm}[Theorem \ref{thm: optimal balan density}, \ref{thm: M-W equality}]
 (1) Let $P\subset\mathbb{R}^{n}$ be a convex body containing $0$
in its interior. When $\inf_{P}\ell_{P}<0$, there is a unique affine
function $\mathfrak{a}$ such that $P\cap\{\mathfrak{a}\geq0\}$ contains
$0$ in its interior and $\ell_{P\cap\{\mathfrak{a}\geq0\}}$ equals
to $\mathfrak{a}$ itself, see (\ref{eq: optimal affine}) for more
specific condition. The convex body $P\cap\{\mathfrak{a}\geq0\}$
is non-uniformly relative D-stable. 

(2) The maximal destabilizer $\mathfrak{b}$ is defined to be $\max\{0,\mathfrak{a}\}$
when $\inf_{P}\ell_{P}<0$, and $\mathfrak{b}=\ell_{P}$ when $\inf_{P}\ell_{P}\geq0$.
Then we have 
\[
\inf\left\{ \frac{D^{\mathrm{NA}}(f)}{\left\Vert f-\bar{f}\right\Vert _{2,P}}\middle|\mathrm{convex}\ f\in L^{2}(P)\right\} =-\left\Vert \mathfrak{b}-V_{P}^{-1}\right\Vert _{2,P}
\]
 and the infimum is attained by $f=\mathfrak{b}$. 

(3) Let $X$ be a toric Fano manifold associated to Delzant polytope
$P$. Then (\ref{eq: M-W Fano}) takes equality, and both sides are
equal to $\left(V_{P}\cdot\mathfrak{b}(0)-1\right)^{1/2}$. 
\end{thm}

The maximal destabilizer $\mathfrak{b}$ is also called the \textit{optimal
balancing density} of $P$, since it is the density function with
the minimum $L^{2}$-norm making the origin is the barycenter of $P$.
It can be determined by solving a polynomial system (\ref{eq: optimal affine}).
The gradient of $\mathfrak{a}$ may not be the scalar multiples of
a lattice vector, see Example \ref{exa: optimal 6d}. In that case,
$P\cap\{\mathfrak{a}\geq0\}$ can not be the moment polytope of a
toric variety. Since $P\cap\{\mathfrak{a}\geq0\}$ satisfies the assumption
of Theorem \ref{Intro thm 2}, it gives us a weak solution on this
subpolytope. In Example \ref{exam: optimal 3fold}, where $X$ is
the Fano 3-fold $\mathbb{P}\left(\mathcal{O}_{\mathbb{P}^{2}}\oplus\mathcal{O}_{\mathbb{P}^{2}}(2)\right)$
which is relatively D-unstable. The part (3) of the above theorem
yields 
\[
\inf_{\omega\in c_{1}(X)}\frac{1}{V}\int_{X}\left(1-e^{h_{\omega}}\right)^{2}\omega^{n}\approx0.6878...
\]

The analogous problem for K-stability had been studied by Sz\'{e}kelyhidi
in \cite{Gabor optimal}. There also exists a unique maximal destabilizer,
but we know a little about its structure yet, e.g. whether it is piecewise
linear or even bounded. This kind of maximal destabilizers can be
seen as the analogue in the manifold situation to the Harder-Narasimhan
filtration of an unstable vector bundle. In addition, since the long
time existence of Calabi flow have not been established yet, the problem
that whether (\ref{eq: M-W general L}) is an equality is still open,
see the recent relevant work \cite{Xia}. 

After the first version of this paper, there are many relevant studies
about Mabuchi solitons \cite{Na,Liyan Zhou,collins,Hisa MS D,Yao19,Han Li}
and Moment-Weight equalities \cite{collins,Xia,Hisa M-W}, etc. Specially
in \cite{collins}, they considered the gradient flow of the Ding
functional $\frac{\partial}{\partial t}\phi_{t}=1-e^{h(\phi_{t})}$.
For the relatively D-unstable toric Fano manifolds, they show that
$e^{h(\phi_{t})}$ (as a function on $P$) will converge to $V_{P}\cdot\mathfrak{b}$
in $L^{2}$-norm. The asymptotic behavior of $(X,\omega_{\phi_{t}})$
is not very understood yet. Such as whether it converges (in the Gromov-Hausdorff
sense) to a Fano variety with a singular Mabuchi soliton, as with
the Hamilton-Tian conjecture predicts for the K\"{a}hler-Ricci flow.
If it does, the limit metric may be related to the weak solution on
$P\cap\{\mathfrak{a}\geq0\}$ given by Theorem \ref{Intro thm 2}.
Since $P\cap\{\mathfrak{a}\geq0\}$ may not be a moment polytope,
we can not expect the limit space to be also a toric variety. We leave
these questions for the future studies.

\subsection*{Notations}

A convex body means a compact convex subset of $\mathbb{R}^{n}$ with
nonempty interior. $N$ is a lattice with a fixed basis, $M$ is the
dual lattice. The coordinates on $N_{\mathbb{R}}$ are denoted by
$(x_{i})$, and $(y_{i})$ for $M_{\mathbb{R}}$. $T_{N}=N\otimes_{\mathbb{Z}}\mathbb{C}^{*}$
is the complex torus, $S_{N}=N\otimes_{\mathbb{Z}}\mathbb{S}^{1}$
is the compact torus. $Q$ is a Fano polytope in $N_{\mathbb{R}}$,
its dual polytope is $P\subset M_{\mathbb{R}}$. The associated toric
Fano variety is denoted by $X_{P}$. We use $\phi$, $\psi$, ...
to denote the psh metrics on line bundles and also the induced convex
functions on $N_{\mathbb{R}}$. Their Legendre duals are denoted by
$\phi^{*}$, $\psi^{*}$, ... The functions on $P$ are denoted by
$f$, $g$, ... The Ricci affine function associated to a convex body
$P$ is always denoted by $\ell_{P}$ (sometimes by $\ell$ when $P$
is clear). 

\subsection*{Acknowledgments}

The author would like to thank An-Min Li, Li Sheng, Ya-Long Shi, and
Satoshi Nakamura for helpful discussions. Specially thanks to Naoto
Yotsutani and Bin Zhou for valuable suggestions and discussions. We
use \textsl{Wolfram Mathematica} to compute the examples and \textsl{GeoGebra}
to make the figure. The author is supported by NSFC (No. 751203123)
and Fundamental Research Funds for the Central Universities (No. 531118010149). 

\section{Mabuchi solitons and the Ding energy}

In this section, we briefly recall the works of Mabuchi \cite{Mabu soliton,Mabu 2 Vector field,Mabu 3 Calabi-Matsu,Mabu 4 Multiplier Herm},
then show that Mabuchi solitons are exactly the critical points of
the Ding energy. 

Let $X$ be a Fano manifold with a smooth reference metric 
\[
\omega_{r}=\sqrt{-1}\sum_{i,j}g_{i\bar{j}}dz^{i}\wedge d\bar{z}^{j}\in2\pi c_{1}(X).
\]
We set $V=\int\omega_{r}^{n}=(2\pi)^{n}c_{1}(X)^{n}$. Let 
\[
\mathcal{H}=\mathcal{H}(X,\omega_{r})=\{u\in C^{\infty}(X,\mathbb{R})\mid\omega_{u}\coloneqq\omega_{r}+i\partial\bar{\partial}u>0\}
\]
be the space of K\"{a}hler potentials. For each $u\in\mathcal{H}$,
there is a unique smooth function $h(u)$ (also denoted by $h_{\omega_{u}}$)
such that 
\[
Ric(\omega_{u})-\omega_{u}=i\partial\bar{\partial}h(u),\ \int_{X}e^{h(u)}\omega_{u}^{n}=V,
\]
where $Ric(\omega_{u})=-i\partial\bar{\partial}\log\omega_{u}^{n}$
is the Ricci curvature form. We call $h(u)$ the \textsl{Ricci potential}
of $\omega_{u}$, which is equal to zero if and only if $\omega_{u}$
is a KE metric. Let $h_{r}$ be the Ricci potential of $\omega_{r}$,
then we have 
\begin{equation}
e^{h(u)}=\frac{V}{\int_{X}e^{-u+h_{r}}\omega_{r}^{n}}\frac{e^{-u+h_{r}}\omega_{r}^{n}}{\omega_{u}^{n}},\ \textrm{for all}\ u\in\mathcal{H}.\label{eq: Ricci poten}
\end{equation}

\subsection{Extremal vector fields}

Let $\textrm{Aut}^{0}(X)$ be the identity component of the automorphism
group of $X$. Its Lie algebra is $\mathfrak{h}(X)$ consisting of
holomorphic vector fields. For any K\"{a}hler metric $\omega$ and
$v\in\mathfrak{h}(X)$, since the Fano manifolds are always simply
connected, there is a unique smooth function $\theta_{v}(\omega):X\rightarrow\mathbb{C}$
such that 
\[
v=\textrm{grad}_{\omega}^{(1,0)}\theta_{v}(\omega)\Leftrightarrow\iota_{v}\omega=i\bar{\partial}\theta_{v}(\omega),\ \textrm{and}\ \int_{X}\theta_{v}(\omega)\omega^{n}=0,
\]
where $\textrm{grad}_{\omega}^{(1,0)}f\coloneqq\sum_{i,j}g^{i\bar{j}}\frac{\partial f}{\partial\bar{z}^{j}}\frac{\partial}{\partial z^{i}}$.
We call $\theta_{v}(\omega)$ the Hamiltonian function of $v$ with
respect to $\omega$. 

We take a maximal compact subgroup $K$ of $\textrm{Aut}^{0}(X)$
with Lie algebra $\mathfrak{k}$. It is unique up to the conjugations.
Let $\omega\in2\pi c_{1}(X)$ be a $K$-invariant K\"{a}hler metric,
then for any $v\in\mathfrak{k}$, its Hamiltonian function $\theta_{v}(\omega)$
is real-valued. We endow the space $C^{\infty}(X,\mathbb{R})$ with
the inner product $\left(f,g\right)_{\omega}=\frac{1}{V}\int fg\omega^{n}$,
and let 
\[
P_{\omega}:C^{\infty}(X,\mathbb{R})\rightarrow\{\theta_{v}(\omega)\mid v\in\mathfrak{k}\}\subset C^{\infty}(X,\mathbb{R})
\]
be the orthogonal projection. 
\begin{defn}
For a maximal compact subgroup $K\subset\textrm{Aut}^{0}(X)$ and
a $K$-invariant K\"{a}hler metric $\omega\in2\pi c_{1}(X)$. The
\textit{extremal vector field} is defined by 
\[
Z_{K}=\textrm{grad}_{\omega}^{(1,0)}P_{\omega}\left(S(\omega)-n\right)\in\mathfrak{k}.
\]
\end{defn}

By \cite{Futaki Mabuchi}, $Z_{K}$ is independent of the choice of
$\omega$. Refer to \cite{Futaki Mabuchi} for the other properties
of $Z_{K}$, such as it belongs to the center of $\mathfrak{k}$ and
its imaginary part generates a $\mathbb{S}^{1}$-action. 
\begin{rem}
\label{rem: extrm vec}In \cite{Futaki Mabuchi}, they first introduce
a bilinear form on $\mathfrak{k}$, 
\[
B_{K}(v,w)\coloneqq\frac{1}{V}\int\theta_{v}(\omega)\theta_{w}(\omega)\omega^{n},\ \textrm{for}\ v,\ w\in\mathfrak{k}.
\]
It is independent of the choice of $\omega$. Then $Z_{K}$ is defined
to be the dual of the Futaki character w.r.t. $B_{K}$, namely it
satisfies 
\[
F(v)=-B_{K}(v,Z_{K}),\ \textrm{for all}\ v\in\mathfrak{k}.
\]
So $Z_{K}$ is independent of the choice of $\omega$, but it depends
on $K$. For another $K'=gKg^{-1}$, $g\in\textrm{Aut}^{0}(M)$, we
have $Z_{K'}=g_{*}Z_{K}$. 
\end{rem}

\subsection{Mabuchi solitons }

In \cite{Mabu soliton}, Mabuchi observed the following interesting
fact. 
\begin{thm}
\label{thm: of Mabuchi}\cite{Mabu soliton} Let $X$ be a Fano manifold,
and $K$ is a maximal compact subgroup of $\mathrm{Aut}^{0}(X)$.
For any $K$-invariant metric $\omega\in2\pi c_{1}(X)$, we have 
\[
P_{\omega}\left(S(\omega)-n\right)=P_{\omega}(1-e^{h_{\omega}}).
\]
Hence $Z_{K}=\mathrm{grad}_{\omega}^{(1,0)}P_{\omega}(1-e^{h_{\omega}})$. 
\end{thm}

Motivated by this observation, Mabuchi introduced the following analogue
of extremal metrics. 
\begin{defn}
\cite{Mabu soliton} For a Fano manifold $X$, a K\"{a}hler metric
$\omega\in2\pi c_{1}(X)$ with Ricci potential $h$ is called a \textsl{Mabuchi
soliton} if $\textrm{grad}_{\omega}^{(1,0)}(1-e^{h})$ is a holomorphic
vector field. 
\end{defn}

When the Futaki invariant for the class $c_{1}(X)$ vanishes, Mabuchi
solitons are KE metrics, hence they are the generalization of KE metrics.
We call them soliton since they give the soliton solutions for the
gradient flow of the Ding functional, see \cite{collins}. 

Mabuchi solitons share many common properties with KE metrics, such
as the uniqueness modulo the action of $\textrm{Aut}^{0}(X)$, and
their isometry groups are maximal compact subgroups of $\textrm{Aut}^{0}(X)$,
see \cite{Mabu 2 Vector field,Mabu 3 Calabi-Matsu,Mabu 4 Multiplier Herm}. 

The potential function of a Mabuchi soliton satisfies a Monge-Amp\`{e}re
type equation. To see this, suppose $\omega$ is a Mabuchi soliton
with isometry group $K$ which is maximal compact. We take a $K$-invariant
reference metric $\omega_{r}$ and express $\omega=\omega_{u}\coloneqq\omega_{r}+i\partial\bar{\partial}u$.
Let $v=\textrm{grad}_{\omega_{u}}^{(1,0)}(1-e^{h(u)})$, it is holomorphic
by definition. Since the Lie derivative $\mathcal{L}_{v}\omega_{u}=i\partial\bar{\partial}(1-e^{h(u)}),$
thus the imaginary part of $v$ preserving $\omega_{u}$, so $v\in\textrm{Lie}(K)$.
By Theorem \ref{thm: of Mabuchi}, $v$ is equal to $Z_{K}$. Then
we have $1-e^{h(u)}=\theta_{Z_{K}}(\omega_{u})=\theta_{Z_{K}}(\omega_{r})+Z_{K}(u)$.
By the formula (\ref{eq: Ricci poten}), after possibly adding a constant
to $u$, we obtain 
\begin{equation}
\left(1-\theta_{Z_{K}}(\omega_{r})-Z_{K}(u)\right)\omega_{u}^{n}=e^{-u+h_{r}}\omega_{r}^{n}.\label{eq: equation for MS}
\end{equation}

\subsection{Critical points of the Ding energy}

We give a variational characterization for Mabuchi solitons. First
we consider the norm of the Ding functional's gradient, which is the
analogue of the Calabi energy in the new GIT model \cite{Dona Ding}. 
\begin{defn}
For a Fano manifold $X$ with a reference metric $\omega_{r}$, the
\textsl{Ding energy} is defined by 
\[
Din(u)\coloneqq\frac{1}{V}\int_{X}\left(1-e^{h(u)}\right)^{2}\omega_{u}^{n},\ \mathrm{for}\ u\in\mathcal{H}(X,\omega_{r}).
\]
\end{defn}

The following observation is the starting point of this paper, it
is analogous to the relation of Calabi energies to extremal metrics.
The proof relies on Futaki's weighted Laplacian operator, see Section
2.4 \cite{Futaki book}. 
\begin{thm}
\label{thm: critical pt of DinE}The critical points of the Ding energy
are exactly the Mabuchi solitons. 
\end{thm}

\begin{proof}
In this proof, $h$ is short for $h(u)$. Suppose $u\in\mathcal{H}(X,\omega_{r})$
is a critical point. Consider a variation $f=\delta u$ of $u$, by
a direct computation with (\ref{eq: Ricci poten}), we have 
\[
\delta h(u)=-\triangle_{u}f-f+\frac{1}{V}\int fe^{h(u)}\omega_{u}^{n},
\]
where $\triangle_{u}=\nabla^{\bar{i}}\nabla_{\bar{i}}$ is the Laplacian
operator of $\omega_{u}$. Thus 
\begin{eqnarray*}
\delta\left(V\cdot Din(u)\right) & = & \int2e^{2h}\cdot\delta h(u)\cdot\omega_{u}^{n}+\int e^{2h}\triangle_{u}f\omega_{u}^{n}\\
 & = & -\int e^{2h}\triangle_{u}f\omega_{u}^{n}-2\int\left(e^{h}-\frac{1}{V}\int e^{2h}\omega_{u}^{n}\right)fe^{h}\omega_{u}^{n}.
\end{eqnarray*}
Let $\triangle_{h}\coloneqq\triangle_{u}+\nabla^{\bar{i}}h\cdot\nabla_{\bar{i}}$
be the Futaki's weighted Laplacian operator (see Page 40 of \cite{Futaki book},
we reverse the sign of $\triangle_{h}$). By a direct computation,
we have $\triangle_{u}e^{2h}=2e^{h}\triangle_{h}e^{h}$. Integrating
by parts, we see the above first term $\int e^{2h}\triangle_{u}f\omega_{u}^{n}=\int2\triangle_{h}e^{h}\cdot fe^{h}\omega_{u}^{n}$.
Hence 
\[
\delta\left(V\cdot Din(u)\right)=-2\int\left(\triangle_{h}e^{h}+e^{h}-\frac{1}{V}\int e^{2h}\omega_{u}^{n}\right)fe^{h}\omega_{u}^{n}.
\]
It vanishes for all smooth function $f$, so we have $\triangle_{h}e^{h}+e^{h}=\frac{1}{V}\int e^{2h}\omega_{u}^{n}$.
It implies that $G\coloneqq e^{h}-\frac{1}{V}\int e^{h}\cdot e^{h}\omega_{u}^{n}$
is an eigenfunction for operator $-\triangle_{h}$, i.e. $-\triangle_{h}G=G$.
By Theorem 2.4.3 \cite{Futaki book}, it follows that $\textrm{grad}_{\omega_{u}}^{(1,0)}e^{h}$
is holomorphic, so $\omega_{u}$ is a Mabuchi soliton. The converse
direction also follows from the above. 
\end{proof}

\section{Toric Fano varieties and the toric pluripotential theory }

Before we start to study the Mabuchi solitons on toric Fano varieties,
in this section we review the invariant pluripotential theory on these
spaces. The references are \cite{Berman Berndtsson,Guedj toric potential}. 

\subsection{\label{subsec:Toric-Fano-Varieties}Toric Fano varieties}

We briefly recall the construction of toric Fano varieties, refer
to \cite{CLS toric,Fulton toric} for an introduction to toric varieties,
and \cite{Debarre} for toric Fano varieties. 

A normal projective variety $X$ (over $\mathbb{C}$) is called a
\textit{Fano variety} if its anti-canonical divisor $-K_{X}$ is $\mathbb{Q}$-Cartier
and ample, moreover $X$ has log terminal singularities. Where $\mathbb{Q}$-Cartier
means that $-rK_{X}$ is a Cartier divisor for some $r\in\mathbb{N}$,
the smallest such $r$ is called the Gorenstein index of $X$. When
$-K_{X}$ is Cartier, we say $X$ is \textit{Gorenstein}. See Section
\ref{subsec: def of klt} for the meaning of log terminal singularities. 

A $n$-dimensional \textit{toric variety} is an algebraic variety
with an effective action by the complex torus $\left(\mathbb{C}^{*}\right)^{n}$
and it has an open dense orbit. All the normal toric varieties can
be constructed from a kind of combinatorial data called \textit{fan}.
Let $N$ be a lattice of rank $n$, $M=\textrm{Hom}_{\mathbb{Z}}(N,\mathbb{Z})$
is the dual lattice. A fan $\Sigma=\{\sigma\}$ in $N_{\mathbb{R}}\coloneqq N\otimes_{\mathbb{Z}}\mathbb{R}$
is a collection of cones such that each cone $\sigma$ is generated
by finite many elements in $N$ and $\sigma\cap\tau\in\Sigma$ for
any $\sigma,\tau\in\Sigma$. In the following, we only consider the
fans given by Fano polytopes. 
\begin{defn}
A polytope $Q\subset N_{\mathbb{R}}$ is called a \textsl{Fano polytope}
if $0\in\mathrm{Int}(Q)$ and all the vertices are primitive elements
in $N$. The collection of cones spanned by the faces of $Q$ (except
the whole $Q$) and plus $\{0\}$ constitute a fan, denoted by $\Sigma_{Q}$. 
\end{defn}

A Fano polytope $Q$ gives rise to a toric Fano variety in the following
way. For each cone $\sigma\in\Sigma_{Q}$, let $\sigma^{\lor}\coloneqq\{y\in M_{\mathbb{R}}\mid y|_{\sigma}\geq0\}$
be the dual cone. Then $\sigma^{\lor}\cap M$ is a finitely generated
semigroup, it induces a finitely generated $\mathbb{C}$-algebra 
\[
\mathbb{C}[\sigma^{\lor}\cap M]\coloneqq\{\sum a_{u}\chi^{u}\mid u\in\sigma^{\lor}\cap M,\ a_{u}\in\mathbb{C}\},
\]
where the sum is finite. Let $U_{\sigma}=\mathrm{Spec}\mathbb{C}[\sigma^{\lor}\cap M]$
be the associated affine variety which has an action by the torus
$T_{N}\coloneqq N\otimes_{\mathbb{Z}}\mathbb{C}^{*}$. For example,
the smallest cone $\{0\}$ gives $U_{\{0\}}=T_{N}$. According to
the inclusion relation between cones, these affine varieties are glued
together to a normal variety $X$, which has a $T_{N}$-action and
an open dense orbit $U_{\{0\}}$.

When we study the K\"{a}hler geometry of $X$, it is more natural
to work with the dual polytope. Let $\{q_{\alpha}\in N\mid\alpha\in\Lambda\}$
be the vertices of $Q$, the dual polytope is defined by
\[
P=Q^{\mathrm{dual}}\coloneqq\{y\in M_{\mathbb{R}}\mid l_{\alpha}(y)\coloneqq\left\langle q_{\alpha},y\right\rangle +1\geq0,\ \textrm{for all}\ \alpha\in\Lambda\}.
\]
We also have $Q=P^{\mathrm{dual}}$. Note that $P$ is not necessarily
a lattice polytope, if it does we call $P$ and $Q$  \textit{reflexive}.
There is an orbit-face correspondence, namely each $k$-dimensional
torus orbit in $X$ is one-to-one corresponding to a $k$-dimensional
face of $P$. For example, the fixed points on $X$ are corresponding
to the vertices of $P$; the invariant prime divisors are corresponding
to the facets (codimension 1) of $P$, then to the vertices of $Q$. 

Let $\{D_{\alpha}\mid\alpha\in\Lambda\}$ be the invariant prime divisors,
then the anti-canonical divisor is given by 
\[
-K_{X}=\sum_{\alpha\in\Lambda}D_{\alpha}.
\]
For $r\in\mathbb{N}$, $-rK_{X}$ is Cartier if and only if $rP$
is a lattice polytope, hence the smallest such $r$ is equal to the
Gorenstein index of $X$. Thus $X$ is Gorenstein ($K_{X}$ is Cartier)
if and only if $P$ is a lattice polytope, i.e. $P$, $Q$ are reflexive.
By Proposition 12 \cite{Debarre}, $X$ has log terminal singularities,
and by Proposition 6.1.10 \cite{CLS toric}, $-K_{X}$ is ample. Hence
$X$ is a toric Fano variety. Moreover, all the toric Fano varieties
can be obtained in this way. They are our research objects in this
paper, so we make a summary. 
\begin{defn}
Let $N$ be a lattice of rank $n$ with the dual lattice $M$. Let
$Q\subset N_{\mathbb{R}}$ be a Fano polytope with the dual polytope
$P\subset M_{\mathbb{R}}$. We denote by $X_{P}$ the toric Fano variety
associated to $Q$ in the above way. Sometimes, we omit the subscript
when $P$ is clear. 
\end{defn}

We say $P$ is a \textit{Delzant polytope} if (1) $P$ is a lattice
polytope; (2) each vertex of $P$ adjoins to $n$ edges, and the primitive
generators of these $n$ edges form a basis of $M$. By Theorem 2.4.3
\cite{CLS toric}, $X_{P}$ is smooth if and only if $P$ is a Delzant
polytope. 

Classifying toric Fano varieties up to isomorphisms is equivalent
to classifying Fano polytopes up to the isomorphisms of the lattice
$N$. In each dimension, there are finite many isomorphism classes
of toric Fano varieties with discrepancy bounded below away from $-1$.
For Gorenstein toric Fano varieties (i.e. reflexive polytopes), the
numbers of isomorphism classes in each dimension are 16($n=2$), 4319,
473800776, ... For smooth toric Fano varieties, the numbers are 5($n=2$),
18, 124, 866, 7622, ... See \cite{Fano polytopes} for a survey of
classifications. 

\subsection{Invariant psh metrics on $-K_{X}$}

Let $X_{P}$ be the toric Fano variety. There is a $r\in\mathbb{N}$
such that $-rK_{X}$ is a line bundle. The $T_{N}$-action has a canonical
lifting to $-rK_{X}$. Let $S_{N}\coloneqq N\otimes_{\mathbb{Z}}\mathbb{S}^{1}\subset T_{N}$
be the compact torus, then we consider the $S_{N}$-invariant metrics
on $-rK_{X}$. In the following discussions, for simplicity, we assume
$r=1$, i.e. $X$ is Gorenstein. 

By choosing a basis $\{e_{i}\}$ of $N$, the lattice $N$ is isomorphic
to $\mathbb{Z}^{n}$, then $N_{\mathbb{R}}\cong\mathbb{R}^{n}$ and
$M_{\mathbb{R}}\cong\mathbb{R}^{n}$. The torus $T_{N}$ is isomorphic
to $\left(\mathbb{C}^{*}\right)^{n}=\{(z_{i})_{i=1}^{n}\mid z_{i}\neq0\}$,
and $S_{N}$ is isomorphic to $\left(\mathbb{S}^{1}\right)^{n}$.
It is convenient to use the logarithmic coordinates on $T_{N}$, let
\[
w_{i}=\log z_{i}^{2}=x_{i}+\sqrt{-1}\theta_{i},\ x_{i}=\log\left|z_{i}\right|^{2},\ 0\leq\theta_{i}<4\pi,
\]
or $z_{i}=e^{\frac{1}{2}w_{i}}=e^{\frac{1}{2}(x_{i}+\sqrt{-1}\theta_{i})}$.
We define a meromorphic section of $-K_{X}$, 
\begin{equation}
s\coloneqq\frac{\partial}{\partial w_{1}}\wedge\cdots\wedge\frac{\partial}{\partial w_{n}}=\frac{1}{2^{n}}z_{1}\cdots z_{n}\frac{\partial}{\partial z_{1}}\wedge\cdots\wedge\frac{\partial}{\partial z_{n}},\label{eq: local frame}
\end{equation}
which is independent of the choice of the basis $\{e_{i}\}$ up to
signs. The associated divisor of $s$ is the anti-canonical divisor
$\sum_{\alpha\in\Lambda}D_{\alpha}$. 

Let $\textrm{PSH}(-K_{X})$ be the space of psh metrics on $-K_{X}$,
refer to \cite{BBEGZ,BEGZ} for a discussion of psh metrics on singular
varieties. When $X$ is singular, we take a toric resolution $\pi:\tilde{X}\rightarrow X$.
Since $X$ is normal, $\textrm{PSH}(-K_{X})$ can be identified with
$\textrm{PSH}(\tilde{X},-\pi^{*}K_{X})$. 

Let $\textrm{PSH}(-K_{X})^{S}$ be the subspace of $S_{N}$-invariant
psh metrics. Given a metric $\phi\in\textrm{PSH}(-K_{X})^{S}$, we
have a function 
\[
-\log\left|s\right|_{\phi}^{2}:T_{N}\rightarrow\mathbb{R}\cup\{-\infty\},
\]
where $s$ is defined by (\ref{eq: local frame}). Abusing the notations
slightly, we denote the above function also by $\phi$. In the coordinates
$(w_{i})$ on $T_{N}$, by the $S_{N}$-invariance, $\phi(w)$ only
depends on $(x_{i})\in\mathbb{R}^{n}$. Hence in the follows, we take
$\phi$ as a function defined on $\mathbb{R}^{n}$, which is convex
since $i\partial\bar{\partial}\phi$ is a positive current. In particular,
this implies the psh metric $\phi$ must be continuous over $T_{N}$. 

Following \cite{Berman Berndtsson}, we introduce two spaces of convex
functions, where they denote them by $\mathcal{P}(\mathbb{R}^{n})$
and $\mathcal{P}_{+}(\mathbb{R}^{n})$. 
\begin{defn}
Let $P\subset M_{\mathbb{R}}$ be a convex body containing $0$ in
its interior. Let $\psi_{P}(x)=\sup_{y\in P}\left\langle x,y\right\rangle $
be the support function of $P$. We define 
\[
\textrm{PSH}(P)\coloneqq\{\phi:\mathbb{R}^{n}\rightarrow\mathbb{R}\ \textrm{convex}\mid\phi\leq\psi_{P}+C\ \textrm{on}\ \mathbb{R}^{n}\ \textrm{for\ some}\ C\},
\]
\[
\textrm{PSH}_{b}(P)\coloneqq\{\phi\in\textrm{PSH}(P)\mid\psi_{P}-C\leq\phi\ \textrm{on}\ \mathbb{R}^{n}\ \textrm{for\ some}\ C\}.
\]
\end{defn}

By \cite{Guedj toric potential} Proposition 3.2, $\textrm{PSH}(-K_{X_{P}})^{S}$
can be identified with $\textrm{PSH}(P)$ by sending a psh metric
$\phi$ on $-K_{X_{P}}$ to the associated convex function on $\mathbb{R}^{n}$
(which is also denoted by $\phi$). The subspace of bounded metrics
$\textrm{PSH}_{b}(-K_{X_{P}})^{S}$ is identified with $\textrm{PSH}_{b}(P)$.
We will further discuss these identifications in Section \ref{sec: Dictionary}. 

The most important metric is the Fubini-Study metric. As before, we
assume $X$ is Gorenstein, otherwise we take $r$ to be sufficiently
divisible. For each $r\geq1$, we have a weight decomposition with
respect to the $T_{N}$-action, 
\[
\textrm{H}^{0}(-rK_{X})=\bigoplus_{p\in rP\cap M}\textrm{span}_{\mathbb{C}}\{z^{p}\cdot s^{\otimes r}\},\ z^{p}\coloneqq z_{1}^{p_{1}}\cdots z_{n}^{p_{n}}.
\]
These eigen-sections induce an invariant psh metric on $-K_{X}$ with
potential 
\[
\phi_{FS,r}=\frac{1}{r}\log\sum_{p\in rP\cap M}\left|z^{p}\right|^{2}=\frac{1}{r}\log\sum_{p\in rP\cap M}e^{\left\langle x,p\right\rangle }.
\]
We see $\phi_{FS,r}$ is smooth and strictly convex on $\mathbb{R}^{n}$.
Moreover, we have 
\[
\psi_{P}\leq\phi_{FS,r}\leq\psi_{P}+\frac{1}{r}\log\left|rP\cap M\right|,
\]
so $\phi_{FS,r}\in\textrm{PSH}_{b}(P)$. On the other hand, each single
eigen-section induces a psh metric, the supremum of them is still
a psh metric, which is equal to 
\[
\sup_{p\in rP\cap M}\frac{1}{r}\log\left|z^{p}\right|^{2}=\sup_{p\in P}\left\langle x,p\right\rangle =\psi_{P}.
\]
Although $\psi_{P}$ is not smooth nor strictly convex, but its Legendre
dual is very simple. In the toric setting, it is convenient to take
$\psi_{P}$ to be the reference metric. 

\subsection{Facts in convex analysis }

We review some relevant facts in convex analysis, the references are
\cite{Figalli,Berman Berndtsson}. 

\subsubsection{Monge-Amp\`{e}re measures}

For a convex function $\phi:\mathbb{R}^{n}\rightarrow\mathbb{R}$
and $x_{0}\in\mathbb{R}^{n}$, the set of \textsl{subgradients} of
$\phi$ at $x_{0}$ is defined by 
\[
\partial\phi(x_{0})\coloneqq\{y\in\mathbb{R}^{n}\mid\phi(x)\geq\phi(x_{0})+\left\langle x-x_{0},y\right\rangle ,\ \textrm{for all}\ x\in\mathbb{R}^{n}\}.
\]
It is a nonempty convex compact set. $\phi$ is differentiable at
$x_{0}$ iff $\partial\phi(x_{0})$ contains only one element, that
is the gradient of $\phi$ at $x_{0}$ denoted by $d\phi(x_{0})$.
If $K\subset\mathbb{R}^{n}$ is a compact set, then $\partial\phi(K)\coloneqq\cup_{x\in K}\partial\phi(x)$
is also compact. 

The \textsl{Monge-Amp\`{e}re measure} $\textrm{MA}(\phi)$ (in the
sense of Alexandrov) is defined by 
\[
\textrm{MA}(\phi)(E)\coloneqq\int_{\partial\phi(E)}dy,
\]
where $E\subset\mathbb{R}^{n}$ is a Borel set and $dy$ denotes the
Lebesgue measure. It is a nonnegative and locally finite Borel measure.
When $\phi$ is smooth and strictly convex, we have $\textrm{MA}(\phi)=\det[\phi_{ij}]dx$.
For general convex functions, $\textrm{MA}(\phi)$ may not be absolutely
continuous with respect to the Lebesgue measure, e.g. $\textrm{MA}(\left|x\right|)=\left|B_{1}(0)\right|\delta_{0}$,
where $\delta_{0}$ is the Dirac measure at the origin. 

We need a modified version of the Monge-Amp\`{e}re measure. For a
convex function $\phi$ on $\mathbb{R}^{n}$, let $P=\overline{\partial\phi(\mathbb{R}^{n})}$.
For a continuous function $0\leq g(y)\in L^{1}(P)$, we define the
$g$-Monge-Amp\`{e}re measure to be
\[
\textrm{MA}_{g}(\phi)(E)\coloneqq\int_{\partial\phi(E)}gdy.
\]
When $\phi$ is smooth and strictly convex, it equals to $g(d\phi)\det[\phi_{ij}]dx$.
When we study the Mabuchi solitons on toric varieties, $P$ is a convex
body, and $g=\ell$ is an affine function such that $P\subset\{\ell\geq0\}$.
In this case, we have $\textrm{MA}_{\ell}(\phi)\leq C\cdot\textrm{MA}(\phi)$
for some $C>0$. In the converse direction, since $\{\ell=0\}$ is
a null set for $dy$, $\textrm{MA}(\phi)$ is also absolutely continuous
with respect to $\textrm{MA}_{\ell}(\phi)$. Both of them are finite
Borel measures, so by the Radon-Nikodym theorem, there is a Borel
measurable function denoted by $\ell(\partial\phi)$ (defined almost
everywhere w.r.t. $\textrm{MA}(\phi)$) such that 
\[
\int_{\partial\phi(E)}\ell dy=\textrm{MA}_{\ell}(\phi)(E)=\int_{E}\ell(\partial\phi)\textrm{MA}(\phi)
\]
for all Borel set $E$. At the points where $\phi$ is differentiable,
we have $\ell(\partial\phi)=\ell(d\phi)$. Moreover, at any point
$x\in\mathbb{R}^{n}$ such that the measure $\left|\partial\phi(x)\right|_{dy}>0$,
by taking $E=\{x\}$, we have 
\[
\ell(\partial\phi)(x)=\frac{\int_{\partial\phi(x)}\ell dy}{\int_{\partial\phi(x)}dy}.
\]

\subsubsection{Legendre duals }

For a convex function $\phi:N_{\mathbb{R}}\rightarrow\mathbb{R}$,
the \textsl{Legendre dual} (or transform) is defined by 
\[
\phi^{*}(y)=\sup_{x\in\mathbb{R}^{n}}\left\langle x,y\right\rangle -\phi(x):M_{\mathbb{R}}\rightarrow\mathbb{R}\cup\{\infty\}
\]
As the supremum of a family of affine functions, $\phi^{*}$ is convex
and lower semi-continuous (lsc). For the support function $\psi_{P}$
of a convex body $P$, we see $\psi_{P}^{*}=0$ on $P$ and $\infty$
elsewhere. For $x\in N_{\mathbb{R}}$ and $y\in M_{\mathbb{R}}$,
we have 
\[
\phi^{*}(y)=\left\langle x,y\right\rangle -\phi(x)\Leftrightarrow y\in\partial\phi(x)\Leftrightarrow x\in\partial\phi^{*}(y).
\]
Thus $\partial\phi(\mathbb{R}^{n})\subset\{\phi^{*}<\infty\}$, the
inclusion may be strict. The Legendre dual of $\phi^{*}$ is $\phi$.
For two convex functions $\phi$ and $\psi$, $\phi\leq\psi$ on $N_{\mathbb{R}}$
if and only if $\phi^{*}\geq\psi^{*}$ on $M_{\mathbb{R}}$. 

Let $P$ be a convex body containing $0$ in its interior. Then $\phi\in\textrm{PSH}(P)$
if and only if $\{\phi^{*}<\infty\}\subset P$, which is equivalent
to $\partial\phi(\mathbb{R}^{n})\subset P$. For $\phi\in\textrm{PSH}(P)$,
we have 
\begin{equation}
\sup_{\mathbb{R}^{n}}\left(\phi-\psi_{P}\right)=-\inf_{P}\phi^{*}=\phi(0),\ \inf_{\mathbb{R}^{n}}\left(\phi-\psi_{P}\right)=-\sup_{P}\phi^{*},\label{eq: sup-inf}
\end{equation}
hence $\left\Vert \phi-\psi_{P}\right\Vert _{\infty,\mathbb{R}^{n}}=\left\Vert \phi^{*}\right\Vert _{\infty,P}$.
It follows that $\phi\in\textrm{PSH}_{b}(P)$ iff $\phi^{*}$ is bounded
above on $P$ ($\phi^{*}$ is always bounded below). 

We say $\phi\in\textrm{PSH}(P)$ is \textit{normalized} (at $0$)
if $\phi\geq\phi(0)=0$, this is equivalent to that $\phi^{*}$ is
normalized at $0\in P$. For any $p\in\partial\phi(0)$, $\phi-\left\langle x,p\right\rangle -\phi(0)$
is normalized. A \textsl{translation} of $\phi$ is defined to be
$\phi_{a}(x)\coloneqq\phi(x+a)$, where $a\in N_{\mathbb{R}}$. A
\textsl{rescaling} of $\phi$ is defined to be $\phi^{(r)}\coloneqq r\phi(\frac{x}{r})$,
where $r>0$. We have $\phi_{a}^{*}=\phi^{*}-\left\langle a,\cdot\right\rangle $
and $\left(\phi^{(r)}\right)^{*}=r\phi^{*}$. 

\subsubsection{Jensen's inequality }

We will apply Jensen's inequality frequently. Let $\mu$ be a probability
Borel measure on $\mathbb{R}^{n}$ with compact support. Let $\mathfrak{c}\in\mathbb{R}^{n}$
be its barycenter. The Jensen inequality says that for any convex
function $g$ defined on the convex hull of $\textrm{supp}\mu$, we
have $\int gd\mu\geq g(\mathfrak{c})$, and the equality holds if
and only if $g$ is affine over the convex hull of $\textrm{supp}\mu$. 

\subsection{\label{sec: Dictionary}The dictionary for psh metrics v.s. convex
functions }

We recall some important classes of psh metrics and the corresponding
classes of convex functions on $\mathbb{R}^{n}$. Refer to \cite{BBEGZ,BEGZ}
for a discussion of psh metrics, and \cite{Berman Berndtsson,Guedj toric potential}
for the psh-convex correspondences. In this section, $X$ is a Gorenstein
toric Fano variety associated to lattice polytope $P$. 

\subsubsection{Plurisubharmonic metrics}

For any $\phi\in\textrm{PSH}(-K_{X})^{S}$, since $\phi-\psi_{P}$
is bounded above on $X$, thus $\phi\in\textrm{PSH}(P)$. Conversely,
given a $\phi\in\textrm{PSH}(P)$, it gives a $S_{N}$-invariant continuous
metric on $-K_{X}|_{T_{N}}$. Since $X\backslash T_{N}$ is an analytic
set, it can be extended to an invariant psh metric on $-K_{X}$. Hence
$\textrm{PSH}(-K_{X})^{S}$ can be identified with $\textrm{PSH}(P)$.
In the following, we always make this identification without mentions. 

\subsubsection{Metrics with full Monge-Amp\`{e}re mass}

For $\phi\in\textrm{PSH}(-K_{X})$, its curvature current is $i\partial\bar{\partial}\phi$.
We can form a non-pluripolar product of $i\partial\bar{\partial}\phi$,
denoted by $\textrm{MA}_{\mathbb{C}}(\phi)=\left(i\partial\bar{\partial}\phi\right)^{n}$,
see \cite{BEGZ} for the definition. It has no mass on analytic subsets.
If $\phi$ is $S_{N}$-invariant, then it is continuous over $T_{N}$,
the non-pluripolar product can be firstly defined on $T_{N}$ via
Bedford-Taylor's theory and then take the zero extension. The total
mass satisfies $\int_{X}\left(i\partial\bar{\partial}\phi\right)^{n}\leq V\coloneqq(2\pi)^{n}c_{1}(X)^{n}$.
When the equality holds, $\phi$ is said to be of full (Monge-Amp\`{e}re)
mass. Bounded metrics always have full mass. Let $\mathcal{E}(-K_{X})$
be the space of psh metrics with full mass, and $\mathcal{E}(-K_{X})^{S}$
be the subspace of the $S_{N}$-invariant such metrics. 

We define the tropicalization map to be
\[
\mathscr{L}:T_{N}\rightarrow\mathbb{R}^{n},\ \mathscr{L}(z_{1},\cdots,z_{n})=\left(\log\left|z_{1}\right|^{2},\cdots,\log\left|z_{1}\right|^{2}\right).
\]
By Lemma 2.2 \cite{Guedj toric potential}, we have 
\begin{equation}
\mathscr{L}_{*}\ \textrm{MA}_{\mathbb{C}}(\phi)|_{T_{N}}=n!(2\pi)^{n}\textrm{MA}(\phi),\label{eq: MA-MAR}
\end{equation}
where $\mathscr{L}_{*}$ is the push-forward for measures. In particular,
if $\phi$ is smooth, take the total mass, we obtain $c_{1}(X)^{n}=n!V_{P}$,
where $V_{P}$ is the Lebesgue volume of $P$. We follow \cite{Berman Berndtsson}
(they use notation $\mathcal{E}_{P}(\mathbb{R}^{n})$) to define 
\[
\mathcal{E}(P)\coloneqq\{\phi\in\textrm{PSH}(P)\mid\int_{\mathbb{R}^{n}}\textrm{MA}(\phi)=V_{P}\}.
\]
Then it can be identified with $\mathcal{E}(-K_{X})^{S}$. In fact,
for any $\phi\in\mathcal{E}(P)$, (\ref{eq: MA-MAR}) implies the
corresponding psh metric has full mass. Conversely, given $\phi\in\mathcal{E}(-K_{X})^{S}$,
since the non-pluripolar product $\textrm{MA}_{\mathbb{C}}(\phi)$
has no mass on $X\backslash T_{N}$, thus by (\ref{eq: MA-MAR}) we
have $\phi\in\mathcal{E}(P)$. 

Moreover, by Theorem 3.6 \cite{Guedj toric potential}, for $\phi\in\textrm{PSH}(P)$,
it belongs to $\mathcal{E}(P)$ iff $\phi^{*}<\infty$ on $\mathrm{Int}(P)$.
In particular, we have $\textrm{PSH}_{b}(P)\subset\mathcal{E}(P)$,
this inclusion is strict, e.g. take $P=[-1,1]$ and $\phi(x)=\left|x\right|-2\sqrt{\left|x\right|+1}$,
we see $\phi\in\mathcal{E}(P)$ but not $\textrm{PSH}_{b}(P)$. In
addition, Lemma 2.7 \cite{Guedj toric potential} says that for any
$\phi\in\mathcal{E}(P)$ and $\epsilon\in(0,1)$, there exists $C_{\epsilon}>0$
such that 
\begin{equation}
\phi\geq(1-\epsilon)\psi_{P}-C_{\epsilon},\ \mathrm{on}\ \mathbb{R}^{n}.\label{eq: proper of full MA}
\end{equation}
It follows that every $\phi\in\mathcal{E}(P)$ is proper on $\mathbb{R}^{n}$. 

As a metric on $-K_{X}$, each $\phi\in\textrm{PSH}(-K_{X})^{S}$
induces a measure 
\[
\Omega_{\phi}\coloneqq i^{n^{2}}\left|\beta\right|_{\phi}^{2}\beta^{*}\wedge\overline{\beta^{*}},
\]
where $\beta$ is a local frame of $-K_{X}$ with dual $\beta^{*}$.
By \cite{BBEGZ} Proposition 1.4, we have $\int_{X}\Omega_{\phi}<\infty$
if $\phi$ has full mass. If we take $\beta$ to be (\ref{eq: local frame}),
then $\Omega_{\phi}|_{T_{N}}=i^{n^{2}}e^{-\phi}dw\wedge d\bar{w}$.
Its push-forward by $\mathscr{L}$ is equal to
\begin{equation}
\mathscr{L}_{*}\ \Omega_{\phi}|_{T_{N}}=(8\pi)^{n}e^{-\phi}dx.\label{eq: mu-exp}
\end{equation}

\subsubsection{Metrics with finite energy}

The Monge-Amp\`{e}re energy 
\[
E:\textrm{PSH}(-K_{X})\rightarrow\mathbb{R}\cup\{-\infty\}
\]
is a functional (up to a constant) such that $\frac{d}{dt}E(\phi_{t})=\frac{1}{V}\int_{X}\dot{\phi}_{t}\left(i\partial\bar{\partial}\phi_{t}\right)^{n}$
for any path $\{\phi_{t}\}$ of smooth metrics. It is non-decreasing
and satisfies $E(\phi+c)=E(\phi)+c$ for any constant $c$. The space
of psh metrics with finite energy is defined to be 
\[
\mathcal{E}^{1}(-K_{X})=\{\phi\in\textrm{PSH}(-K_{X})\mid E(\phi)>-\infty\}.
\]
It is contained in $\mathcal{E}(-K_{X})$. Let $\mathcal{E}^{1}(-K_{X})^{S}$
be the subspace of the $S_{N}$-invariant such metrics. 

In the toric setting, we consider the restriction of $E$ on $\textrm{PSH}(-K_{X})^{S}\cong\textrm{PSH}(P)$.
If we normalize it by requiring $E(\psi_{P})=0$, then by Proposition
2.9 \cite{Berman Berndtsson}, we have 
\begin{equation}
E(\phi)=-\frac{1}{V_{P}}\int_{P}\phi^{*}dy,\ \phi\in\textrm{PSH}(P).\label{eq: MA energy}
\end{equation}
It follows that $\mathcal{E}^{1}(-K_{X})^{S}$ can be identified with
\[
\mathcal{E}^{1}(P)\coloneqq\{\phi\in\textrm{PSH}(P)\mid\int_{P}\phi^{*}dy<\infty\}.
\]
We see $\mathcal{E}^{1}(P)\subset\mathcal{E}(P)$. In fact, if $\phi\in\mathcal{E}^{1}(P)$,
then $\phi^{*}<\infty$ a.e. on $P$. Since $\phi^{*}$ is convex,
it implies $\phi^{*}<\infty$ on $\mathrm{Int}(P)$, hence $\phi\in\mathcal{E}(P)$. 

We need a modified version of $E$. Let $g\geq0$ be a continuous
function on $P$ with $\int_{P}gdy=1$, then we define 
\begin{equation}
E_{g}(\phi)\coloneqq-\int_{P}\phi^{*}\cdot gdy,\ \textrm{for}\ \phi\in\mathcal{E}^{1}(P).\label{eq: g-MA energy}
\end{equation}
It is exactly the restriction of the $g$-modified Monge-Amp\`{e}re
energy defined in \cite{Cpx optimal trans} Lemma 2.14. 

\subsubsection{Bounded psh metrics}

By (\ref{eq: sup-inf}), the space of the invariant bounded psh metrics
$\textrm{PSH}_{b}(-K_{X})^{S}$ can be identified with $\textrm{PSH}_{b}(P)$.
Proposition 4.1 \cite{Guedj toric potential} says that the $S_{N}$-invariant
bounded psh metrics are automatically continuous on $X$. 

\subsubsection{Smooth metrics with positive curvature }

We assume $X$ is smooth, equivalently $P$ is a Delzant polytope.
Let $\mathcal{H}(-K_{X})^{S}$ be the space of the invariant smooth
metrics on $-K_{X}$ with positive curvature. Let $G(y)=\sum_{\alpha\in\Lambda}l_{\alpha}\log l_{\alpha}$
be the Guillemin model symplectic potential on $P$. Let $\mathcal{H}(P)\subset\textrm{PSH}_{b}(P)$
be the subspace consisting of smooth and strictly convex $\phi$ such
that $\phi^{*}-G$ can be smoothly extended to a neighborhood of $P$.
Guillemin shows that $\mathcal{H}(-K_{X})^{S}$ can be identified
with $\mathcal{H}(P)$. For any $\phi\in\mathcal{H}(P)$, $d\phi:N_{\mathbb{R}}\rightarrow\mathrm{Int}(P)$
is a diffeomorphism. Moreover, $d\phi\circ\mathscr{L}$ can be extended
to a smooth map $\mathfrak{m}_{\phi}:X\rightarrow P$ called the \textit{moment
map} induced by $\phi$. The push-forward of $\textrm{MA}_{\mathbb{C}}(\phi)$
by $\mathfrak{m}_{\phi}$ is $n!(2\pi)^{n}dy$. 
\begin{rem}
In a summary, the following spaces of psh metrics 
\[
\mathcal{H}(-K_{X})^{S}\subset\textrm{PSH}_{b}(-K_{X})^{S}\subset\mathcal{E}^{1}(-K_{X})^{S}\subset\mathcal{E}(-K_{X})^{S}\subset\textrm{PSH}(-K_{X})^{S}
\]
can be identified with $\mathcal{H}(P)\subset\textrm{PSH}_{b}(P)\subset\mathcal{E}^{1}(P)\subset\mathcal{E}(P)\subset\textrm{PSH}(P)$
respectively. Note that except $\mathcal{H}(P)$, all these spaces
are preserved by rescaling $\phi(x)\mapsto r\phi(\frac{x}{r})$. 
\end{rem}

\subsection{Toric geodesics }

In the space of smooth K\"{a}hler metrics, smooth geodesic segments
with respect to Mabuchi's $L^{2}$-metric does not necessarily exist.
Hence some weak geodesics are introduced, such as the finite energy
geodesics by Darvas \cite{Darvas}. In toric setting, we consider
the geodesics of invariant metrics in $\mathcal{E}^{1}(P)$, called
the \textsl{toric geodesics}. They are easy to be described by the
Legendre transformation, see \cite{Ross Nystrom,Song Zelditch} for
more discussions. 

For any $\phi_{0},\phi_{1}\in\mathcal{E}^{1}(P)$, the toric geodesic
segment connecting them is given by 
\[
\phi_{t}=\left((1-t)\phi_{0}^{*}+t\phi_{1}^{*}\right)^{*}\in\mathcal{E}^{1}(P),\ t\in[0,1].
\]
The toric geodesic rays started from $\phi_{0}\in\mathcal{E}^{1}(P)$
are given by 
\begin{equation}
\phi_{t}\coloneqq\left(\phi_{0}^{*}+tf\right)^{*}\in\mathcal{E}^{1}(P),\ t\geq0,\label{eq: toric geo ray}
\end{equation}
where $f:\mathrm{Int}(P)\rightarrow\mathbb{R}$ is an integrable convex
function. It can be verified that $\phi_{t}(x)$ is convex in $(x,t)$.
By (\ref{eq: g-MA energy}), $E_{g}$ is affine along toric geodesic
rays. 

\subsection{Toric Ding functionals and Ricci potentials }

For $\phi\in\mathcal{E}^{1}(-K_{X})$, since the measure $\Omega_{\phi}$
has finite total mass, we define $L(\phi)\coloneqq-\log\int_{X}\Omega_{\phi}$.
Recall the Ding functional is defined by 
\[
D(\phi)=-E(\phi)+L(\phi):\mathcal{E}^{1}(-K_{X})\rightarrow\mathbb{R}.
\]
Consider its restriction on $\mathcal{E}^{1}(-K_{X})^{S}\cong\mathcal{E}^{1}(P)$,
by (\ref{eq: MA energy}) and (\ref{eq: mu-exp}), we have 
\begin{equation}
D(\phi)=\frac{1}{V_{P}}\int_{P}\phi^{*}dy-\log\int_{\mathbb{R}^{n}}e^{-\phi}dx:\mathcal{E}^{1}(P)\rightarrow\mathbb{R}.\label{eq: toric Ding}
\end{equation}
For a toric geodesic ray $\{\phi_{t}\}$ given by (\ref{eq: toric geo ray}),
the Pr\'{e}kopa-Leindler inequality implies that $-\log\int_{\mathbb{R}^{n}}e^{-\phi_{t}}$
is convex in $t$. Since $E(\phi_{t})$ is affine, we know $D(\phi_{t})$
is convex, this is a special case of Berndtsson's convexity \cite{Berndtsson}. 

Assume $X$ is smooth. For $\phi\in\mathcal{H}(-K_{X})^{S}$, denote
its Ricci potential by $h(\phi)$, then we have 
\begin{equation}
e^{h(\phi)}=\frac{V}{\int_{X}\Omega_{\phi}}\cdot\frac{\Omega_{\phi}}{\textrm{MA}_{\mathbb{C}}(\phi)}=\left(\frac{V_{P}}{\int_{\mathbb{R}^{n}}e^{-\phi}}\cdot\frac{e^{-\phi}}{\det[\phi_{ij}]}\right)\circ\mathscr{L},\ \textrm{on}\ T_{N}.\label{eq: toric Ricci poten}
\end{equation}

\section{Limit slopes of the toric Ding functional }

In this section, we derive an explicit formula for the limit slopes
of the Ding functional along toric geodesic rays. When the geodesic
rays induced from toric test-configurations, we recompute it by Berman's
lct formula \cite{Berman}. Two results coincide with each other. 

\subsection{Limit slopes along toric geodesic rays }

Let $\{\phi_{t}\}$ be a toric geodesic ray given by (\ref{eq: toric geo ray}).
The Monge-Amp\`{e}re energy $E(\phi_{t})$ is affine in $t$ with
the slope $-\frac{1}{V_{P}}\int_{P}f$. For the second term of the
Ding functional, $L(\phi_{t})$ is convex in $t$. Taking the derivative,
it leads us to consider the asymptotic behavior of the following measure
on $\mathbb{R}^{n}$, 
\[
\mu_{t}\coloneqq\mathscr{L}_{*}\left(\frac{\Omega_{\phi_{t}}}{\int_{X}\Omega_{\phi_{t}}}\right)=\frac{e^{-\phi_{t}}dx}{\int_{\mathbb{R}^{n}}e^{-\phi_{t}}},\ \textrm{when}\ t\rightarrow\infty.
\]
We need to rescale it to obtain a nontrivial limit. Let $\nu_{t}$
be the push-forward of $\mu_{t}$ by the shrinking transformation
$x\mapsto\frac{1}{t}x$ on $\mathbb{R}^{n}$. 
\begin{thm}
\label{thm: lim slope of toric geo ray}Let $P$ be a convex body
containing $0$ in its interior. Let $0\leq t\mapsto\phi_{t}\in\mathcal{E}^{1}(P)$
be a toric geodesic ray given by (\ref{eq: toric geo ray}) and suppose
$\phi_{0}\in\textrm{PSH}_{b}(P)$. Then we have 
\begin{equation}
\lim_{t\rightarrow\infty}\frac{1}{t}L(\phi_{t})=-f(0).\label{eq: limit slope of L}
\end{equation}
Then the limit slope of the Ding functional is equal to
\[
\lim_{t\rightarrow\infty}\frac{1}{t}D(\phi_{t})=-f(0)+\frac{1}{V_{P}}\int_{P}fdy.
\]
Moreover, the probability measure $\nu_{t}$ weakly converges to a
measure $\nu_{\infty}$ which is supported on $\{x\in\mathbb{R}^{n}\mid f^{*}(x)=\inf_{\mathbb{R}^{n}}f^{*}\}$. 
\end{thm}

\begin{proof}
Let $f_{t}=\frac{1}{t}\phi_{0}^{*}+f$, it converges to $f$ as $t\rightarrow\infty$.
Then $\phi_{t}(x)=(tf_{t})^{*}=tf_{t}^{*}(\frac{x}{t})$. Let $g_{t}=f_{t}-f_{t}(0)$
and $g=f-f(0)$. Since $g_{t},\ g\in L^{1}(P)$, thus $g_{t}^{*}$,
$g^{*}\in\mathcal{E}^{1}(P)$, so they are proper over $\mathbb{R}^{n}$.
We have 
\[
g_{t}\rightarrow g\coloneqq f-f(0),\ \inf_{\mathbb{R}^{n}}g_{t}^{*}=-g_{t}(0)=0,\ \inf_{\mathbb{R}^{n}}g^{*}=-g(0)=0.
\]
Rewrite the integral 
\begin{equation}
\int_{\mathbb{R}^{n}}e^{-\phi_{t}}dx=t^{n}\int e^{-tf_{t}^{*}}dx=t^{n}e^{tf(0)+\phi_{0}^{*}(0)}\int e^{-tg_{t}^{*}}dx.\label{eq: integral ident}
\end{equation}
By (\ref{eq: sup-inf}) and $\phi_{0}\in\textrm{PSH}_{b}(P)$, we
have 
\begin{equation}
\sup_{\mathbb{R}^{n}}\left|g_{t}^{*}-g^{*}\right|=\sup_{P}\left|g_{t}-g\right|=\frac{C}{t},\label{eq: control g and gt}
\end{equation}
where $C>0$ depends on $\phi_{0}$. It follows that $e^{-C}\int e^{-tg^{*}}\leq\int e^{-tg_{t}^{*}}\leq e^{C}\int e^{-tg^{*}}$.
Thus it suffices to estimate $\int e^{-tg^{*}}$. We will show that
there exists $\delta,\ C'>0$ such that 
\begin{equation}
\delta t^{-n}\leq\int_{\mathbb{R}^{n}}e^{-tg^{*}}dx\leq C',\ \mathrm{for}\ t\geq1.\label{eq: decay estim of integral}
\end{equation}
The upper bound is trivial, since $g^{*}\geq0$, $\int e^{-tg^{*}}\leq\int e^{-g^{*}}$. 

For the lower bound, we suppose that $g^{*}$ attains the minimum
value $0$ at $x_{0}$, we claim there is a $A>0$ such that 
\begin{equation}
g^{*}(x)\leq A\left|x-x_{0}\right|,\ \textrm{on}\ \mathbb{R}^{n}.\label{eq: upper b of g*}
\end{equation}
Firstly, since convex functions are always locally Lipschitz, there
exists $r,\ A_{1}>0$ such that 
\[
g^{*}(x)\leq A_{1}\left|x-x_{0}\right|,\ \textrm{on}\ B_{r}(x_{0}).
\]
Secondly, there is a $C>0$ such that $g^{*}\leq\psi_{P}+C$ on $\mathbb{R}^{n}$.
Hence when $\left|x-x_{0}\right|\geq r$, we have 
\[
g^{*}(x)\leq\sup_{y\in P}\left\langle x,y\right\rangle +C\leq\sup_{y\in P}\left|y\right|\cdot\left|x-x_{0}\right|+C_{1}\leq A_{2}\left|x-x_{0}\right|,
\]
where $A_{2}=\sup_{y\in P}\left|y\right|+C_{1}r^{-1}$. Now (\ref{eq: upper b of g*})
follows if we take $A=\max\{A_{1},A_{2}\}$. It gives the lower bound
of (\ref{eq: decay estim of integral}), 
\[
\int e^{-tg^{*}}dx\geq\int e^{-At\left|x-x_{0}\right|}dx=\delta t^{-n},\ \delta\coloneqq\int_{\mathbb{R}^{n}}e^{-A\left|x\right|}dx.
\]
Finally, (\ref{eq: limit slope of L}) follows from (\ref{eq: integral ident})
and (\ref{eq: decay estim of integral}). 

For the last statement, it is easy to see 
\[
\nu_{t}=\frac{e^{-tg_{t}^{*}}dx}{\int_{\mathbb{R}^{n}}e^{-tg_{t}^{*}}}.
\]
Let $\triangle\coloneqq\{f^{*}=\inf_{\mathbb{R}^{n}}f^{*}\}=\{g^{*}=0\}$,
which is a compact convex set. It suffices to show that for any closed
subset $F\subset\mathbb{R}^{n}$ (bounded or not) such that $F\cap\triangle=\emptyset$,
we always have $\nu_{t}(F)\rightarrow0$. In fact, set $\gamma=\inf_{F}g^{*}>0$,
by (\ref{eq: control g and gt}) and (\ref{eq: decay estim of integral}),
we have 
\[
\nu_{t}(F)\leq e^{2C}\frac{\int_{F}e^{-tg^{*}}}{\int_{\mathbb{R}^{n}}e^{-tg^{*}}}\leq e^{2C}\delta^{-1}t^{n}e^{-t\gamma}\int_{F}e^{-t\left(g^{*}-\gamma\right)}\rightarrow0,
\]
note $\int_{F}e^{-t\left(g^{*}-\gamma\right)}$ is uniformly bounded
above since the integrand is decreasing in $t$.
\end{proof}

\subsection{Test configurations and Berman's lct formula }

We review Berman's lct formula for the limit slopes of the Ding functional
along geodesic rays induced from test configurations. 

Let $X$ be a Fano variety. A \textsl{test configuration} $(\mathcal{X},\mathcal{L})$
for $(X,-K_{X})$ is a proper normal variety $\mathcal{X}$ with a
$\mathbb{Q}$-line bundle $\mathcal{L}$ and a morphism $\pi:\mathcal{X}\rightarrow\mathbb{P}^{1}$
satisfying: there is a linearized $\mathbb{C}^{*}$-action on $(\mathcal{X},\mathcal{L})$
such that $\pi$ is equivariant with respect to the usual $\mathbb{C}^{*}$-action
on $\mathbb{P}^{1}$; the restricted family $(\mathcal{X},\mathcal{L})|_{\mathbb{P}^{1}\backslash\{0\}}$
is equivariantly isomorphic to $(X\times\mathbb{P}^{1}\backslash\{0\},-\pi_{1}^{*}K_{X})$
equipped with a trivial $\mathbb{C}^{*}$-action (i.e. only acts on
$\mathbb{P}^{1}\backslash\{0\}$). Let $\mathcal{X}_{\tau}\coloneqq\pi^{-1}(\tau)$
and $\mathcal{L}_{\tau}\coloneqq\mathcal{L}|_{\mathcal{X}_{\tau}}$. 

For a test-configuration $(\mathcal{X},\mathcal{L})$ and a bounded
psh metric $\phi_{0}$ on $-K_{X}$, via an envelop construction (see
\cite{Berman} 2.4.), we obtain a geodesic ray $0\leq t\mapsto\phi_{t}\in\textrm{PSH}_{b}(-K_{X})$
started from $\phi_{0}$. We call $\{\phi_{t}\}$ the geodesic ray
induced from $(\mathcal{X},\mathcal{L})$, and is also called Phong-Sturm's
geodesic ray. 

For the Monge-Amp\`{e}re energy $E(\phi)$, by Theorem 3.6 \cite{BHJ2},
we have 
\begin{equation}
\lim_{t\rightarrow\infty}\frac{1}{t}E(\phi_{t})=E^{\textrm{NA}}(\mathcal{X},\mathcal{L})\coloneqq\frac{\mathcal{L}^{n+1}}{(n+1)\left(K_{X}^{-1}\right)^{n}},\label{eq: NA MA energy}
\end{equation}
where $\mathcal{L}^{n+1}$ is the intersection number and $E^{\textrm{NA}}$
is called the non-Archimedean (NA) Monge-Amp\`{e}re energy. 

\subsubsection{\label{subsec: def of klt}Log discrepancies and singularity types}

Let $(Y,B)$ be a \textsl{pair}, i.e. a normal variety $Y$ with a
$\mathbb{Q}$-Weil divisor $B$ such that $K_{(Y,B)}\coloneqq K_{Y}+B$
is $\mathbb{Q}$-Cartier. For a birational morphism $f:Y'\rightarrow Y$,
there is a unique divisor $B'$ on $Y'$ such that $K_{Y'}+B'\sim f^{*}K_{(Y,B)}$
and $f_{*}B'=B$. The \textit{relatively canonical divisor} is defined
as $K_{Y'/(Y,B)}\coloneqq K_{Y'}-f^{*}K_{(Y,B)}\sim-B'$. 

A \textit{divisorial valuation} $v$ on $Y$ is a valuation on the
function field of $Y$ with the form $v=c\cdot\mathrm{ord}_{F}$,
where $c>0$ and $\mathrm{ord}_{F}$ is the vanishing order along
a prime divisor $F\subset Y'$ on some birational model $Y'\rightarrow Y$.
The \textsl{log discrepancy} of $v$ on the pair $(Y,B)$ is defined
by
\[
A_{(Y,B)}(v)\coloneqq c\left(1-\mathrm{ord}_{F}B'\right),
\]
where $B'$ is defined above. A pair $(Y,B)$ is said to be \textit{sub
log canonical} (sublc) if $A_{(Y,B)}(v)\geq0$ for any divisorial
valuation $v$ on $Y$. In particular, sublc implies the coefficients
of $B$ are $\leq1$. A pair $(Y,B)$ is said to be \textit{log canonical}
(lc) if it is sublc and $B\geq0$. A pair $(Y,B)$ is said to be \textit{Kawamata
log terminal} (klt) if $A_{(Y,B)}(v)>0$ for all $v$ and $B\geq0$.
When $Y$ is $\mathbb{Q}$-Gorenstein (i.e. $K_{Y}$ is $\mathbb{Q}$-Cartier),
we say $Y$ has log terminal singularities if the pair $(Y,0)$ is
klt. 

We will use the following property of toric varieties, see \cite{Kollar Santa-Cruz}
Proposition 3.7. For any toric variety $X$ with the torus action
by $T$, let $B=X\backslash T$ be the sum of all invariant prime
divisors, then the pair $(X,B)$ is log canonical. Since $B\sim-K_{X}$,
if $K_{X}$ is $\mathbb{Q}$-Cartier, then $X$ has log terminal singularities. 

\subsubsection{\label{subsec: Berman lct formula}Berman's lct formula}

Let $(\mathcal{X},\mathcal{L})$ be a test-configuration for $(X,-K_{X})$.
Since $\mathcal{L}\cong-K_{\mathcal{X}/\mathbb{P}^{1}}$ on $\mathcal{X}\backslash\mathcal{X}_{0}$,
there exists a unique $\mathbb{Q}$-Weil divisor $D$ supported on
$\mathcal{X}_{0}$ such that $\mathcal{L}+K_{\mathcal{X}/\mathbb{P}^{1}}\sim-D$.
In particular, $K_{\mathcal{X}}+D$ is $\mathbb{Q}$-Cartier. We define
\[
L^{\textrm{NA}}(\mathcal{X},\mathcal{L})\coloneqq\textrm{lct}_{(\mathcal{X},D)}(\mathcal{X}_{0})-1,
\]
where $\textrm{lct}_{(\mathcal{X},D)}(\mathcal{X}_{0})\coloneqq\sup\{c\in\mathbb{Q}\mid(\mathcal{X},D+c\mathcal{X}_{0})\ \mathrm{is\ sublc}\}$
is called the \textsl{log canonical threshold} (lct). The \textit{Berman-Ding
invariant} is defined as 
\[
D^{\textrm{NA}}(\mathcal{X},\mathcal{L})\coloneqq-E^{\textrm{NA}}(\mathcal{X},\mathcal{L})+L^{\textrm{NA}}(\mathcal{X},\mathcal{L}).
\]
It is also called the non-Archimedean Ding functional. 
\begin{thm}
\label{thm: Berman's lct}\cite{Berman} Let $X$ be a Fano variety,
and $(\mathcal{X},\mathcal{L})$ be a test-configuration for $(X,-K_{X})$.
For any $\phi_{0}\in\textrm{PSH}_{b}(-K_{X})$, let $t\mapsto\phi_{t}\in\textrm{PSH}_{b}(-K_{X})$
be the induced geodesic ray started from $\phi_{0}$. Then we have
\[
\lim_{t\rightarrow\infty}\frac{1}{t}L(\phi_{t})=L^{\mathrm{NA}}(\mathcal{X},\mathcal{L}),\ \lim_{t\rightarrow\infty}\frac{1}{t}D(\phi_{t})=D^{\mathrm{NA}}(\mathcal{X},\mathcal{L}).
\]
\end{thm}

\subsection{Applying Berman's lct formula to toric test configurations}

Let $X$ be a toric Fano variety associated to polytope $P$. A \textsl{toric
test-configuration} $(\mathcal{X},\mathcal{L})$ for $(X,-K_{X})$
means the total space $\mathcal{X}$ is also a toric variety with
a $T_{N}\times\mathbb{C}^{*}$-action, where the first factor acts
on each fiber of $\mathcal{X}\rightarrow\mathbb{P}^{1}$ and coincides
with the $T_{N}$-action on $X$, the second factor gives the structure
$\mathbb{C}^{*}$-action of test-configurations. By \cite{Dona toric},
toric test-configurations can be built from the rational PL convex
functions on $P$. We recall this construction. Let 
\[
f(y)=\max_{\beta\in\Gamma}\sum_{k=1}^{n}b_{\beta}^{k}y_{k}+b_{\beta}^{0},\ \textrm{where}\ b_{\beta}^{k},\ b_{\beta}^{0}\in\mathbb{Q}.
\]
Assume there are no redundant affine functions in the above definition
of $f$. Take a large integer $L$, we define 
\[
\square=\{(y,z)\in M_{\mathbb{R}}\oplus\mathbb{R}\mid y\in P,\ 0\leq z\leq L-f(y)\},
\]
which is a $(n+1)$-dimensional polytope. We can image it to be a
chalet over $P$ with the roof given by $L-f$. 

We determine the primitive normal vectors of the facets of $\square$.
Take a sufficiently divisible $r\in\mathbb{N}$ such that $r\square$
is a lattice polytope and $rb_{\beta}^{0}\in\mathbb{Z}$ for all $\beta\in\Gamma$.
Let $m_{\beta}$ be the least common multiple of the denominators
of rational numbers $\{b_{\beta}^{k}\}_{k=1}^{n}$. Then for $\beta\in\Gamma$,
\[
n_{\beta}=(-m_{\beta}b_{\beta}^{1},\cdots,-m_{\beta}b_{\beta}^{n},-m_{\beta})\in\mathbb{Z}^{n+1}\cong N\oplus\mathbb{Z}
\]
are the inner normal vectors of the roof facets of $\square$. It
is easy to see $n_{\beta}$ are primitive. Set $d_{\beta}\coloneqq-m_{\beta}rb_{\beta}^{0}+m_{\beta}rL\in\mathbb{Z}$.
Then the defining affine functions for $r\square$ are 
\begin{eqnarray}
\textrm{Roof:} &  & \left\langle n_{\beta},(y,z)\right\rangle +d_{\beta}\geq0,\ \textrm{for all}\ \beta\in\Gamma;\nonumber \\
\textrm{Wall:} &  & \left\langle (q_{\alpha},0),(y,z)\right\rangle +r\geq0,\ \textrm{for all}\ \alpha\in\Lambda;\label{eq: Defining toric TC}\\
\textrm{Floor:} &  & z\geq0.\nonumber 
\end{eqnarray}

The lattice polytope $r\square$ induces a toric variety $\mathcal{X}$
with a line bundle $r\mathcal{L}$, where $\mathcal{L}$ is a $\mathbb{Q}$-line
bundle. Moreover, $(\mathcal{X},\mathcal{L})$ has a linearized $T_{N}\times\mathbb{C}^{*}$-action.
There is an equivariant morphism $\pi:\mathcal{X}\rightarrow\mathbb{P}^{1}$
induced from the projection $N\oplus\mathbb{Z}\rightarrow\mathbb{Z}$.
In a summary, $(\mathcal{X},\mathcal{L})$ is a toric test-configuration
for $(X,-K_{X})$ given by datum $(P,f,L)$. If we change $L$ to
$L+c$, then $\mathcal{L}$ would be $\mathcal{L}\otimes\pi^{*}\mathcal{O}(c)$.
Next we apply Berman's lct formula to compute its Berman-Ding invariant. 
\begin{thm}
\label{thm: NA Ding of toric TC}Let $X$ be a toric Fano variety
associated to polytope $P$. Let $(\mathcal{X},\mathcal{L})$ be the
toric test-configuration for $(X,-K_{X})$ given by datum $(P,f,L)$,
where $f$ is a rational PL convex function on $P$ and $L$ is a
large integer. Then we have $L^{\mathrm{NA}}(\mathcal{X},\mathcal{L})=L-f(0)$
and $E^{\mathrm{NA}}(\mathcal{X},\mathcal{L})=V_{P}^{-1}\int_{P}(L-f)dy$.
The Berman-Ding invariant is equal to 
\[
D^{\mathrm{NA}}(\mathcal{X},\mathcal{L})=-f(0)+\frac{1}{V_{P}}\int_{P}fdy.
\]
\end{thm}

\begin{proof}
Let $\{D_{\beta}\}_{\beta\in\Gamma}$, $\{D_{\alpha}\}_{\alpha\in\Lambda}$
and $D_{\infty}$ be the invariant prime divisors on $\mathcal{X}$
corresponding to the (roof, wall and floor) facets of $\square$,
where $D_{\infty}=\pi^{-1}(\infty)\cong X_{P}$ is the fiber over
$\infty$. By the defining functions (\ref{eq: Defining toric TC})
of $\square$, we know $\mathcal{L}$ is given by a $\mathbb{Q}$-Cartier
divisor $\sum_{\alpha\in\Lambda}D_{\alpha}+\sum_{\beta\in\Gamma}r^{-1}d_{\beta}D_{\beta}$.
The canonical divisor is $-K_{\mathcal{X}}=\sum_{\alpha}D_{\alpha}+\sum_{\beta}D_{\beta}+D_{\infty}$.
Moreover, Proposition  6.2.7 \cite{CLS toric} tells us the pulling-back
divisor is given by

\[
-\pi^{*}K_{\mathbb{P}^{1}}=\sum_{\beta\in\Gamma}m_{\beta}D_{\beta}+D_{\infty},\ \mathcal{X}_{0}=\pi^{-1}(0)=\sum_{\beta\in\Gamma}m_{\beta}D_{\beta}.
\]
It follows that 
\[
-K_{\mathcal{X}/\mathbb{P}^{1}}=-K_{\mathcal{X}}+\pi^{*}K_{\mathbb{P}^{1}}=\sum_{\alpha\in\Lambda}D_{\alpha}+\sum_{\beta\in\Gamma}\left(1-m_{\beta}\right)D_{\beta}.
\]
Now we take 
\[
D\coloneqq\sum_{\beta\in\Gamma}\left(1-m_{\beta}-\frac{1}{r}d_{\beta}\right)D_{\beta}=\sum_{\beta\in\Gamma}\left(1+(b_{\beta}^{0}-L-1)m_{\beta}\right)D_{\beta},
\]
it satisfies $-D\sim\mathcal{L}+K_{\mathcal{X}/\mathbb{P}^{1}}$ and
$\textrm{Supp}D\subset\mathcal{X}_{0}.$ Next we compute 
\[
\textrm{lct}_{(\mathcal{X},D)}(\mathcal{X}_{0})\coloneqq\sup\{c\in\mathbb{Q}\mid(\mathcal{X},D+c\mathcal{X}_{0})\ \textrm{is\ sublc}\}.
\]
Let 
\[
D_{c}\coloneqq D+c\mathcal{X}_{0}=\sum_{\beta\in\Gamma}\left(1+(b_{\beta}^{0}-L-1+c)m_{\beta}\right)D_{\beta}.
\]
If $(\mathcal{X},D_{c})$ is sublc, then each coefficient of $D_{c}$
is $\leq1$. It implies 
\[
c\leq L+1-\max_{\beta\in\Gamma}b_{\beta}^{0}=L+1-f(0).
\]
On the other hand, by Proposition 3.7 \cite{Kollar Santa-Cruz}, we
know that as a toric variety, the pair $(\mathcal{X},-K_{\mathcal{X}})$
is always log canonical, so $(\mathcal{X},-K_{\mathcal{X}})$ is sublc.
For any rational $c<L+1-f(0)$, we have $D_{c}\leq-K_{\mathcal{X}}$
(by comparing their coefficients). Then by Property 3.4.1.5 in \cite{Kollar Santa-Cruz}
(monotonicity of discrepancy), $(\mathcal{X},-K_{\mathcal{X}})$ is
sublc implies $(\mathcal{X},D_{c})$ is sublc. Thus we have 
\[
L^{\textrm{NA}}(\mathcal{X},\mathcal{L})=\textrm{lct}_{(\mathcal{X},D)}(\mathcal{X}_{0})-1=L-f(0).
\]
For the term $E^{\textrm{NA}}(\mathcal{X},\mathcal{L})$, since we
have 
\[
\mathcal{L}^{n+1}=(n+1)!\mathrm{vol}(\square)=(n+1)!\int_{P}(L-f)dy,
\]
put this into (\ref{eq: NA MA energy}), also note $c_{1}(X)^{n}=n!V_{P}$.
Then the result follows. 
\end{proof}
Since the Berman-Ding invariant is equal to the limit slope by Theorem
\ref{thm: Berman's lct}, thus the above result should match with
Theorem \ref{thm: lim slope of toric geo ray}. In fact, for the toric
test-configuration given by datum $(P,f,L)$, the induced toric geodesic
ray is 
\begin{equation}
\phi_{t}=\left(\phi_{0}^{*}+t(f-L)\right)^{*}.\label{eq: ind toric ray f-L}
\end{equation}

\begin{rem}
In \cite{Boucksom Jonsson}, Boucksom-Jonsson studied the asymptotics
of a family of volume forms associated to a degeneration, of which
test-configurations are special cases. We compare our asymptotic estimate
for the total mass $\int_{X}\Omega_{\phi_{t}}$ with that of Theorem
A in \cite{Boucksom Jonsson}. 

Consider the geodesic ray (\ref{eq: ind toric ray f-L}), backing
to the proof of Theorem \ref{thm: lim slope of toric geo ray}, by
(\ref{eq: mu-exp}) and (\ref{eq: integral ident}) (replace $f$
by $f-L$), we have 
\[
\int_{X}\Omega_{\phi_{t}}=(8\pi)^{n}\int_{\mathbb{R}^{n}}e^{-\phi_{t}}dx=(8\pi)^{n}e^{\phi_{0}^{*}(0)}\cdot e^{t\left(f(0)-L\right)}t^{n}\int_{\mathbb{R}^{n}}e^{-tg_{t}^{*}}dx.
\]
Recall $\int e^{-tg_{t}^{*}}$ can be controlled by $\int e^{-tg^{*}}$
and vice versa. In this case, $f$ is piecewise linear, we can make
a better estimate than (\ref{eq: decay estim of integral}), 
\[
\int_{\mathbb{R}^{n}}e^{-tg^{*}}\sim\frac{1}{t^{n-d}},\ d=\dim\{f^{*}=\inf_{\mathbb{R}^{n}}f^{*}\}.
\]
It follows that there is a $C>0$ such that 
\[
C^{-1}e^{-tL^{\textrm{NA}}(\mathcal{X},\mathcal{L})}t^{d}\leq\int_{X}\Omega_{\phi_{t}}\leq Ce^{-tL^{\textrm{NA}}(\mathcal{X},\mathcal{L})}t^{d},
\]
for $t\gg1$. This matches with the estimate in Theorem A \cite{Boucksom Jonsson}.
Actually, if we rename the complex parameter $t$ in \cite{Boucksom Jonsson}
by $\tau$, then our $t=-\log\left|\tau\right|^{2}$. The invariant
$\kappa_{\min}$ in \cite{Boucksom Jonsson} is equal to $L^{\textrm{NA}}(\mathcal{X},\mathcal{L})=L-f(0)$. 
\end{rem}

\subsection{\label{subsec: D vs K}D-stability and K-stability}

In this section, we compare Ding stability with K-stability for toric
Fano varieties. Let $X$ be a toric Fano variety associated to polytope
$P$. For a toric test-configuration $(\mathcal{X},\mathcal{L})$
for $(X,-K_{X})$ given by datum $(P,f,L)$, Donaldson \cite{Dona toric}
computed the so called the Donaldson-Futaki invariant of $(\mathcal{X},\mathcal{L})$,
\[
M^{\textrm{NA}}(\mathcal{X},\mathcal{L})=\frac{1}{V_{P}}\left(\int_{\partial P}fd\sigma-n\int_{P}fdy\right),
\]

where $d\sigma$ is a measure on $\partial P$ determined by the normal
vectors of the facts of $P$, it satisfies $d\sigma(\partial P)=nV_{P}$.
Later, it is realized that this is actually the limit slope of the
Mabuchi functional along the induced geodesic ray, see \cite{BHJ2,Li Chi}.
Thus it is also called the non-Archimedean Mabuchi functional, as
the notation indicates. 

For a bounded convex function $f$ on $P$, we define 
\[
D^{\textrm{NA}}(f)\coloneqq-f(0)+\frac{1}{V_{P}}\int_{P}fdy,\ M^{\textrm{NA}}(f)\coloneqq\frac{1}{V_{P}}\left(\int_{\partial P}fd\sigma-n\int_{P}fdy\right).
\]

\begin{defn}
For a toric Fano variety $X_{P}$ associated to polytope $P$. We
say $X_{P}$ is \textit{D-stable} if $D^{\mathrm{NA}}(f)\geq0$ for
all rational PL convex $f$ on $P$, and the equality holds only for
affine $f$. \textit{K-stable} is defined in the same way by replacing
$D^{\mathrm{NA}}$ by $M^{\mathrm{NA}}$. Note the full name of D-stability
(K-stability) should be Ding-polystability (K-polystability) with
respect to toric test-configurations. 
\end{defn}

We have an elementary inequality for these two invariants. See \cite{BHJ1}
7.8 for a comparison for general varieties. 
\begin{prop}
For any bounded convex function $f$ on $P$, we have 
\begin{equation}
D^{\mathrm{NA}}(f)\leq M^{\mathrm{NA}}(f).\label{eq: Ding vs M}
\end{equation}
The equality holds if and only if $f$ is radially affine, i.e. $f(sz)=(1-s)f(0)+sf(z)$
for any $s\in[0,1]$ and $z\in\partial P$. 
\end{prop}

\begin{proof}
Consider the difference 
\[
M^{\textrm{NA}}(f)-D^{\textrm{NA}}(f)=\frac{1}{V_{P}}\left(\int_{\partial P}\left(f-f(0)\right)d\sigma-(n+1)\int_{P}\left(f-f(0)\right)dy\right).
\]
We define a radially affine function $\bar{f}$ on $P$, for any $s\in[0,1]$
and $z\in\partial P$, let 
\[
\bar{f}(sz)=(1-s)f(0)+sf(z).
\]
Since $f$ is convex, we have $f\leq\bar{f}$. Next we decompose $P$
into the union of cones $\{P_{\alpha}\mid\alpha\in\Lambda\}$, each
of them has the vertex $0$ and takes a facet of $P$ as the base.
For each $P_{\alpha}$, we show that 
\begin{equation}
\int_{P_{\alpha}}\left(\bar{f}-f(0)\right)dy=\frac{1}{n+1}\int_{\left(\partial P\right)_{\alpha}}\left(f-f(0)\right)d\sigma,\label{eq:integral relation}
\end{equation}
where $\left(\partial P\right)_{\alpha}\subset\partial P$ is the
base of $P_{\alpha}$. In fact, after adding a constant, we assume
$f(0)=0$. Consider a map 
\[
\rho:\left(\partial P\right)_{\alpha}\times(0,1]\rightarrow P_{\alpha}\backslash\{0\},\ \rho(z,s)=sz.
\]
By the definition of the boundary measure $d\sigma$, we see $\rho^{*}dy=s^{n-1}d\sigma\wedge ds$.
Hence we have 
\[
\int_{P_{\alpha}}\bar{f}dy=\int_{\left(\partial P\right)_{\alpha}\times(0,1]}sf(z)\rho^{*}dy=\int_{\left(\partial P\right)_{\alpha}}f(z)d\sigma\cdot\int_{0}^{1}s^{n}ds=\frac{1}{n+1}\int_{\left(\partial P\right)_{\alpha}}fd\sigma.
\]
Then taking the sum of (\ref{eq:integral relation}) over $\alpha\in\Lambda$,
we obtain 
\[
(n+1)\int_{P}\left(\bar{f}-f(0)\right)dy=\int_{\partial P}\left(f-f(0)\right)d\sigma.
\]
Combing this with $M^{\textrm{NA}}-D^{\textrm{NA}}$, we obtain 
\[
M^{\textrm{NA}}(f)-D^{\textrm{NA}}(f)=\frac{n+1}{V_{P}}\int_{P}\left(\bar{f}-f\right)dy\geq0.
\]
\end{proof}
\begin{cor}
Let $X_{P}$ be the toric Fano variety given by a polytope $P$, then
the following are equivalent to each other: (1) $X_{P}$ is D-stable;
(2) $X_{P}$ is K-stable; (3) the barycenter of $P$ is the origin. 
\end{cor}

\begin{proof}
(1) to (2) directly follows from (\ref{eq: Ding vs M}). (2) to (3),
we know $M^{\textrm{NA}}(f)=0$ for all affine $f$. Then the above
Proposition implies $D^{\textrm{NA}}(y_{i})=M^{\textrm{NA}}(y_{i})=0$,
thus the barycenter is the origin. (3) to (1) follows from Jensen's
inequality. 
\end{proof}
The seminar work \cite{Wang-Zhu} shows that a toric Fano manifold
$X_{P}$ admits a KE metric if and only if the barycenter of $P$
is the origin. This condition is further equivalent to the Ding-polystability
(also K-polystability) of $X_{P}$. It is a special case of the Yau-Tian-Donaldson
correspondence. In the next section, we extend this correspondence
to Mabuchi solitons. 

\section{Relative Ding stability}

In this section, for toric Fano varieties, we introduce the stability
notion corresponding to Mabuchi solitons. We first introduce the modified
Ding functional $D_{\ell}$ which takes Mabuchi solitons as its critical
points, then use its limit slopes $D_{\ell}^{\textrm{NA}}$ to define
the relative Ding stability. 

\subsection{Modified Ding functionals}

Let $X$ be a toric Fano variety associated to  polytope $P$. We
take $K\supset S_{N}$ to be a maximal compact subgroup of $\textrm{Aut}^{0}(X)$,
and a $K$-invariant smooth metric $\phi$ on $-K_{X}$ with curvature
$\omega=i\partial\bar{\partial}\phi$. For each $v\in\mathfrak{s}\coloneqq\mathrm{Lie}(S_{N})$,
there is an affine function $f_{v}$ on $P$ such that 
\[
\theta_{v}(\omega)=f_{v}\circ\mathfrak{m}_{\phi},\ \int_{P}f_{v}dy=0.
\]
Recall $\theta_{v}(\omega)$ is the Hamiltonian of $v$ (w.r.t. $\omega$)
and $\mathfrak{m}_{\phi}=d\phi$ is the moment map induced by $\phi$.
We call $f_{v}$ the \textsl{affine Hamiltonian} for $v$, which is
independent of the choice of $\phi$. 

The Futaki-Mabuchi's bilinear form (see Remark \ref{rem: extrm vec})
is given by 
\[
B_{K}(v,w)=\frac{1}{V_{P}}\int_{P}f_{v}f_{w}dy,\ \textrm{for}\ v,w\in\mathfrak{s}.
\]
Let $\theta_{P}$ be the affine Hamiltonian of the extremal vector
field $Z_{K}$. Since $Z_{K}$ belongs to the center of $\mathfrak{k}$
(see \cite{Futaki Mabuchi}), and $S_{N}$ is a maximal torus in $K$
(by Corollary 4.7 \cite{Cox}), we know $Z_{K}\in\mathfrak{s}$. Recall
$Z_{K}$ is characterized by the condition: $F(v)=-B_{K}(v,Z_{K})$
for all $v\in\mathfrak{s}$. For rational $v$, the Futaki invariant
$-F(v)$ is equal to $M^{\textrm{NA}}(f_{v})=D^{\textrm{NA}}(f_{v})$,
since that $f_{v}$ induces a product toric test-configuration. Hence
the condition for $Z_{K}$ follows that 
\begin{equation}
-f(0)+\frac{1}{V_{P}}\int_{P}fdy\eqqcolon D^{\textrm{NA}}(f)=\frac{1}{V_{P}}\int_{P}f\theta_{P}dy,\ \textrm{for all affine}\ f.\label{eq: define extremal affine}
\end{equation}
We can determine $\theta_{P}$ by this condition. Actually, for any
subset $P\subset\mathbb{R}^{n}$ with nonempty interior, applying
the Riesz representation theorem to the space of affine functions
with $L^{2}$-inner product, we know there exists a unique affine
function $\theta_{P}$ satisfying (\ref{eq: define extremal affine}).
Note that $\int_{P}\theta_{P}dy=0$ by taking $f\equiv1$. 
\begin{defn}
(Ricci affine function $\ell_{P}$) Let $P\subset\mathbb{R}^{n}$
be a convex body. The unique affine function $\theta_{P}$ satisfying
(\ref{eq: define extremal affine}) is called the \textit{extremal
affine function} associated to $P$. Let $\ell_{P}\coloneqq V_{P}^{-1}\left(1-\theta_{P}\right)$,
it satisfies 
\[
-f(0)+\int_{P}f\ell_{P}\ dy=0,\ \textrm{for all affine}\ f.
\]
We call $\ell_{P}$ the \textit{Ricci affine function} associated
to $P$, sometimes denote it by $\ell$ for simplicity. Note that
$\int_{P}\ell dy=1$ and $\ell(0)=\int_{P}\ell^{2}dy>0$ by taking
$f=\ell$. 
\end{defn}

We explain the name of $\ell_{P}$. If the metric $\phi$ is a Mabuchi
soliton, it should satisfy $1-e^{h(\phi)}=\theta_{P}(d\phi)$, so
$e^{h(\phi)}=V_{P}\cdot\ell(d\phi)$, this explains the name. Moreover,
by the formula (\ref{eq: toric Ricci poten}) for $e^{h(\phi)}$,
we see $\phi$ satisfies a Monge-Amp\`{e}re type equation, 
\[
\ell(d\phi)\textrm{MA}(\phi)=\frac{e^{-\phi}}{\int_{\mathbb{R}^{n}}e^{-\phi}},\ \textrm{on}\ \mathbb{R}^{n}.
\]
Next we modify the Ding functional such that its critical points satisfy
the above equation. 
\begin{defn}
Let $P\subset\mathbb{R}^{n}$ be a convex body such that $0\in\mathrm{Int}(P)$,
and $g(y)$ is a bounded continuous function on $P$ with $\int_{P}g=1$. 

(1) We define the $g$-modified Ding functional $D_{g}:\mathcal{E}^{1}(P)\rightarrow\mathbb{R}$
by
\[
D_{g}(\phi)\coloneqq-\log\int_{\mathbb{R}^{n}}e^{-\phi}-E_{g}(\phi),
\]
where $E_{g}(\phi)=-\int_{P}\phi^{*}gdy$ is the $g$-modified Monge-Amp\`{e}re
energy. When $g=\ell_{P}$, we call $D_{\ell}$ the \textsl{modified
Ding functional}. Note that $D_{\ell}$ is invariant under the transformation
$\phi\mapsto\phi(\cdot+a)+c$, where $a\in\mathbb{R}^{n}$ and $c\in\mathbb{R}$. 

(2) We define the $g$-modified Berman-Ding invariant by
\[
D_{g}^{\textrm{NA}}(f)\coloneqq-f(0)+\int_{P}f\ gdy,
\]
where $f$ is an integrable convex function on $P$. When $g=\ell_{P}$,
we call $D_{\ell}^{\textrm{NA}}$ the relative Berman-Ding invariant,
it vanishes for affine functions. 
\end{defn}

The motivation of this definition is as follows. For (1), by Proposition
2.13 \cite{Berman Berndtsson}, if $\phi$ is a critical point of
$D_{g}$, then it satisfies 
\[
\textrm{MA}_{g}(\phi)=\frac{e^{-\phi}}{\int_{\mathbb{R}^{n}}e^{-\phi}},\ \textrm{on}\ \mathbb{R}^{n}
\]
in the sense of Alexandrov. For (2), $D_{g}^{\textrm{NA}}$ is the
limit slope of $D_{g}$ along the geodesic ray (\ref{eq: toric geo ray}).
Finally, since $\phi^{*}(0)=-\inf_{\mathbb{R}^{n}}\phi$, we can rewrite
$D_{g}$ as follows, 
\[
D_{g}(\phi)=-\log\int_{\mathbb{R}^{n}}e^{-\left(\phi-\inf_{\mathbb{R}^{n}}\phi\right)}+D_{g}^{\textrm{NA}}(\phi^{*}).
\]
This shows a close relation of $D_{g}$ to its non-Archimedean counterpart. 

\subsection{Relative Ding stability }

As with Ding stability, we define the relative Ding stability in terms
of $D_{\ell}^{\textrm{NA}}$. Moreover, to build a YTD type correspondence,
we also need to study the coercivity of $D_{\ell}$, which is closely
related to the existence of Mabuchi solitons. So in this section,
we discuss the relations between stability, coercivity and the positivity
of $\ell_{P}$. Our proof of the Proposition \ref{prop: D-stable}
and \ref{prop: almost cycle} will use the partial coercivity of $D_{\ell}$,
whose proof is a little cumbersome, so we postpone it to the next
section (Theorem \ref{thm: Partial coercive}). The proof of Theorem
\ref{thm: Partial coercive} can be read before the following discussions
for stability. 

\subsubsection{The relative D-(un)stability}
\begin{defn}
\label{def: rD-stable and unstable}Let $P\subset\mathbb{R}^{n}$
be a polytope containing $0$ in its interior. 

(1) We say $P$ is \textit{relatively D-stable} if $D_{\ell}^{\textrm{NA}}(f)\geq0$
for any rational PL convex function $f$ on $P$, and the equality
holds only when $f$ is affine. 

(2) We say $P$ is \textit{relatively D-unstable} if $D_{\ell}^{\textrm{NA}}(f)<0$
for some rational PL convex function $f$. 

For a toric Fano variety $X_{P}$, we say it is relatively D-(un)stable
if $P$ is. Note that the full name of this notion should be the relative
Ding-polystability with respect to toric test-configurations. 
\end{defn}

We have not define relative D-semistability (i.e. $D_{\ell}^{\textrm{NA}}(f)\geq0$
for all rational PL convex $f$), since it coincides with relative
D-stability, as the below proposition shows. 
\begin{prop}
\label{prop: D-stable}Let $P\subset\mathbb{R}^{n}$ be a polytope
containing $0$ in its interior. The following are equivalent to each
other:
\begin{enumerate}
\item $D_{\ell}$ is bounded from below on $\mathcal{E}^{1}(P)$.
\item $D_{\ell}^{\mathrm{NA}}(f)\geq0$ for any rational PL convex function
$f$ on $P$.
\item $\ell_{P}\geq0$ on $P$. 
\end{enumerate}
When $\ell_{P}\geq0$ on $P$, if a convex $f$ satisfies $D_{\ell}^{\mathrm{NA}}(f)=0$,
then $f$ is affine. Hence $P$ is relatively D-(un)stable if and
only if $\inf_{P}\ell_{P}\geq0$ ($<0$).
\end{prop}

\begin{proof}
$(1)\Rightarrow(2)$. For a rational PL convex $f$, take the limit
slope of $D_{\ell}$ along the geodesic ray $\phi_{t}=\left(tf\right)^{*}$,
then $D_{\ell}\geq-C$ follows that $D_{\ell}^{\textrm{NA}}(f)\geq0$. 

$(2)\Rightarrow(3)$. If it is not true, $P^{-}\coloneqq P\cap\{\ell\leq0\}$
would be a sub-polytope and $0\notin P^{-}$, since $\ell(0)>0$.
Let $f\coloneqq\max(-\ell,0)$, by $(2)$ we have $D_{\ell}^{\mathrm{NA}}(f)\geq0$.
However, $D_{\ell}^{\textrm{NA}}(f)=\int_{P^{-}}f\ell=-\int_{P^{-}}\ell^{2}<0$,
it is a contradiction. 

$(3)\Rightarrow(1)$ follows from Theorem \ref{thm: Partial coercive}.
The next statement follows from the equality condition of Jensen's
inequality. For the last equivalence relation, when $\inf_{P}\ell_{P}<0$,
we have seen $f\coloneqq\max(-\ell,0)$ satisfies $D_{\ell}^{\textrm{NA}}(f)<0$.
Making a perturbation, we can find a rational PL convex $f'$ such
that $D_{\ell}^{\textrm{NA}}(f')<0$, so by definition, $P$ is relatively
D-unstable. 
\end{proof}
\begin{rem}
Note that in the Definition \ref{def: rD-stable and unstable} of
relative D-stability, we only use rational PL convex functions. Then
by the above proposition, it automatically implies $D_{\ell}^{\mathrm{NA}}(f)>0$
for all non-affine convex $f$. For relative K-stability, it is not
clear whether these two versions of stability are equivalent to each
other. 
\end{rem}

\subsubsection{The reduced J-functionals}

Recall in the toric setting, the J-functional is given by 
\[
J(\phi)\coloneqq\sup_{\mathbb{R}^{n}}(\phi-\psi_{P})-E(g)=\frac{1}{V_{P}}\int_{P}\phi^{*}dy-\inf_{P}\phi^{*},\ \phi\in\mathcal{E}^{1}(P).
\]
Note the first term of the usual J-functional is $\int_{X}(\phi-\psi_{P})\omega_{0}^{n}$,
and these two versions can be controlled by each other. In order to
define coercivity in the present of a Lie group action, Hisamoto \cite{Hisa toric}
introduced the reduced J-functional, which has been appeared in \cite{Zhou-Zhu Adv}
with a different form. 
\begin{defn}
\label{def:def of reduced J}(1) Let $P$ be a convex body containing
$0$ in its interior, and $g\geq0$ is a bounded continuous function
on $P$ with $\int_{P}g=1$. We set $J_{g}(\phi)\coloneqq\int_{P}\phi^{*}gdy-\inf_{P}\phi^{*}$,
called the $g$-modified J-functional. Let $\phi_{a}(x)\coloneqq\phi(x+a)$
be a translation of $\phi$. The $g$-modified \textsl{reduced J-functional}
is defined by 
\begin{equation}
J_{g}^{T}(\phi)\coloneqq\inf_{a\in\mathbb{R}^{n}}J_{g}(\phi_{a})=\inf\left\{ \int_{P}(\phi^{*}-l)\ gdy-\inf_{P}(\phi^{*}-l)\middle|\textrm{affine}\ l\right\} ,\label{eq: def reduced J}
\end{equation}
for $\phi\in\mathcal{E}^{1}(P)$. In above, we use $\phi_{a}^{*}=\phi^{*}-\left\langle a,\cdot\right\rangle $.
When $g\equiv V_{P}^{-1}$, we simply denote $J^{T}=J_{V_{P}^{-1}}^{T}$. 

(2) For an integrable convex function $f$ on $P$, we set $J(f)\coloneqq\frac{1}{V_{P}}\int_{P}fdy-\inf_{P}f$,
and 
\[
J^{T}(f)\coloneqq\inf\left\{ \frac{1}{V_{P}}\int_{P}(f-l)dy-\inf_{P}(f-l)\middle|\textrm{affine}\ l\right\} .
\]
\end{defn}

The below proposition gives a more explicit formula for $J_{g}^{T}$. 
\begin{prop}
\label{prop: property of reduced J}Let $P$ be a convex body containing
$0$ in its interior. 

(1) Let $g\geq0$ is a bounded continuous function on $P$ with $\int_{P}g=1$.
Let $\mathfrak{c}_{g}$ be the barycenter of $(P,gdy)$, then for
any $\phi\in\mathcal{E}^{1}(P)$, we have 
\begin{equation}
J_{g}^{T}(\phi)=\int_{P}\phi^{*}g-\phi^{*}(\mathfrak{c}_{g})\geq0.\label{eq: formula J^T}
\end{equation}
Moreover, in the definition (\ref{eq: def reduced J}), $l$ attains
the infimum if and only if $\nabla l\in\partial\phi^{*}(\mathfrak{c}_{g})$.
If a convex function $f\in L^{1}(P)$ satisfies $J^{T}(f)=0$, then
$f$ is affine.

(2) There exists $C>0$ only depending on $P$, such that $C^{-1}\int_{P}\phi^{*}\leq J^{T}(\phi)\leq C\int_{P}\phi^{*}$
for all $\phi\in\mathcal{E}^{1}(P)$ normalized at $0$, i.e. satisfying
$\phi\geq\phi(0)=0$. 
\end{prop}

\begin{proof}
1. By the definition, 
\begin{eqnarray*}
J_{g}^{T}(\phi) & = & \inf_{a\in\mathbb{R}^{n}}\left[\int_{P}\left(\phi^{*}-\left\langle a,y\right\rangle \right)gdy-\inf_{P}\left(\phi^{*}-\left\langle a,y\right\rangle \right)\right]\\
 & = & \inf_{a\in\mathbb{R}^{n}}\left[\int_{P}\phi^{*}g-\left\langle a,\mathfrak{c}_{g}\right\rangle +\phi(a)\right]\\
 & = & \int_{P}\phi^{*}g-\sup_{a\in\mathbb{R}^{n}}\left[\left\langle a,\mathfrak{c}_{g}\right\rangle -\phi(a)\right]=\int_{P}\phi^{*}g-\phi^{*}(\mathfrak{c}_{g})\geq0.
\end{eqnarray*}
The last $\geq$ is due to Jensen's inequality. We see $l$ attains
the infimum in (\ref{eq: def reduced J}) iff $a=\nabla l$ attains
the above supremum, so if and only if $\mathfrak{c}_{g}\in\partial\phi(a)$,
equivalently $a\in\partial\phi^{*}(\mathfrak{c}_{g})$. The last statement
follows from the equality condition of Jensen's inequality. 

2. Firstly, we take a continuous function $g$ such that $\int_{P}g=1$,
$C^{-1}\leq g\leq C$ and the barycenter of $(P,gdy)$ is $0$, where
$C>0$ only depends on $P$. 

By definition, $J_{g}(\phi)=\int_{P}\left(\phi^{*}-\inf_{P}\phi^{*}\right)gdy$
for $\phi\in\mathcal{E}^{1}(P)$, so 
\[
V_{P}C^{-1}\cdot J(\phi)\leq J_{g}(\phi)\leq V_{P}C\cdot J(\phi).
\]
Replace $\phi$ by $\phi_{a}$ and take the infimum over $a\in\mathbb{R}^{n}$,
we obtain 
\[
V_{P}C^{-1}\cdot J^{T}(\phi)\leq J_{g}^{T}(\phi)\leq V_{P}C\cdot J^{T}(\phi),\ \textrm{for all}\ \phi\in\mathcal{E}^{1}(P).
\]
If $\phi$ is normalized at $0$, we have $\phi^{*}\geq\phi^{*}(0)=0$,
then $J_{g}^{T}(\phi)=\int_{P}\phi^{*}g$ by (\ref{eq: formula J^T}).
Finally, by $C^{-1}\int_{P}\phi^{*}\leq\int_{P}\phi^{*}g\leq C\int_{P}\phi^{*}$,
the desired inequality follows. 
\end{proof}
\begin{rem}
By (\ref{eq: formula J^T}), we see the limit slope of $J^{T}$ along
the geodesic ray (\ref{eq: toric geo ray}) is equal to $J^{T}(f)$,
thus the non-Archimedean counterpart of $J^{T}$ coincides with itself.
There is another interesting coincidence: when the barycenter of $(P,dy)$
is $0$, $J^{T}(f)$ formally coincides with the Berman-Ding invariant
$D^{\textrm{NA}}(f)$. 
\end{rem}

\subsubsection{The uniformly relative D-stability}
\begin{defn}
\label{def: Uniform rD-stable}Let $P$ be a convex body containing
$0$ in its interior. 

(1) For a $\delta>0$, we say the modified Ding functional $D_{\ell}$
is \textit{$\delta$-coercive} if there is a constant $C_{\delta}$
such that 
\[
D_{\ell}(\phi)\geq\delta J^{T}(\phi)-C_{\delta},\ \textrm{for all}\ \phi\in\mathcal{E}^{1}(P).
\]
We call the supremum of such $\delta$ the\textit{ coercivity threshold}
of $D_{\ell}$. We say $D_{\ell}$ is \textit{coercive} if it is $\delta$-coercive
for some $\delta>0$.

(2) For a $\delta>0$, we say $P$ is \textit{$\delta$-uniformly}
relative D-stable if 
\[
D_{\ell}^{\textrm{NA}}(f)\geq\delta J^{T}(f),
\]
for all rational PL convex $f$ on $P$ (then for all integrable convex
$f$ by approximation). We call the supremum of such $\delta$ the
\textit{stability threshold} of $P$. We say $P$ is \textsl{uniformly
relative D-stable} if it is \textit{$\delta$-}uniformly relative
D-stable for some $\delta>0$. For a toric Fano variety $X_{P}$,
we say it is uniformly relative D-stable if $P$ is. 

(3) We say $P$ is \textsl{non-uniformly} relative D-stable if it
is relatively D-stable (Definition \ref{def: rD-stable and unstable})
but not uniformly. 
\end{defn}

\begin{rem}
Since $D_{\ell}$ is invariant by adding affine functions, thus by
Proposition \ref{prop: property of reduced J} (2), the coercivity
of $D_{\ell}$ is equivalent to the existence of some $\delta,C>0$
such that $D_{\ell}(\phi)\geq\delta\int_{P}\phi^{*}-C$ for all normalized
$\phi\in\mathcal{E}^{1}(P)$. 
\end{rem}

\begin{prop}
\label{prop: almost cycle}Let $P\subset\mathbb{R}^{n}$ be a polytope
containing $0$ in its interior. Then

(1) If $D_{\ell}$ is $\delta$-coercive, then $P$ is $\delta$-uniformly
D-stable. 

(2) If $P$ is $\delta$-uniformly D-stable, then $\ell_{P}\geq\delta\cdot V_{P}^{-1}$
on $P$. 

(3) If $\inf_{P}\ell_{P}>0$, then for any $0<\epsilon<1$, $D_{\ell}$
is $(\epsilon V_{P}\cdot\inf_{P}\ell_{P})$-coercive. 
\end{prop}

\begin{proof}
(1) Taking the limit slopes along geodesic rays. By (\ref{eq: formula J^T}),
the limit slopes of $J^{T}$ are given by itself. 

(2) By (\ref{eq: formula J^T}), $\delta$-uniformly D-stable means
that 
\begin{equation}
-f(0)+\int_{P}f\ell\geq\delta\left(\frac{1}{V_{P}}\int_{P}fdy-f(\mathfrak{c}_{P})\right),\label{eq: delta D-stable}
\end{equation}
for all PL convex $f$, where $\mathfrak{c}_{P}$ is the barycenter
of $(P,dy)$. If the conclusion is not true, then $P_{\delta}\coloneqq P\cap\{\ell\leq\delta\cdot V_{P}^{-1}\}$
would be a polytope with nonempty interior. We can find a nonzero
PL convex function $f\geq0$ such that $\mathrm{supp}f\subset P_{\delta}$
and $f(0)=f(\mathfrak{c}_{P})=0$, then 
\[
\int_{P}f(\ell-\frac{\delta}{V_{P}})dy=\int_{P_{\delta}}f(\ell-\frac{\delta}{V_{P}})dy<0.
\]
This contradicts with (\ref{eq: delta D-stable}). 

(3) By Theorem \ref{thm: Partial coercive}, for any $\epsilon\in(0,1)$,
there exists a $C_{\epsilon}$ such that $D_{\ell}(\phi)\geq\epsilon J_{\ell}^{T}(\phi)-C_{\epsilon}$
for all $\phi\in\mathcal{E}^{1}(P)$. Note that $J_{\ell}^{T}(\phi)$
is equal to 
\[
\inf_{a\in\mathbb{R}^{n}}\int_{P}\left(\phi_{a}^{*}-\inf_{P}\phi_{a}^{*}\right)\ell\geq\inf_{P}\ell\cdot V_{P}\cdot\inf_{a\in\mathbb{R}^{n}}\frac{1}{V_{P}}\int_{P}\left(\phi_{a}^{*}-\inf_{P}\phi_{a}^{*}\right)=\inf_{P}\ell\cdot V_{P}\cdot J^{T}(\phi).
\]
Thus $D_{\ell}(\phi)\geq\epsilon V_{P}\cdot\inf_{P}\ell\cdot J^{T}(\phi)-C_{\epsilon}$
for all $\phi\in\mathcal{E}^{1}(P)$. 
\end{proof}
Note the statements in the above proposition almost constitute a cycle.
Immediately, we have 
\begin{cor}
\label{cor: unif D-stable}Let $P$ be a polytope containing $0$
in its interior. Then the following are equivalent to each other. 

(1) $D_{\ell}$ is coercive. (2) $P$ is uniformly relative D-stable.
(3) $\inf_{P}\ell_{P}>0$. 

Moreover, if they hold, the coercivity threshold and the stability
threshold both equal to $V_{P}\cdot\inf_{P}\ell_{P}$. 
\end{cor}

Combining Proposition \ref{prop: D-stable} and Corollary \ref{cor: unif D-stable},
we see relative Ding stability can be completely detected by the sign
of $\inf_{P}\ell_{P}$. In a summary, we have 
\begin{cor}
\label{cor: summary of rD-stability}Let $X_{P}$ be a toric Fano
variety associated to  polytope $P$ with Ricci affine function $\ell_{P}$,
then 
\begin{enumerate}
\item $X_{P}$ is relatively D-stable if and only if $\inf_{P}\ell_{P}\geq0$. 
\item $X_{P}$ is (non-)uniformly relative D-stable if and only if $\inf_{P}\ell_{P}>0$
$(=0)$. 
\item $X_{P}$ is relatively D-unstable if and only if $\inf_{P}\ell_{P}<0$. 
\end{enumerate}
\end{cor}

\subsubsection{\label{subsec: rD vs rK}Relative D-stability v.s. relative K-stability}

For a toric Fano Manifold $X_{P}$, the existence of extremal metrics
in class $c_{1}(X)$ is conjectured to be equivalent to (uniformly)
relative K-stability, see \cite{Gabor PHD}. We will not state the
various versions of relative K-stability, which are not clear to be
equivalent to each other. We only state a relation of it to relative
D-stability. 

Recall the relative Donaldson-Futaki invariant (for $c_{1}(X)$) is
given by 
\[
M_{\theta}^{\textrm{NA}}(f)\coloneqq\frac{1}{V_{P}}\left(\int_{\partial P}fd\sigma-\int_{P}(n+\theta_{P})fdy\right).
\]
By inequality (\ref{eq: Ding vs M}), we have $M_{\theta}^{\textrm{NA}}(f)\geq D_{\ell}^{\textrm{NA}}(f)$
for all bounded convex function $f$ on $P$. 

If $X_{P}$ is relatively D-stable, i.e. $\inf_{P}\ell_{P}\geq0$,
then $M_{\theta}^{\textrm{NA}}(f)\geq0$ for all bounded convex $f$
and the equality holds only when $f$ is affine (by the equality condition
for Jensen's inequality), which is called \textsl{strongly} relative
K-stable in \cite{Yotsutani ZhouB}. 

If $X_{P}$ is uniformly relative D-stable, i.e. $\inf_{P}\ell_{P}>0$,
then there is a $\delta>0$ such that $M_{\theta}^{\textrm{NA}}(f)\geq\delta J^{T}(f)$
for all bounded convex $f$, which is called uniformly relative K-stable
in \cite{Hisa toric}. 

This gives a simple criterion for determining the K-stability, actually
has been obtained in \cite{Zhou-Zhu Adv,Yotsutani ZhouB} without
mention Ding stability. An interesting point is that Mabuchi solitons
and extremal metrics seem to have no connection from the view of their
equations. Moreover, this criterion is not necessary, e.g. $\mathbb{P}\left(\mathcal{O}_{\mathbb{P}^{2}}\oplus\mathcal{O}_{\mathbb{P}^{2}}(2)\right)$
is relative D-unstable (see Example \ref{exam: optimal 3fold}), but
it is relatively K-stable since it admits an extremal metric in $c_{1}(X)$
(according to a computation by Yotsutani \cite{private comm} using
the criterion in \cite{Apostolov}. 

\subsection{Examples }
\begin{example}
\label{exa: 2 polygons}In dimension 2, modulo the action of $\mathrm{GL}_{2}(\mathbb{Z})$,
there are 16 isomorphism classes of reflexive Fano polygons, see Figure
\ref{fig: 2-polygons} for their dual polygons $P$. They give all
Gorenstein toric Fano surfaces, which are orbifolds with isolated
singularities corresponding to the vertices marked by solid $\diamondsuit$.
In each pair: \{8A,4A\}, \{9,3\}, \{8B,4B\}, \{7A,5A\},\{4C,8C\} and
\{7B,5B\}, two polygons are dual to each other up to a $\mathrm{GL}_{2}(\mathbb{Z})$-transformation.
All the rests are self-dual. 

The crossing lines in the figure are the set $\{\ell_{P}=0\}$. When
$\inf_{P}\ell_{P}>0$, these lines are not drawn. The dashed line
in 6D will be explained in Example \ref{exa: optimal 6d}. We see
\{4C,8C\} are non-uniformly relative D-stable; \{7B,5B,6C,6D\} are
relatively D-unstable; and all the rests are uniformly relative D-stable.
The barycenter of \{3,9,4A,8A,6A\} is the origin, hence they are D-stable.
We also note that two dual polygons have same relative Ding stability.
It is natural to ask whether this holds in each dimension. 
\begin{figure}[h]
\includegraphics[scale=0.3]{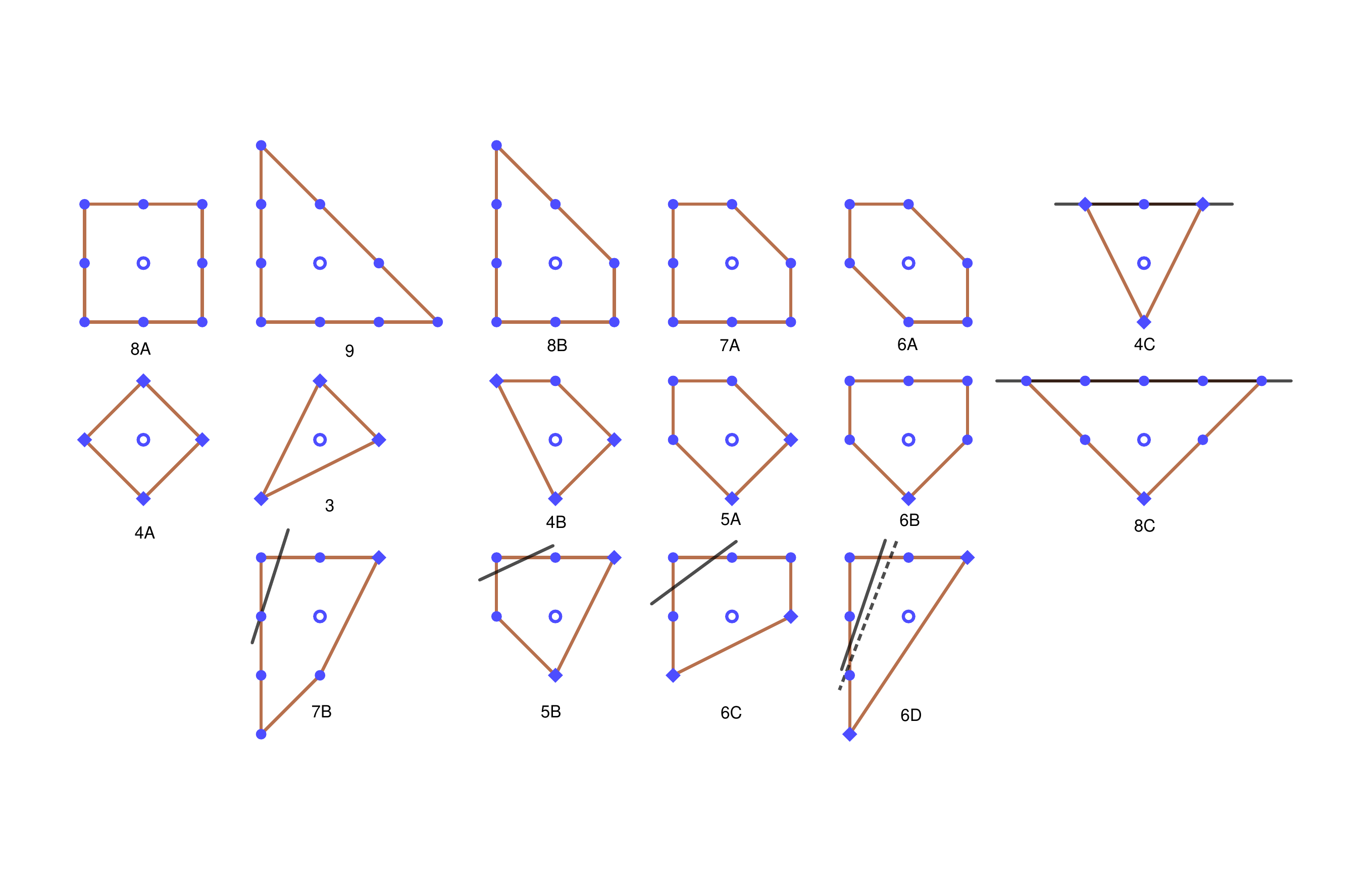}

\caption{\label{fig: 2-polygons}The polygon $P$ of all Gorenstein toric Fano
surfaces.}
\end{figure}
\end{example}

\begin{example}
In dimension 3, there are 18 isomorphism classes of smooth toric Fano
3-folds, see Table-1 in \cite{Yotsutani 4fold} (first version) for
a list of them, and the labels used below. In \cite{Yotsutani ZhouB},
Yotsutani and Zhou found some criterions to determine the relative
K-stability of toric varieties. They computed $\theta_{P}$ for these
3-folds, see Table 2 in \cite{Yotsutani ZhouB}. Since $\ell_{P}=V_{P}^{-1}\left(1-\theta_{P}\right)$,
by their computation results, we can determine the relative Ding stability
of these 3-folds. 

In 18 classes, $\{\mathbb{P}^{3},\mathcal{B}_{4},\mathcal{C}_{3},\mathcal{C}_{5},\mathcal{F}_{1}\}$
have vanishing Futaki invariant, hence they are uniformly relative
D-stable and admit KE metrics; $\{\mathcal{B}_{2},\mathcal{B}_{3},\mathcal{C}_{1},\mathcal{C}_{4},\mathcal{E}_{3},\mathcal{E}_{4},\mathcal{F}_{2}\}$
are uniformly relative D-stable with nonvanishing Futaki invariant;
the rests $\{\mathcal{B}_{1},\mathcal{C}_{2},\mathcal{D}_{1},\mathcal{D}_{2},\mathcal{E}_{1},\mathcal{E}_{2}\}$
are relatively D-unstable. 

For the classification in dimension 4, see \cite{Yotsutani 4fold}. 
\end{example}

\begin{rem}
In dimension $\leq4$, there is no smooth toric Fano manifold that
is non-uniformly relative D-stable. We do not know whether such manifolds
exist in higher dimension. We also do not know whether there exists
a rational polytope $P$ such that $P\cap\{\ell_{P}=0\}$ is a face
of $P$ with codimension $>1$. 
\end{rem}

\subsection{\label{sec: Exist result for unifm}Existence results in the uniformly
stable case}

In the case of $\inf_{P}\ell>0$, the equation $\ell(d\phi)\textrm{MA}(\phi)=e^{-\phi}$
had been studied in \cite{Berman Berndtsson}. Theorem 1.1 in \cite{Berman Berndtsson}
yields, 
\begin{thm}
\cite{Berman Berndtsson} (1) Let $P$ be a convex body containing
$0$ in its interior and suppose that $\inf_{P}\ell_{P}>0$. Then
there exists a smooth and strictly convex function $\phi\in\textrm{PSH}_{b}(P)$
satisfying 
\[
\ell(d\phi)\textrm{MA}(\phi)=e^{-\phi},\ \textrm{on}\ \mathbb{R}^{n},
\]
 and $d\phi:\mathbb{R}^{n}\rightarrow\mathrm{Int}(P)$ is a diffeomorphism.
Such $\phi$ is unique up to translations. 

(2) Let $X_{P}$ be a toric Fano variety associated to  polytope $P$
with $\inf_{P}\ell_{P}>0$. The above solution $\phi$ gives a $S_{N}$-invariant
weak Mabuchi soliton on $X_{P}$, in the sense that its potential
is continuous on $X_{P}$ and smooth on $T_{N}\subset X_{P}$. 
\end{thm}

When $X_{P}$ is smooth (i.e. $P$ is Delzant), we can obtain a smooth
Mabuchi soliton via the continuation method in \cite{Wang-Zhu}. 
\begin{thm}
\cite{Wang-Zhu} Let $X$ be a toric Fano manifold associated to  Delzant
polytope $P$ with $\inf_{P}\ell_{P}>0$, then $X$ admits a smooth
$S_{N}$-invariant Mabuchi soliton. 
\end{thm}

\begin{proof}
The proof is almost a duplication of \cite{Wang-Zhu}, so we just
sketch it. Consider the following equations parameterized by $t\in[0,1]$
(compare with (\ref{eq: equation for MS})), 
\[
\left(1-\theta_{Z}(\omega_{u_{t}})\right)\omega_{u_{t}}^{n}=e^{h_{r}-tu_{t}}\omega_{r}^{n},
\]
where the reference metric $\omega_{r}$ is $S_{N}$-invariant. Let
$\mathcal{I}\subset[0,1]$ be the set of $t$ such that the equation
is solvable. Our aim is to show $1\in\mathcal{I}$. $0\in\mathcal{I}$
is shown in the Appendix of \cite{Liyan Zhou} and the openness of
$\mathcal{I}$ follows from the implicit function theorem. All that
remains is to show $\mathcal{I}$ is closed, it suffices to establish
a uniform estimate for $u_{t}$. By the $S_{N}$-invariance, the above
equation can be reduced to 
\[
\ell(d\phi_{t})\textrm{MA}(\phi_{t})=e^{-t\phi_{t}-(1-t)\phi_{r}},\ \textrm{on}\ \mathbb{R}^{n},
\]
where $\phi_{t}$ and $\phi_{r}$ are the corresponding convex functions
of $\omega_{u_{t}}$ and $\omega_{r}$, we have $u_{t}=\phi_{t}-\phi_{r}$.
Since the term $\ell(d\phi_{t})$ is bounded above and below (away
from zero), the next arguments in \cite{Liyan Zhou} also work in
our setting, then we obtain a uniform estimate for $u_{t}$. 
\end{proof}
\begin{rem}
Above theorem is also covered by Theorem 1.5 in \cite{Liyan Zhou},
where they considered the more general group compactifications. As
in \cite{Shiyl Zhu}, the above theorem can be extended to toric Fano
orbifolds, then it will include the uniformly relative D-stable Fano
surfaces in Example \ref{exa: 2 polygons}. 
\end{rem}

\section{Partial coercivity}

In this section, when $\ell_{P}\geq0$ on $P$, we establish a partial
coercive inequality for $D_{\ell}$. 
\begin{thm}
\label{thm: Partial coercive} (Partial coercivity) Let $P\subset\mathbb{R}^{n}$
be a polytope containing $0$ in its interior. Suppose that $\ell_{P}\geq0$
on $P$. Then for any $\epsilon\in(0,1)$ there is a constant $C_{\epsilon}>0$
depending on $P$ and $\epsilon$ such that
\begin{equation}
D_{\ell}(\phi)\geq\epsilon\int_{P}\phi^{*}\ell-C_{\epsilon}\label{eq: partial coercive}
\end{equation}
for all $\phi\in\mathcal{E}^{1}(P)$ such that $\phi\geq\phi(0)=0$.
Then we have $D_{\ell}(\phi)\geq\epsilon J_{\ell}^{T}(\phi)-C_{\epsilon}$
for all $\phi\in\mathcal{E}^{1}(P)$. In particular, $D_{\ell}$ is
bounded from below. 
\end{thm}

Its proof needs a subsolution of the equation for Mabuchi solitons. 
\begin{lem}
\label{lem: subsolut}(Existence of Subsolution) Let $P\subset\mathbb{R}^{n}$
be a polytope containing $0$ in its interior. Suppose that $\ell_{P}\geq0$
on $P$. Then there is a smooth and strictly convex $\phi\in\textrm{PSH}_{b}(P)$
such that $C\cdot\ell_{P}(d\phi)\mathrm{MA}(\phi)\geq e^{-\phi}$
on $\mathbb{R}^{n}$ for some constant $C>0$. 
\end{lem}

\begin{proof}
Let $\varphi(x)=\log\sum_{p}e^{\left\langle x,p\right\rangle }$,
where the sum is over all vertices of $P$. It can be directly checked
that $\varphi$ is smooth, strictly convex and $\psi_{P}\leq\varphi(x)\leq\psi_{P}+C_{P}$,
on $\mathbb{R}^{n}$ for some constant $C_{P}$. We show that $\phi(x)\coloneqq R\varphi(\frac{x}{R})$
could satisfy the desired conditions when $R$ is sufficiently large. 

Firstly, $\phi$ is still in $\textrm{PSH}_{b}(P)$ and we have 
\[
\ell(d\phi)(Rx)=\frac{\sum_{p}\ell(p)e^{\left\langle x,p\right\rangle }}{\sum_{p}e^{\left\langle x,p\right\rangle }},\ \mathrm{MA}(\phi)(Rx)=\frac{1}{R^{n}}\mathrm{MA}(\varphi)(x)
\]
and $e^{-\phi(Rx)}=\left(\sum_{p}e^{\left\langle x,p\right\rangle }\right)^{-R}$.
In order to make $C\cdot\ell(d\phi)\mathrm{MA}(\phi)\geq e^{-\phi}$
holds on $\mathbb{R}^{n}$, by taking the logarithm, it is sufficient
to make 
\begin{equation}
\log\sum_{p}\ell(p)e^{\left\langle x,p\right\rangle }+\log\mathrm{MA}(\varphi)(x)+(R-1)\varphi(x)-n\log R\label{eq: lemma quan}
\end{equation}
is bounded from below. If we can bound the first two terms from below
by some $-C_{1}\psi_{P}-C_{2}$ (where $C_{i}>0$ only depends on
$P$), then by $\varphi\geq\psi_{P}$, we can take $R$ to be sufficiently
large so that the whole (\ref{eq: lemma quan}) is bounded from below. 

(1) For the first term, we have 
\[
\sum_{p}\ell(p)e^{\left\langle x,p\right\rangle }=\sum_{\ell(p)>0}\ell(p)e^{\left\langle x,p\right\rangle }\ge\delta\sum_{\ell(p)>0}e^{\left\langle x,p\right\rangle },
\]
where the second sum is over all vertices $p$ such that $\ell(p)>0$.
Note that $\log\sum_{\ell(p)>0}e^{\left\langle x,p\right\rangle }\geq\psi_{Q}$,
where $\psi_{Q}(x)=\sup_{y\in Q}\left\langle x,y\right\rangle $ and
$Q$ is the convex hull of all vertices $p$ of $P$ such that $\ell(p)>0$.
Taking a $K\gg1$ such that $-Q\subset K\cdot P$, then $\psi_{Q}\geq-K\psi_{P}$.
Finally, we have 
\[
\log\sum_{p}\ell(p)e^{\left\langle x,p\right\rangle }\geq-K\psi_{P}+\log\delta,\ \textrm{on}\ \mathbb{R}^{n}.
\]

(2) For the second term: $\log\mathrm{MA}(\varphi)$, we compute
\[
\varphi_{ij}(x)=\sum_{p}(p_{i}-\bar{p}_{i})(p_{j}-\bar{p}_{j})\frac{e^{\left\langle x,p\right\rangle }}{\sum_{p}e^{\left\langle x,p\right\rangle }},\ \bar{p}_{i}=\frac{\sum_{p}p_{i}e^{\left\langle x,p\right\rangle }}{\sum_{p}e^{\left\langle x,p\right\rangle }}.
\]
Our aim is to bound $\log\det[\varphi_{ij}]$ from below by some $-C_{1}\psi_{P}-C_{2}$.
First we can take the factor $\sum_{p}e^{\left\langle x,p\right\rangle }$
(in $\varphi_{ij}$) out of the determinant. Next we give a lower
bound for the eigenvalues of matrix $H_{ij}(x)\coloneqq\sum_{p}(p_{i}-\bar{p}_{i})(p_{j}-\bar{p}_{j})e^{\left\langle x,p\right\rangle }$.
For any $\xi\in\mathbb{R}^{n}$, we have 
\[
\sum_{i,j}H_{ij}\xi_{i}\xi_{j}=\sum_{p}\left\langle p-\bar{p},\xi\right\rangle ^{2}e^{\left\langle x,p\right\rangle }\geq e^{\min_{p}\left\langle x,p\right\rangle }\sum_{p}\left\langle p-\bar{p},\xi\right\rangle ^{2},
\]
where $\bar{p}=(\bar{p}_{1},\cdots,\bar{p}_{n})$ depends on $x$.
We claim there is a $\gamma>0$ depending on $P$ such that $\sum_{p}\left\langle p-\bar{p},\xi\right\rangle ^{2}\geq\gamma\left|\xi\right|^{2}$
for all $x,\ \xi\in\mathbb{R}^{n}$. Consider a continuous function
$Q(\eta,\xi)=\sum_{p}\left\langle p-\eta,\xi\right\rangle ^{2}$ defined
on compact space $Y\coloneqq P\times\mathbb{S}^{n-1}$. Then we must
have $\inf_{Y}Q>0$, otherwise $Q=0$ at some point in $Y$, but this
never happens unless $P$ is contained in a hyperplane. Note that
$\bar{p}\in P$, then the claim follows from $\inf_{Y}Q>0$. 

It follows that $\sum_{i,j}H_{ij}\xi_{i}\xi_{j}\geq e^{\min_{p}\left\langle x,p\right\rangle }\gamma\left|\xi\right|^{2}$
for all $\xi\in\mathbb{R}^{n}$. This gives a lower bound for eigenvalues,
so $\det H\geq e^{n\min_{p}\left\langle x,p\right\rangle }\gamma^{n}$,
namely 
\[
\log\det H(x)\geq n\min_{p}\left\langle x,p\right\rangle +n\log\gamma,\ \textrm{for all}\ x\in\mathbb{R}^{n}.
\]
 Finally, we take a $K>0$ such that $-P\subset K\cdot P$, then 
\[
\min_{p}\left\langle x,p\right\rangle =-\sup_{p}\left\langle x,-p\right\rangle =-\psi_{-P}(x)\geq-K\psi_{P}.
\]
Now we bound $\log\det[\varphi_{ij}]$ from below by some $-C_{1}\psi_{P}-C_{2}$. 
\end{proof}
\begin{rem}
This lemma is expected to hold for general convex bodies, but we can
not take $\varphi(x)=\log\int_{P}e^{\left\langle x,y\right\rangle }dy$
and go through the same proof, since it may not belong to $\textrm{PSH}_{b}(P)$,
e.g. $P=[-1,1]$. 
\end{rem}

Once we have a subsolution, (\ref{eq: partial coercive}) follows
from the same method in Theorem 2.16 \cite{Berman Berndtsson}. 
\begin{proof}
(of Theorem \ref{thm: Partial coercive})

Step 1. By the above lemma, there is a smooth and strictly convex
$\phi_{0}\in\textrm{PSH}_{b}(P)$ such that $\int_{\mathbb{R}^{n}}e^{-\phi_{0}}=1$
and 
\[
C\ \ell_{P}(d\phi_{0})\mathrm{MA}(\phi_{0})\geq e^{-\phi_{0}},\ \textrm{on}\ \mathbb{R}^{n}.
\]
Firstly we show there exists $C_{1}>0$ depending on $\phi_{0}$ such
that
\begin{equation}
-\log\int e^{-\phi}dx\geq CE_{\ell}(\phi)-C_{1},\label{eq: rough coercive}
\end{equation}
for all $\phi\in\mathcal{E}^{1}(P)$ with $\phi(0)=0$. 

Since $\phi_{0}\in\textrm{PSH}_{b}(P)$, suppose $\psi_{P}\leq\phi_{0}+A$
for a constant $A$. By (\ref{eq: sup-inf}), $\phi(0)=0$ implies
$\phi\leq\psi_{P}$, hence $\phi-A\leq\phi_{0}$. Let $\{\phi_{t}\}_{0\leq t\leq1}$
be the toric geodesic segment connecting $\phi_{0}$ with $\phi_{1}\coloneqq\phi-A$.
Since $\phi_{0}\geq\phi_{1}$, we have $\frac{\partial\phi_{t}}{\partial t}|_{t=0}\leq0$.
Then by the convexity of $D_{\ell}$ and $D_{\ell}(\phi_{1})=D_{\ell}(\phi)$,
we have 
\begin{eqnarray*}
D_{\ell}(\phi)-D_{\ell}(\phi_{0}) & \geq & \frac{d}{dt}D_{\ell}(\phi)|_{t=0}=\int_{\mathbb{R}^{n}}\frac{\partial\phi_{t}}{\partial t}|_{t=0}\left(e^{-\phi_{0}}-\ell(d\phi_{0})\mathrm{MA}(\phi_{0})\right)\\
 & \geq & (C-1)\int_{\mathbb{R}^{n}}\frac{\partial\phi_{t}}{\partial t}|_{t=0}\ell(d\phi_{0})\mathrm{MA}(\phi_{0})\\
 & = & (C-1)\left(E_{\ell}(\phi-A)-E_{\ell}(\phi_{0})\right).
\end{eqnarray*}
The last equality since $E_{\ell}$ is affine along $\{\phi_{t}\}$.
Since $D_{\ell}(\phi)=-\log\int e^{-\phi}-E_{\ell}(\phi)$, (\ref{eq: rough coercive})
follows. 

Step 2. For any normalized $\phi\in\mathcal{E}^{1}(P)$ and $0<r<1$,
we put $\phi^{(r)}\coloneqq r\phi(\frac{x}{r})$ into (\ref{eq: rough coercive}),
since $E_{\ell}(\phi^{(r)})=-\int_{P}\phi^{(r)*}\ell=rE_{\ell}(\phi)$,
it yields 
\[
-\log\int e^{-r\phi}dx\geq rCE_{\ell}(\phi)-C_{1}+n\log r.
\]
Since $\phi\geq0$, we have 
\[
D_{\ell}(\phi)\geq-\log\int e^{-r\phi}dx-E_{\ell}(\phi)=(1-rC)\int_{P}\phi^{*}\ell-C(r).
\]
Finally, when we take $r$ to be sufficiently small, then (\ref{eq: partial coercive})
follows. 

For any $\phi\in\mathcal{E}^{1}(P)$, take $a\in\mathbb{R}^{n}$ such
that $\tilde{\phi}(x)=\phi(x+a)-\phi(a)$ is normalized at $0$. Then
we put $\tilde{\phi}$ into (\ref{eq: partial coercive}), and note
that $D_{\ell}(\tilde{\phi})=D_{\ell}(\phi)$ and 
\[
J_{\ell}^{T}(\phi)=J_{\ell}^{T}(\tilde{\phi})=\int_{P}(\tilde{\phi})^{*}\ell,
\]
by (\ref{eq: formula J^T}), the desired inequality follows. 
\end{proof}

\section{The non-uniformly stable case}

In this section, we deal with the borderline case: $\inf_{P}\ell_{P}=0$.
In Example \ref{exa: 2 polygons}, polygon 4C and 8C are such cases.
There does not exist smooth solutions in this case, since the equation
(\ref{eq: equation for MS}) implies $1-\theta_{Z}(\omega)>0$ on
$X$, so $\inf_{P}\ell>0$. However, as the below example shows, there
may exist some weak solutions with mild singularities. 
\begin{example}
(Mabuchi soliton with conic singularities) \label{exa: 8c conic MS}
The polygon 8C in Example \ref{exa: 2 polygons} is given by 
\[
P=\{(y_{1},y_{2})\mid-y_{2}\geq-1,y_{1}+y_{2}\geq-1,-y_{1}+y_{2}\geq-1\}
\]
with vertices $(0,-1)$, $(2,1)$, $(-2,1)$. The associated toric
Fano surface $X_{P}$ is $\mathbb{P}(1,1,2)$, a weighted projective
plane. It has only one singular point $[0,0,1]$, which is corresponding
to the vertex $(0,-1)$. The Ricci affine function is $\ell_{P}=\frac{3}{8}-\frac{3}{8}y_{2}$. 

We construct an explicit solution for Mabuchi solitons. Let 
\[
\phi_{FS}(x)\coloneqq\log\sum_{\textrm{vertex}\ p}e^{\left\langle x,p\right\rangle }=\log\left(e^{-x_{2}}+e^{2x_{1}+x_{2}}+e^{-2x_{1}+x_{2}}\right)
\]
and $\phi(x)\coloneqq4\phi_{FS}(\frac{x}{4})$. We see $\phi\in\textrm{PSH}_{b}(P)$
is smooth and strictly convex on $\mathbb{R}^{2}$. By Wolfram Mathematica,
we find
\begin{equation}
\ell(d\phi)\mathrm{MA}(\phi)=3e^{-\phi}.\label{eq:8c solution}
\end{equation}
Thus it gives a Mabuchi soliton at least on the dense torus $T_{N}$.
To investigate its behavior near the boundary divisors, we need to
know its Legendre dual. It is not easy to compute it directly, instead,
we define
\begin{eqnarray}
v(y) & = & \mathbf{2}(1-y_{2})\log(1-y_{2})+(1+y_{1}+y_{2})\log(1+y_{1}+y_{2})\label{eq:Guillemin solu 8c}\\
 &  & +(1-y_{1}+y_{2})\log(1-y_{1}+y_{2})\nonumber 
\end{eqnarray}
and $\psi=v^{*}$. Since $v$ is bounded on $P$, so $\psi\in\textrm{PSH}_{b}(P)$.
A direct computation (on $v$) shows that $\ell(d\psi)\mathrm{MA}(\psi)=\frac{3}{64}e^{-\psi}$.
By the uniqueness of the solution for (\ref{eq:8c solution}), there
is $a\in\mathbb{R}^{n}$ and $c\in\mathbb{R}$ such that $\psi(x)=\phi(x+a)+c$,
thus $\phi^{*}$ equals to $v$ up to an affine function (given by
$a$ and $c$). 

Note (\ref{eq:Guillemin solu 8c}) is exactly the Guillemin's model
symplectic potential except the coefficient \textbf{2}. Let $\triangle\subset X$
be the prime divisor corresponding to the facet $\{\ell=0\}=\{y_{2}=1\}$.
Note $\triangle$ does not contain the singular point. By Proposition
3.1 \cite{SongJ WangXW}, we know $\phi$ has conic singularities
along $\triangle$ with angle $\pi$. Consider the equation on $X_{P}$,
\[
\left(1-\theta_{Z}(\omega_{u})\right)\omega_{u}^{n}=e^{-u+h_{r}}\omega_{r}^{n},
\]
where $\omega_{r}$ is a smooth reference metric and $\omega_{u}=i\partial\bar{\partial}\phi$,
$\phi$ is the solution above. Comparing with smooth solutions, a
main difference is that the term $\omega_{u}^{n}$ has a pole along
$\triangle$ and $1-\theta_{Z}=\ell(d\phi)$ vanishes along $\triangle$,
so they together make the RHS is bounded from above and below (away
from zero). Taking $i\partial\bar{\partial}$ on both sides, $Ric(\omega_{u})$
contributes a current $\pi[\triangle]$, but $i\partial\bar{\partial}\log\left(1-\theta_{Z}(\omega_{u})\right)$
also contributes a current $\pi[\triangle]$ by the Poincar\'{e}-Lelong
equation, so two currents are canceled by each other. 

We can also write down an explicit solution for polygon 4C. For more
examples, take $a\in\mathbb{N}$, then 
\[
P_{a}\coloneqq\{(y_{1},y_{2})\mid-y_{2}\geq-1,ay_{1}+y_{2}\geq-1,-ay_{1}+y_{2}\geq-1\}
\]
is non-uniformly stable. Since $\ell=C_{a}(1-y_{2})$, where $C_{a}>0$
is a constant. 
\end{example}

Next we establish a general existence result in the case of $\inf_{P}\ell_{P}=0$
($\inf_{P}\ell_{P}>0$ is also included). The framework of our proof
is same to \cite{Berman Berndtsson}, via the variational method and
relying on the partial coercivity (\ref{eq: partial coercive}), but
we need some new techniques to remedy the absence of a positive lower
bound for $\ell_{P}$. 
\begin{thm}
\label{thm: Exist nonunifm}Let $P\subset\mathbb{R}^{n}$ be a polytope
containing $0$ in its interior. Suppose that $\ell_{P}\geq0$ on
$P$. Then there exists $\phi\in\textrm{PSH}_{b}(P)$ satisfying 
\[
\mathrm{MA}_{\ell}(\phi)=\frac{e^{-\phi}}{\int_{\mathbb{R}^{n}}e^{-\phi}},\ \textrm{on}\ \mathbb{R}^{n}
\]
in the sense of Alexandrov. The Legendre dual $\phi^{*}$ is H\"{o}lder
continuous on $P$ for any exponent $\gamma\in(0,1)$. In the case
of $\dim P\cap\{\ell=0\}\leq\frac{n}{2}$ (including the case that
the intersection is empty), we know $\phi$ is smooth and strictly
convex on $\mathbb{R}^{n}$, and $d\phi$ gives a diffeomorphism from
$\mathbb{R}^{n}$ to $\mathrm{Int}(P)$. 
\end{thm}

\begin{proof}
(1) Since $D_{\ell}$ is bounded below and invariant under translations,
there is a sequence of normalized functions $\{\phi_{j}\}\subset\mathcal{E}^{1}(P)$
such that $D_{\ell}(\phi_{j})\rightarrow\inf_{\mathcal{E}^{1}(P)}D_{\ell}$,
where ``normalized'' means $\phi_{j}\geq\phi_{j}(0)=0$ or equivalently
$\phi_{j}^{*}\geq\phi_{j}^{*}(0)=0$. The partial coercivity (\ref{eq: partial coercive})
implies $\int_{P}\phi_{j}^{*}\ell\leq C$. By a compactness result
(\cite{Berman Berndtsson} Proposition 2.8), after passing to a subsequence,
we assume $\phi_{j}$ converges to a $\phi\in\textrm{PSH}(P)$ uniformly
over compact subsets in $\mathbb{R}^{n}$. Note $\phi$ is also normalized. 

Since $\phi^{*}\leq\liminf_{j}\phi_{j}^{*}$, by Fatou's lemma, we
have 
\begin{equation}
\int_{P}\phi^{*}\ell\leq\int_{P}\liminf_{j}\phi_{j}^{*}\ell\leq\liminf_{j}\int_{P}\phi_{j}^{*}\ell\leq C.\label{eq: int ul<C}
\end{equation}
This implies $\phi^{*}<\infty$ a.e. on $P$. Since $\phi^{*}$ is
convex, it must have $\phi^{*}<\infty$ on $\mathrm{Int}(P)$, so
$\phi\in\mathcal{E}(P)$. In particular, $\phi$ is proper on $\mathbb{R}^{n}$,
thus $\int_{\mathbb{R}^{n}}e^{-\phi}<\infty$. Note that we can not
assert $\phi\in\mathcal{E}^{1}(P)$ yet. 

(2) We show $\lim_{j}\int e^{-\phi_{j}}=\int e^{-\phi}$. In \cite{Berman Berndtsson}
(2.17), this follows from a Skoda type uniform estimate. Here we give
another proof. 

We claim that $\{\phi_{j}^{*}\}$ are uniformly bounded above on $\frac{1}{2}P$.
As a finite-valued convex function on $\frac{1}{2}P$, $\phi_{j}^{*}$
can attain the maximum at some vertex of $\frac{1}{2}P$. At each
vertex $v$ of $\frac{1}{2}P$, we take a ball $B_{r}(v)\subset\subset\frac{2}{3}P$.
Suppose that $\ell\geq\delta>0$ on $\frac{2}{3}P$. By Jensen's inequality,
\[
\phi_{j}^{*}(v)\leq\frac{1}{\left|B_{r}\right|}\int_{B_{r}(v)}\phi_{j}^{*}\leq\frac{1}{\left|B_{r}\right|\delta}\int_{\frac{2}{3}P}\phi_{j}^{*}\ell\leq\frac{1}{\left|B_{r}\right|\delta}\int_{P}\phi_{j}^{*}\ell\leq\frac{C}{\left|B_{r}\right|\delta}.
\]
Thus we have $\sup_{\frac{1}{2}P}\phi_{j}^{*}\leq C$, this is the
claim. It follows that $\phi_{j}\geq\frac{1}{2}\psi_{P}-C$ on $\mathbb{R}^{n}$. 

For any $\epsilon>0$, since $\int_{B_{R}^{c}}e^{-\phi_{j}}\leq e^{C}\int_{B_{R}^{c}}e^{-\frac{1}{2}\psi_{P}}$,
there exists a $R(\epsilon)>0$ such that $\int_{B_{R}^{c}}e^{-\phi_{j}}<\epsilon$
for all $j\geq1$ and $\int_{B_{R}^{c}}e^{-\phi}<\epsilon$. Since
\[
\left|\int_{\mathbb{R}^{n}}e^{-\phi_{j}}-\int_{\mathbb{R}^{n}}e^{-\phi}\right|\leq\int_{B_{R}}\left|e^{-\phi_{j}}-e^{-\phi}\right|+\int_{B_{R}^{c}}e^{-\phi_{j}}+\int_{B_{R}^{c}}e^{-\phi}
\]
and $\phi_{j}$ uniformly converges to $\phi$ on $B_{R}$, we obtain
$\lim_{j}\int e^{-\phi_{j}}=\int e^{-\phi}$. 

The Skoda type estimate in \cite{Berman Berndtsson} (2.17) can be
also obtained by the above method. 

(3) Combing (\ref{eq: int ul<C}) with (2), we have 
\begin{equation}
D_{\ell}(\phi)\leq\liminf_{j}D_{\ell}(\phi_{j})=\inf_{\mathcal{E}^{1}(P)}D_{\ell}\label{eq: lsc of D}
\end{equation}
Consider a perturbation $\phi+t\eta$, where $\eta$ is a compact
supported smooth function on $\mathbb{R}^{n}$. Since it may be not
convex, we take the envelop 
\[
\phi_{t}(x)\coloneqq\sup\{\psi(x)\mid\textrm{PSH}(P)\ni\psi\leq\phi+t\eta\}=(\phi+t\eta)^{**},
\]
 where $**$ is taking the Legendre dual twice, see \cite{Berman Berndtsson}
section 2.6. Since $\eta$ is bounded, we see $\left|\phi_{t}-\phi\right|$
is also bounded on $\mathbb{R}^{n}$. Since $\phi\in\mathcal{E}(P)$
and $\int_{P}\phi^{*}\ell<\infty$, thus $\phi_{t}\in\mathcal{E}(P)$
and $\int_{P}\phi_{t}^{*}\ell<\infty$.  

Although we do not know $\phi_{t}\in\mathcal{E}^{1}(P)$, but we still
have $\inf_{\mathcal{E}^{1}(P)}D_{\ell}\leq D_{\ell}(\phi_{t})$.
Indeed, by the regularization Lemma 2.2 \cite{Berman Berndtsson},
there exists a sequence $\varphi_{t}^{k}\in\textrm{PSH}_{b}(P)$ decreasing
to $\phi_{t}$. Then by Fatou's lemma, we have $\lim_{k}D_{\ell}(\varphi_{t}^{k})=D_{\ell}(\phi_{t})$.
Since $D_{\ell}(\varphi_{t}^{k})\geq\inf_{\mathcal{E}^{1}(P)}D_{\ell}$,
the inequality follows. 

By (\ref{eq: lsc of D}), we have 
\[
D_{\ell}(\phi)\leq D_{\ell}(\phi_{t})\leq-\log\int e^{-(\phi+t\eta)}-E_{\ell}(\phi_{t})\eqqcolon h(t).
\]
Note that $h(0)=D_{\ell}(\phi)$, so $h(t)$ attains a minimum at
$t=0$. By a slight modification of Proposition 2.13 \cite{Berman Berndtsson}
(replacing $dy$ by $\ell dy$ in its proof), it follows that $E_{\ell}(\phi_{t})=-\int_{P}(\phi+t\eta)^{*}\ell$
is differentiable at $t=0$ with derivative $\int_{\mathbb{R}^{n}}\eta\mathrm{MA}_{\ell}(\phi)$.
Hence $h'(0)=0$ for all $\eta$, it implies that $\phi$ satisfies
\[
\mathrm{MA}_{\ell}(\phi)=\frac{e^{-\phi}}{\int e^{-\phi}}
\]
in the sense of Alexandrov. By adding a constant, we assume $\int e^{-\phi}=1$.
It means that for any Borel subset $E\subset\mathbb{R}^{n}$, we have
\begin{equation}
\int_{\partial\phi(E)}\ell=\int_{E}e^{-\phi}.\label{eq: Alexandrov-ell}
\end{equation}

(4) We show $\phi\in\textrm{PSH}_{b}(P)$. To simplify the notations,
we denote $u=\phi^{*}$. By (\ref{eq: int ul<C}), we know $\int_{P}u\ell<\infty$.
Moreover, for any $q\geq1$ we have 
\begin{equation}
\int_{P}\left|du\right|^{q}\ell=\int_{\mathbb{R}^{n}}\left|x\right|^{q}\mathrm{MA}_{\ell}(\phi)=\int_{\mathbb{R}^{n}}\left|x\right|^{q}e^{-\phi}<\infty.\label{eq: du^ql<C}
\end{equation}
The last finiteness since $\phi\in\mathcal{E}(P)$ and (\ref{eq: proper of full MA}).
Note the LHS integral is well-defined, since $u$ is convex hence
differentiable a.e. 

We state a weighted Poincar\'{e} inequality due to Chua and Wheeden
\cite{W Poincare} Theorem 1.1, which says for a bounded convex domain
$\Omega\subset\mathbb{R}^{n}$ with a weight function $w=\rho^{k}$
(where $k\geq0$, $\rho\geq0$ is concave on $\Omega$) and any $1\leq q\leq p<\infty$,
there exists a constant $B_{p,q}$ such that 
\[
\left(\int_{\Omega}\left|f-\frac{\int_{\Omega}fw}{\int_{\Omega}w}\right|^{q}wdx\right)^{1/q}\le B_{p,q}\cdot\left(\int_{\Omega}w\right)^{\frac{1}{q}-\frac{1}{p}}\mathrm{diam}(\Omega)\left(\int_{\Omega}\left|df\right|^{p}wdx\right)^{1/p}
\]
for any Lipschitz function $f$ on $\Omega$. Where $B_{p,q}$ only
depends on $p,q$. 

For any $\delta\in(0,1)$, $u=\phi^{*}$ is Lipschitz over $\delta\cdot P\subset P$.
We apply the above inequality to datum $(\Omega,\rho,k,p,q,f)=(\delta P,\ell,1,q,q,u)$,
then there is a constant $C_{q}$ (only depending on $q$) such that
\[
\int_{\delta P}\left|u-\frac{\int_{\delta P}u\ell}{\int_{\delta P}\ell}\right|^{q}\ell\le C_{q}\cdot\mathrm{diam}(\delta P)^{q}\int_{\delta P}\left|du\right|^{q}\ell
\]
holds for any $\delta\in(0,1)$. By (\ref{eq: du^ql<C}), the RHS
is uniformly bounded, let $\delta\rightarrow1$ we know $\int_{P}u^{q}\ell<\infty$
for any $q\geq1$. 

Now we can get rid of $\ell.$ For any $r\geq1$, by the H\"{o}lder
inequality, 
\[
\int_{P}u^{r}=\int_{P}u^{r}\ell^{\frac{1}{p}}\ell^{-\frac{1}{p}}\le\left(\int_{P}u^{pr}\ell\right)^{\frac{1}{p}}\left(\int_{P}\ell^{-\frac{p'}{p}}\right)^{\frac{1}{p'}},
\]
where $p,p'>1$ are conjugate to each other. We take $\frac{p'}{p}=p'-1$
is sufficiently small such that $\int_{P}\ell^{-p'/p}<\infty$. Then
it follows that $\int_{P}u^{r}<\infty$ for all $r\geq1$, in particular
$\phi\in\mathcal{E}^{1}(P)$. Applying H\"{o}lder's inequality again
to (\ref{eq: du^ql<C}), we can show $\int_{P}\left|du\right|^{r}<\infty$
for all $r\geq1$. By the Sobolev imbedding theorem for domains with
Lipschitz boundary, we conclude that $u=\phi^{*}$ is H\"{o}lder
continuous over $P$ for any exponent $\gamma\in(0,1)$. In particular,
$\phi^{*}$ is bounded on $P$, so $\phi\in\textrm{PSH}_{b}(P)$. 

(5) Smoothness of $\phi$. When $\inf_{P}\ell>0$, since $\phi$ satisfies
$\mathrm{MA}(\phi)=\frac{e^{-\phi}}{\ell(\partial\phi)}$ (in the
sense of Alexandrov) and the RHS is locally bounded above and from
below away from $0$, the smoothness of $\phi$ directly follows from
the regularity results of Caffarelli, as in \cite{Berman Berndtsson}
section 2.12. When $\inf_{P}\ell=0$, we can not obtain a local upper
bound for $\mathrm{MA}(\phi)$ directly. But we can do this under
the assumption: $\dim P\cap\{\ell=0\}\leq\frac{n}{2}$ (including
$\inf_{P}\ell>0$), then we can show $\phi$ is smooth. 

Let $K\subset\mathbb{R}^{n}$ be a compact set, then $\partial\phi(K)\subset P$
is also compact. We will show that $\partial\phi(K)\cap\{\ell=0\}$
must be empty if $\dim P\cap\{\ell=0\}\leq\frac{n}{2}$. If this is
not true, there is $p\in\partial\phi(x_{0}),x_{0}\in K$ such that
$\ell(p)=0$. Let 
\[
S(p)\coloneqq\{\xi\in N_{\mathbb{R}}\mid\left\langle \xi,p\right\rangle \geq\left\langle \xi,y\right\rangle \ \textrm{for all}\ y\in P\},
\]
which is a cone of dimension $\geq n-\dim P\cap\{\ell=0\}$. We claim
$\phi$ is affine over the set $x_{0}+S(p)$. 

In fact, for any nonzero $\xi\in S(p)$, consider function $f(t)\coloneqq\phi(x_{0}+t\xi)$,
which is convex in $t\in\mathbb{R}$. By the definition of sub-differentials,
it is easy to see $\left\langle \xi,p\right\rangle \in\partial f(0)$.
Suppose $f$ is differentiable at $t_{0}>0$ (almost everywhere),
then for any $q\in\partial\phi(x_{0}+t_{0}\xi)\subset P$, we have
$f'(t_{0})=\left\langle \xi,q\right\rangle \leq\left\langle \xi,p\right\rangle $.
The monotonicity of slopes of $f$ implies $f$ must be affine on
$[0,\infty)$, i.e. $\phi(x_{0}+t\xi)=\phi(x_{0})+t\left\langle \xi,p\right\rangle $
for all $t\geq0$. The claim follows, more specifically we have 
\begin{equation}
\phi(x_{0}+\xi)=\phi(x_{0})+\left\langle \xi,p\right\rangle ,\ \textrm{for all}\ \xi\in S(p).\label{eq: phi is affine S}
\end{equation}

Next we show $\mathrm{MA}(\phi)$ is locally bounded from below. Let
$F\subset\mathbb{R}^{n}$ be any bounded subset, suppose $\phi\leq C$
on $F$ and $\ell\leq C'$ on $P$, then for any Borel subset $E\subset F$,
by (\ref{eq: Alexandrov-ell}), we have 
\[
C'\mathrm{MA}(\phi)(E)\geq\int_{\partial\phi(E)}\ell=\int_{E}e^{-\phi}\geq e^{-C}\left|E\right|,
\]
thus $\mathrm{MA}(\phi)\geq\frac{1}{C'}e^{-C}\ dx$ on $F$. 

Now we use the assumption: $\dim P\cap\{\ell=0\}\leq\frac{n}{2}$,
it implies $\dim S(p)\geq\frac{n}{2}$. Let 
\[
\psi(x)\coloneqq\phi(x)-\phi(x_{0})-\left\langle x-x_{0},p\right\rangle ,
\]
then $\psi\geq0$ since $p\in\partial\phi(x_{0})$ and $\psi$ vanishes
on $x_{0}+S(p)$ since (\ref{eq: phi is affine S}). Taking a small
ball $B$ intersecting $x_{0}+S(p)$, then $\mathrm{MA}(\psi)\geq\lambda\ dx$
on $B$ for some $\lambda>0$. But Lemma 2.3 in \cite{Mooney} says,
if $\mathrm{MA}(\psi)\geq\lambda\ dx$ on a ball, $\psi$ can not
vanish along a subspace of dimension $\geq\frac{n}{2}$. We obtain
a contradiction, thus $\partial\phi(K)\cap\{\ell=0\}$ must be empty
for any compact $K\subset\mathbb{R}^{n}$. 

Now we can locally bound $\mathrm{MA}(\phi)$ from both sides. For
any $R\gg1$, take $K=\Omega_{R}\coloneqq\{\phi\leq R\}$, which is
a compact convex set. Since we already show $\partial\phi(K)$ does
not intersect $\{\ell=0\}$, so $\inf_{\partial\phi(K)}\ell>0$. For
any Borel set $E\subset\Omega_{R}$, we have 
\[
\inf_{\partial\phi(K)}\ell\cdot\mathrm{MA}(\phi)(E)\leq\int_{\partial\phi(E)}\ell=\int_{E}e^{-\phi}\leq e^{-\min\phi}\left|E\right|.
\]
Hence we have $\lambda dx\leq\mathrm{MA}(\phi)\leq\Lambda dx$ on
$\Omega_{R}$ for some $\lambda,\ \Lambda>0$. By Corollary 4.11 \cite{Figalli},
we know $\phi$ is strictly convex on $\Omega_{R}$, thus on whole
$\mathbb{R}^{n}$. Note that strict convexity implies $\partial\phi(K)\cap\partial P=\emptyset$
for all compact $K\subset\mathbb{R}^{n}$. If this is not true, $\phi$
would be affine over a cone as the above $x_{0}+S(p)$. 

With the strict convexity and the bounds for $\mathrm{MA}(\phi)$
from both sides, by Caffarelli's interior $C^{1,\alpha}$ estimates
(see Corollary 4.21 \cite{Figalli}), we know $\phi\in C_{\mathrm{loc}}^{1,\alpha}(\Omega_{R})$
for some $\alpha\in(0,1)$ depending on $R$. Thus the RHS of $\mathrm{MA}(\phi)=\frac{e^{-\phi}}{\ell(\partial\phi)}$
belongs to $C^{0,\alpha}$ locally, then by Caffarelli's interior
$C^{2,\alpha}$ estimates (see Corollary 4.43 \cite{Figalli}), we
have $\phi\in C_{\mathrm{loc}}^{2,\alpha}(\Omega_{R})$. Finally,
the smoothness of $\phi$ follows from a standard bootstrapping argument. 
\end{proof}
\begin{rem}
\label{rem: Modified Guillemin}The solution $\phi$ obtained above
is expected to be also smooth when $\dim P\cap\{\ell=0\}>\frac{n}{2}$.
On the other hand, we make a conjecture about the asymptotics of $\phi^{*}$
near the boundary set $P\cap\{\ell=0\}$. 

Assume $P$ is a Delzant polytope. For a vertex on $\{\ell=0\}$,
by changing the coordinate, $P$ is locally given by 
\[
\{y\mid y_{i}\geq-1,\ 1\leq i\leq n\},\ \textrm{and}\ \ell(y)=\sum_{\alpha=1}^{p}c_{\alpha}(y_{\alpha}+1),\ c_{\alpha}>0.
\]
Consider a model potential $u(y)=\sum_{i=1}^{n}(y_{i}+1)\log(y_{i}+1)+\lambda\cdot\ell\log\ell$.
By a computation, in order to make the two sides of 
\[
\ell\cdot\det[u_{ij}]^{-1}=e^{u-y\cdot du}
\]
have matching vanish factors along $P\cap\{\ell=0\}$, $\lambda$
must be equal to $\ell(0)^{-1}$. Hence we conjecture $\phi^{*}$
satisfies a modified Guillemin condition, namely 
\[
\phi^{*}=\frac{1}{\ell(0)}\ell\cdot\log\ell+\sum_{\alpha\in\Lambda}l_{\alpha}\log l_{\alpha}+w,
\]
recall $l_{\alpha}(y)=\left\langle q_{\alpha},y\right\rangle +1$
are the defining functions of the facets of $P$, where $w(y)$ can
be extended to a smooth function on a neighborhood of $P$. Note this
form explains the coefficient \textbf{2} (so the cone angle $\pi$)
in (\ref{eq:Guillemin solu 8c}). 
\end{rem}

\section{The unstable case and moment-weight equalities}

\subsection{Toric moment-weight inequalities}

In this section, we derive a toric version of (\ref{eq: M-W Fano})
from Jensen's inequality. 
\begin{defn}
Let $P\subset\mathbb{R}^{n}$ be a convex body containing $0$ in
its interior. We define 
\[
\mathfrak{B}_{P}=\left\{ 0\leq b(y)\in L^{2}(P)\mid\int_{P}bdy=1;\ \int_{P}y_{i}\ bdy=0,\ \textrm{for}\ 1\leq i\leq n\right\} .
\]
We call $b(y)$ is a \textit{balancing density} on $P$, since the
barycenter of $(P,bdy)$ is the origin. 
\end{defn}

\begin{prop}
\label{prop: h is balancing}Let $X$ be a smooth toric Fano variety
associated to  Delzant polytope $P$. Then for any $\phi\in\mathcal{H}(-K_{X})^{S}$,
we have $V_{P}^{-1}e^{h(\phi)}(d\phi^{*})\in\mathfrak{B}_{P}$. 
\end{prop}

\begin{proof}
Recall $h(\phi)$ is the Ricci potential given by (\ref{eq: toric Ricci poten}).
Here we take $h(\phi)$ as a function defined on $N_{\mathbb{R}}$,
so $e^{h(\phi)}(d\phi^{*})$ is a function on $P$. By (\ref{eq: toric Ricci poten}),
we see $\int_{P}e^{h(\phi)}(d\phi^{*})dy=V_{P}$. For $1\leq i\leq n$,
we have ($C$ is a constant)
\[
\int_{P}y_{i}e^{h(\phi)}(d\phi^{*})dy=C\int_{\mathbb{R}^{n}}\partial_{i}\phi\ e^{-\phi}dx=-C\int_{\mathbb{R}^{n}}\partial_{i}(e^{-\phi})dx.
\]
Since $\phi(x)\geq\delta\left|x\right|-C'$ for some $\delta$, $C'>0$,
applying the divergence formula on the ball $B_{R}$ and then let
$R\rightarrow\infty$, we see the last integral is equal to zero. 
\end{proof}
In the below, we use the $L^{2}$-norm $\left\Vert f\right\Vert _{2}^{2}\coloneqq\int_{P}\left|f\right|^{2}dy$
and its variant $\left\Vert f\right\Vert _{2,d\bar{y}}^{2}\coloneqq\frac{1}{V_{P}}\int_{P}\left|f\right|^{2}dy$.
For simplicity, we will write $e^{h(\phi)}(d\phi^{*})$ as $e^{h(\phi)}$,
which is a function on $P$. 
\begin{prop}
(Toric Moment-Weight inequality) Let $X$ be a smooth toric Fano variety
associated to  Delzant polytope $P$. 

(1) Let $\bar{f}=\frac{1}{V_{P}}\int_{P}fdy$, then we have 
\begin{equation}
\inf_{\phi\in\mathcal{H}(P)}\left\Vert e^{h(\phi)}-1\right\Vert _{2,d\bar{y}}\geq\sup_{\textrm{convex}\ f}\frac{-D^{\mathrm{NA}}(f)}{\left\Vert f-\bar{f}\right\Vert _{2,d\bar{y}}},\label{eq: toric M-W ineq}
\end{equation}
where the supremum is taking over all non-constant convex functions
on $P$. 

(2) Let $\pi:L^{2}(P,dy)\rightarrow\left\{ l|_{P}\mid l\ \textrm{is affine}\right\} $
be the $L^{2}$-orthogonal projection, then we have 
\begin{equation}
\inf_{\phi\in\mathcal{H}(P)}\left\Vert e^{h(\phi)}-V_{P}\ell\right\Vert _{2,d\bar{y}}\geq\sup_{\textrm{convex}\ f}\frac{-D_{\ell}^{\mathrm{NA}}(f)}{\left\Vert f-\pi(f)\right\Vert _{2,d\bar{y}}},\label{eq: relative toric M-W ineq}
\end{equation}
where the supremum is taking over all non-affine convex functions
on $P$. 
\end{prop}

\begin{proof}
Let $g(y)$ be a continuous function on $P$ with $\int_{P}g=1$,
we will take $g=V_{P}^{-1}$ and $\ell_{P}$. For any $\phi\in\mathcal{H}(P)$
and convex function $f$ on $P$, we have 
\[
\int_{P}f\left(V_{P}^{-1}e^{h(\phi)}-g\right)dy=-f(0)+\frac{1}{V_{P}}\int_{P}fe^{h(\phi)}dy-D_{g}^{\textrm{NA}}(f)\geq-D_{g}^{\textrm{NA}}(f).
\]
The last step follows from Jensen's inequality with Proposition \ref{prop: h is balancing}. 

Now take $g=V_{P}^{-1}$, note that $\int_{P}f\left(e^{h(\phi)}-1\right)dy=\int_{P}\left(f-\bar{f}\right)\left(e^{h(\phi)}-1\right)dy$.
By the Cauchy-Schwarz inequality, we have 
\[
\int_{P}\left(f-\bar{f}\right)\left(e^{h(\phi)}-1\right)\frac{dy}{V_{P}}\leq\left(\int_{P}\left|f-\bar{f}\right|^{2}\frac{dy}{V_{P}}\right)^{1/2}\cdot\left(\int_{P}\left|e^{h(\phi)}-1\right|^{2}\frac{dy}{V_{P}}\right)^{1/2}.
\]
Then (\ref{eq: toric M-W ineq}) follows from the above two inequalities. 

Take $g=\ell_{P}$, by Proposition \ref{prop: h is balancing}, we
have $\int_{P}\pi(f)\left(e^{h(\phi)}-V_{P}\ell\right)dy=0$. Then
(\ref{eq: relative toric M-W ineq}) follows in the same way. 
\end{proof}

\subsection{The maximal destabilizers}

Let $P$ be a convex body containing $0$ in its interior. We define
\[
\mathcal{C}_{2}(P)\coloneqq\{f:P\rightarrow\mathbb{R}\cup\{\infty\}\ \textrm{is convex}\mid f\in L^{2}(P)\}.
\]
Note that for any $f\in\mathcal{C}_{2}(P)$, since it is convex and
$f<\infty$ almost everywhere, it must have $f<\infty$ on $\textrm{Int}(P)$. 

Let us consider the optimization problem on the RHS of (\ref{eq: toric M-W ineq}),
we need to minimize the functional 
\begin{equation}
W(f)\coloneqq\frac{D^{\textrm{NA}}(f)}{\left\Vert f-\bar{f}\right\Vert _{2}},\ \ f\in\mathcal{C}_{2}(P).\label{eq: W(f)}
\end{equation}
When the denominator is zero, we set $W=0$. If there is a minimizer,
then we take it as the maximal destabilizer for Ding stability. In
\cite{Gabor optimal}, Sz\'{e}kelyhidi had considered a similar problem
for K-stability. As in \cite{Gabor optimal}, we first try to minimize
\[
W_{\ell}(f)\coloneqq\frac{D_{\ell}^{\textrm{NA}}(f)}{\left\Vert f-\pi(f)\right\Vert _{2}},\ f\in\mathcal{C}_{2}(P).
\]
The advantage is that $W_{\ell}$ is invariant under addition by affine
functions. 
\begin{lem}
\label{lem: inf of Wl}Let $P\subset\mathbb{R}^{n}$ be a convex body
containing $0$ in its interior. Suppose that $\inf_{P}\ell_{P}<0$.
Then there is a unique minimizer for $W_{\ell}$, denoted by $\hat{f}\in\mathcal{C}_{2}(P)$
such that 
\begin{equation}
\pi(\hat{f})=0\ \textrm{and}\ \int_{P}\hat{f}^{2}=\hat{f}(0).\label{eq: normalize minimizer}
\end{equation}
The above two conditions are equivalent to require $\int_{P}\hat{f}l=0$
for all affine $l$ and $D_{\ell}^{\mathrm{NA}}(\hat{f})=-\left\Vert \hat{f}\right\Vert _{2}^{2}$. 
\end{lem}

\begin{proof}
Since $W_{\ell}$ is invariant under addition of affine functions,
we only need to minimize it over the normalized functions in $\mathcal{C}_{2}(P)$,
i.e. satisfying $f\geq f(0)=0$. 

Step 1: show $W_{\ell}$ is bounded below. By Lemma 7 \cite{Gabor optimal},
there exists $C_{P}>0$ only depending on $P$ such that $\left\Vert f\right\Vert _{2}\leq C_{P}\left\Vert f-\pi(f)\right\Vert _{2}$
for all normalized $f\in\mathcal{C}_{2}(P)$. It follows that 
\[
D_{\ell}^{\textrm{NA}}(f)=\int_{P}f\ell\geq-\left\Vert l\right\Vert _{2}\left\Vert f\right\Vert _{2}\geq-C_{P}\left\Vert l\right\Vert _{2}\left\Vert f-\pi(f)\right\Vert _{2},
\]
thus $W_{\ell}\geq-C_{P}\left\Vert l\right\Vert _{2}$. 

Step 2: Let $\{f_{k}\}\subset\mathcal{C}_{2}(P)$ be a sequence of
normalized convex functions such that $W_{\ell}(f_{k})\rightarrow\inf_{\mathcal{C}_{2}(P)}W_{\ell}<0$.
Rescaling them such that $\int_{P}f_{k}^{2}=1$, by H\"{o}lder's
inequality, $\int_{P}f_{k}$ are bounded above. By Corollary 5.2.5
\cite{Dona toric}, there is a subsequence (still denoted by $f_{k}$)
converges to a function $\hat{f}:P\rightarrow\mathbb{R}\cup\{\infty\}$
uniformly over compact subsets of $\mathrm{Int}(P)$. Clearly, $\hat{f}$
is also a normalized convex function. By Fatou's lemma, we have $\int_{P}\hat{f}^{2}\leq1$
thus $\hat{f}\in\mathcal{C}_{2}(P)$. 

Step 3: show $\hat{f}$ is a minimizer. By the weak compactness in
$L^{2}(P)$, after passing to a subsequence, we can assume $f_{k}$
weakly converges in $L^{2}(P)$. The weak limit must be $\hat{f}$.
It implies 
\[
D_{\ell}^{\textrm{NA}}(f_{k})=\int_{P}f_{k}\ell\rightarrow\int_{P}\hat{f}\ell=D_{\ell}^{\textrm{NA}}(\hat{f}),
\]
thus 
\begin{equation}
\left\Vert f_{k}-\pi(f_{k})\right\Vert _{2}=\frac{D_{\ell}^{\textrm{NA}}(f_{k})}{W_{\ell}(f_{k})}\rightarrow\frac{D_{\ell}^{\textrm{NA}}(\hat{f})}{\inf W_{\ell}}.\label{eq: minimize relation}
\end{equation}
The weak convergence also implies $\pi(f_{k})\rightarrow\pi(\hat{f})$
in the space of affine functions. 

If $D_{\ell}^{\textrm{NA}}(\hat{f})=0$, (\ref{eq: minimize relation})
implies $\left\Vert f_{k}-\pi(f_{k})\right\Vert _{2}\rightarrow0$,
but this contradicts with $\int_{P}f_{k}^{2}=1$ and $\left\Vert f_{k}\right\Vert _{2}\leq C_{P}\left\Vert f_{k}-\pi(f_{k})\right\Vert _{2}$.
Hence $D_{\ell}^{\textrm{NA}}(\hat{f})<0$ and so $\lim_{k}\left\Vert f_{k}-\pi(f_{k})\right\Vert _{2}>0$.
By (\ref{eq: minimize relation}), we have 
\[
\frac{D_{\ell}^{\textrm{NA}}(\hat{f})}{\lim_{k}\left\Vert f_{k}-\pi(f_{k})\right\Vert _{2}}=\inf W_{\ell}.
\]
Once we show 
\begin{equation}
\lim_{k}\left\Vert f_{k}-\pi(f_{k})\right\Vert _{2}\geq\left\Vert \hat{f}-\pi(\hat{f})\right\Vert _{2},\label{eq: step3 inequ}
\end{equation}
then it implies $W_{\ell}(\hat{f})\leq\inf W_{\ell}$, since $\hat{f}\in\mathcal{C}_{2}(P)$,
thus $\hat{f}$ is a minimizer. In fact, for any compact subset $K\subset\textrm{Int}(P)$,
since $f_{k}$ uniformly converges to $\hat{f}$ on $K$, we have
\[
\left\Vert \hat{f}-\pi(\hat{f})\right\Vert _{2,K}=\lim_{k}\left\Vert f_{k}-\pi(f_{k})\right\Vert _{2,K}\leq\lim_{k}\left\Vert f_{k}-\pi(f_{k})\right\Vert _{2},
\]
where $\left\Vert \cdot\right\Vert _{2,K}$ is the $L^{2}$-norm over
$K$. Then we take $\sup_{K}$ of the LHS, (\ref{eq: step3 inequ})
follows. Hence $\hat{f}$ is a minimizer. 

Next we adjust $\hat{f}$ to satisfy (\ref{eq: normalize minimizer}).
First we replace $\hat{f}$ by $\hat{f}-\pi(\hat{f})$, which is still
a minimizer. After that, we have $\pi(\hat{f})=0$, equivalently $\int_{P}\hat{f}l=0$
for all affine $l$. Thus $\hat{f}(0)=-D_{\ell}^{\textrm{NA}}(\hat{f})>0$.
Finally, we rescale $\hat{f}$ such that $\int_{P}\hat{f}^{2}=\hat{f}(0)$.
Now we obtain a minimizer $\hat{f}$ satisfying (\ref{eq: normalize minimizer}).
The statement for the equivalence of two normalization conditions
is trivial. 

Step 4: the uniqueness. Suppose $\hat{g}\in\mathcal{C}_{2}(P)$ is
another minimizer satisfying (\ref{eq: normalize minimizer}), then
$W_{\ell}(\hat{f})=W_{\ell}(\hat{g})$ and 
\[
W_{\ell}(\hat{g})=\frac{-\hat{g}(0)}{\left\Vert \hat{g}\right\Vert _{2}}=-\left\Vert \hat{g}\right\Vert _{2}.
\]
It implies $\left\Vert \hat{f}\right\Vert _{2}=\left\Vert \hat{g}\right\Vert _{2}$
and $D_{\ell}^{\textrm{NA}}(\hat{f})=D_{\ell}^{\textrm{NA}}(\hat{g})<0$.
Since $\hat{g}\neq\hat{f}$, $\hat{g}$ can not be proportional to
$\hat{f}$, thus we have $\left\Vert \hat{f}+\hat{g}\right\Vert _{2}<2\left\Vert \hat{f}\right\Vert _{2}$.
It follows that 
\[
W_{\ell}(\hat{f}+\hat{g})=\frac{2D_{\ell}^{\textrm{NA}}(\hat{f})}{\left\Vert \hat{f}+\hat{g}\right\Vert _{2}}<\frac{2D_{\ell}^{\textrm{NA}}(\hat{f})}{2\left\Vert \hat{f}\right\Vert _{2}}=W_{\ell}(\hat{f})=\inf_{\mathcal{C}_{2}(P)}W_{\ell}.
\]
It is a contradiction, since $\hat{f}+\hat{g}\in\mathcal{C}_{2}(P)$. 
\end{proof}
Next we study the structure of the minimizer $\hat{f}$ and give it
more characterizations. 
\begin{thm}
\label{thm: optimal balan density}Let $P\subset\mathbb{R}^{n}$ be
a convex body containing $0$ in its interior. Suppose that $\inf_{P}\ell_{P}<0$.
Let $\hat{f}$ be the unique minimizer for $W_{\ell}$ satisfying
(\ref{eq: normalize minimizer}), which is obtained in the above lemma.
Then we have (1)-(4) below. 

(1) Let $\mathfrak{b}\coloneqq\hat{f}+\ell_{P}$, then $\mathfrak{b}\geq0$,
and $\int_{P}\mathfrak{b}=1$, $\int_{P}y_{i}\mathfrak{b}=0$ for
$1\leq i\leq n$, hence $\mathfrak{b}\in\mathfrak{B}_{P}$. Moreover,
$\int_{P}\mathfrak{b}^{2}=\mathfrak{b}(0)$. 

(2) $\mathfrak{b}$ is the unique element in $\mathfrak{B}_{P}$ with
the minimum $L^{2}$-norm. 

(3) For any convex $b\in\mathfrak{B}_{P}$, we have $\int_{P}b^{2}\geq b(0)$.
The equality holds if and only if $b=\mathfrak{b}$. 

(4) $\mathfrak{b}$ is a simple convex function, i.e. $\mathfrak{b}=\max\{0,\mathfrak{a}\}$,
where $\mathfrak{a}$ is the unique affine function such that 
\begin{equation}
\int_{P\cap\{\mathfrak{a}>0\}}\mathfrak{a}=1,\int_{P\cap\{\mathfrak{a}>0\}}y_{i}\mathfrak{a}=0,\ \textrm{for}\ 1\leq i\leq n.\label{eq: optimal affine}
\end{equation}

(5) If $\inf_{P}\ell_{P}\geq0$, then $\ell_{P}$ is the unique element
in $\mathfrak{B}_{P}$ with the minimum $L^{2}$-norm. 
\end{thm}

\begin{proof}
(1) For any $h\in\mathcal{C}_{2}(P)$, we have $f_{t}\coloneqq\hat{f}+th\in\mathcal{C}_{2}(P)$
for $t\geq0$. Since $\hat{f}$ is a minimizer for $W_{\ell}$, thus
$\frac{d}{dt}W_{\ell}(f_{t})|_{t=0}\geq0$. Combing with (\ref{eq: normalize minimizer}),
it follows that 
\begin{equation}
-h(0)+\int_{P}h(\ell+\hat{f})\geq0,\ \textrm{for all}\ h\in\mathcal{C}_{2}(P).\label{eq: variational inequ}
\end{equation}
Let $\mathfrak{b}\coloneqq\ell+\hat{f}$, then (\ref{eq: variational inequ})
is equivalent to $D_{\mathfrak{b}}^{\textrm{NA}}(h)\geq0$ for all
$h\in\mathcal{C}_{2}(P)$. In particular, $D_{\mathfrak{b}}^{\textrm{NA}}(l)=0$
for all affine $l$. Then we have 
\[
D_{\mathfrak{b}}^{\textrm{NA}}(\mathfrak{b})=D_{\mathfrak{b}}^{\textrm{NA}}(\hat{f})=-\hat{f}(0)+\int_{P}\hat{f}\mathfrak{b}=-\hat{f}(0)+\int_{P}\hat{f}^{2}=0,
\]
where we used (\ref{eq: normalize minimizer}). This is exactly $\int_{P}\mathfrak{b}^{2}=\mathfrak{b}(0)$.
Finally, we show $\mathfrak{b}\geq0$. If it is not true, let $\mathfrak{b}^{+}=\max\{\mathfrak{b},0\}\in\mathcal{C}_{2}(P)$,
we have 
\[
D_{\mathfrak{b}}^{\textrm{NA}}(\mathfrak{b}^{+})=-\mathfrak{b}(0)+\int_{P}\mathfrak{b}\mathfrak{b}^{+}=-\int_{P}\mathfrak{b}^{2}+\int_{\{\mathfrak{b}\geq0\}}\mathfrak{b}^{2}=-\int_{\{\mathfrak{b}<0\}}\mathfrak{b}^{2}<0,
\]
note that $\mathfrak{b}$ is continuous on $\mathrm{Int}(P)$. This
contradicts to (\ref{eq: variational inequ}), so we have $\mathfrak{b}\geq0$. 

(2) For any $b\in\mathfrak{B}_{P}$, by Jensen's inequality with measure
$bdy$, we have 
\[
\left\Vert \mathfrak{b}\right\Vert ^{2}=\mathfrak{b}(0)\leq\int_{P}\mathfrak{b}\cdot bdy\leq\left\Vert \mathfrak{b}\right\Vert \cdot\left\Vert b\right\Vert ,
\]
thus $\left\Vert \mathfrak{b}\right\Vert \leq\left\Vert b\right\Vert $.
Hence $\mathfrak{b}$ has the minimum $L^{2}$-norm. If $\left\Vert \mathfrak{b}\right\Vert =\left\Vert b\right\Vert $,
then $b$ must be proportional to $\mathfrak{b}$. Since $\int_{P}b=\int_{P}\mathfrak{b}=1$,
it forces $b=\mathfrak{b}$. 

(3) For any convex $b\in\mathfrak{B}_{P}$, by Jensen's inequality
with measure $bdy$, we have $\int_{P}b^{2}=\int_{P}b\cdot bdy\geq b(0)$.
If the equality holds, again by Jensen's inequality, 
\[
\left\Vert b\right\Vert ^{2}=b(0)\leq\int_{P}b\cdot\mathfrak{b}dy\leq\left\Vert b\right\Vert \cdot\left\Vert \mathfrak{b}\right\Vert ,
\]
thus $\left\Vert b\right\Vert \leq\left\Vert \mathfrak{b}\right\Vert $.
Since $\mathfrak{b}$ has the minimum $L^{2}$-norm, (2) tells us
$b=\mathfrak{b}$. 

(4) Note that $\mathfrak{b}dy$ is a positive Radon measure with barycenter
$0$. Since $\mathfrak{b}$ is convex, the condition $\int_{P}\mathfrak{b}\cdot\mathfrak{b}dy=\mathfrak{b}(0)$
means that Jensen's inequality with measure $\mathfrak{b}dy$ takes
equality. By the equality condition of Jensen's inequality, $\mathfrak{b}$
must be affine over the convex hull of set $\{\mathfrak{b}>0\}$. 

We claim that $\{\mathfrak{b}>0\}$ is convex. In fact, for any $y_{1},y_{2}\in\{\mathfrak{b}>0\}$,
connecting them by a line segment $L$. Since $L$ is contained in
the convex hull of $\{\mathfrak{b}>0\}$, so $\mathfrak{b}$ must
be affine along $L$, this implies $L\subset\{\mathfrak{b}>0\}$,
thus $\{\mathfrak{b}>0\}$ is convex. In particular, $\{\mathfrak{b}>0\}$
is connected. The only possibility is $\mathfrak{b}=\max\{0,\mathfrak{a}\}$
for some affine function $\mathfrak{a}$. Moreover, since $\mathfrak{b}\in\mathfrak{B}_{P}$,
so $\mathfrak{a}$ should satisfy (\ref{eq: optimal affine}). 

For the uniqueness of $\mathfrak{a}$, if there is an affine function
$l$ satisfying (\ref{eq: optimal affine}) (replacing $\mathfrak{a}$
by $l$), then we have $\int_{P\cap\{l>0\}}l^{2}=l(0)$. Then $b\coloneqq\max\{0,l\}\in\mathfrak{B}_{P}$
and satisfies $\int_{P}b^{2}=b(0)$, by (3) we know $b=\mathfrak{b}$,
so $l=\mathfrak{a}$. 

(5) The proof is same to (2). 
\end{proof}
The characteristic condition (\ref{eq: optimal affine}) for $\mathfrak{a}$
can be rephrased in a simple way, that is the Ricci affine function
associated to $P\cap\{\mathfrak{a}\geq0\}$ is $\mathfrak{a}$ itself.
It gives a polynomial system with $n+1$ equations for the coefficients
of $\mathfrak{a}$, the above theorem ensures that it admits a unique
solution. 
\begin{defn}
Let $P\subset\mathbb{R}^{n}$ be a convex body containing $0$ in
its interior, and $\ell_{P}$ is the associated Ricci affine function.
We define the \textit{optimal balancing density }of $P$ to be 
\[
\mathfrak{b}=\max\{0,\mathfrak{a}\},\ \textrm{when}\ \inf_{P}\ell_{P}<0;\ \ \mathfrak{b}=\ell_{P},\ \textrm{when}\ \inf_{P}\ell_{P}\geq0.
\]
Where $\mathfrak{a}$ is the unique affine function satisfying (\ref{eq: optimal affine}).
By Theorem \ref{thm: optimal balan density}, in any case, $\mathfrak{b}$
is the unique element in $\mathfrak{B}_{P}$ with the minimum $L^{2}$-norm.
We always have $\int_{P}\mathfrak{b}^{2}=\mathfrak{b}(0)$. 
\end{defn}

\begin{cor}
\label{cor: inf of W}Let $P\subset\mathbb{R}^{n}$ be a convex body
containing $0$ in its interior. Then we have 
\[
\inf\left\{ W(f)\mid f\in\mathcal{C}_{2}(P)\right\} =-\left\Vert V_{P}^{-1}-\mathfrak{b}\right\Vert _{2},
\]
and the infimum is attained by $\mathfrak{b}$. Recall $W(f)$ is
defined by (\ref{eq: W(f)}). 
\end{cor}

\begin{proof}
For any $f\in\mathcal{C}_{2}(P)$, since $\mathfrak{b}\in\mathfrak{B}_{P}$,
by Jensen's inequality we have 
\[
D^{\textrm{NA}}(f)=D_{\mathfrak{b}}^{\textrm{NA}}(f)+\int_{P}f(V_{P}^{-1}-\mathfrak{b})dy\geq-\left\Vert f-\bar{f}\right\Vert _{2}\cdot\left\Vert V_{P}^{-1}-\mathfrak{b}\right\Vert _{2}.
\]
Thus $W(f)\geq-\left\Vert V_{P}^{-1}-\mathfrak{b}\right\Vert _{2}$.
On the other hand, taking $f=\mathfrak{b}$, since $\int_{P}\mathfrak{b}^{2}=\mathfrak{b}(0)$,
we have $D^{\textrm{NA}}(\mathfrak{b})=-\left\Vert V_{P}^{-1}-\mathfrak{b}\right\Vert _{2}^{2}$.
Thus $W(\mathfrak{b})=-\left\Vert V_{P}^{-1}-\mathfrak{b}\right\Vert _{2}$,
all statements follow. 
\end{proof}
\begin{example}
\label{exa: optimal 6d}For the unstable polygon 6D in Example \ref{exa: 2 polygons},
\[
P=\{(y_{1},y_{2})\mid y_{1}>-1,\ y_{2}<1,\ 3y_{1}-2y_{2}<1\}.
\]
The associated toric Fano surface is $X_{P}=\mathbb{P}(1,2,3)$. By
a computation, we have $\ell_{P}=\frac{2}{3}y_{1}-\frac{2}{9}y_{2}+\frac{5}{9}$.
To determine the affine function $\mathfrak{a}$, firstly we need
to know the rough location of $\{\mathfrak{a}=0\}$. For this, we
use the following iteration scheme. Let $\ell_{0}\coloneqq\ell_{P}$,
then for any $k\geq1$, let $\ell_{k}$ be the Ricci affine function
of the polytope $P\cap\{\ell_{k-1}\geq0\}$. After several steps,
when the set $\{\ell_{k}=0\}$ is approaching to be stable, we solve
the polynomial system (\ref{eq: optimal affine}). By computer, we
find 
\[
\mathfrak{a}=\frac{1}{162}\left(67+26\sqrt{7}\right)y_{1}-\frac{2}{243}\left(29+4\sqrt{7}\right)y_{2}+\frac{7}{486}\left(29+4\sqrt{7}\right).
\]
In Figure \ref{fig: 2-polygons}, the set $\{\mathfrak{a}=0\}$ is
showed by a dashed line. Note that the slope of this line is an irrational
number, so $P\cap\{\mathfrak{a}\geq0\}$ can not be the moment polytope
of some toric variety. 
\end{example}

\begin{example}
\label{exam: optimal 3fold}Consider smooth toric Fano 3-fold $X=\mathbb{P}\left(\mathcal{O}_{\mathbb{P}^{2}}\oplus\mathcal{O}_{\mathbb{P}^{2}}(2)\right)$
(with label $\mathcal{B}_{1}$ in \cite{Yotsutani 4fold}), which
is associated to polytope 
\[
P=\{y\in\mathbb{R}^{3}\mid y_{1}\ge-1,\ y_{2}\ge-1,\ y_{1}+y_{2}+2y_{3}\leq1,\ -1\leq y_{3}\leq1\}
\]
with vertices $(0,-1,1),(-1,0,1),(-1,-1,1),(4,-1,-1),(-1,4,-1),(-1,-1,-1)$.
It is also the blowup of $\mathbb{P}(1,1,1,2)$ at the unique singularity
point $[0,0,0,1]$. By a computation, we have $V_{P}=\frac{31}{3}$
and $\ell_{P}=\frac{60}{349}y_{3}+\frac{57}{349}$. Thus $(-K_{X})^{3}=3!V_{P}=62$.
We see $X$ is relatively D-unstable. 

Next we determine the optimal balancing density. By the symmetry of
$P$, we assume $\mathfrak{a}=a(y_{3}+b)$. The condition (\ref{eq: optimal affine})
requires $\int_{P\cap\{y_{3}>-b\}}y_{i}(y_{3}+b)dy=0$, for $i=1,2,3$.
By solving any one equation, we obtain 
\[
b=\frac{1}{6}\left(30\left(9+\sqrt{51}\right)\right)^{1/3}+\frac{1}{6}\left(30\left(9-\sqrt{51}\right)\right)^{1/3}-1\approx0.9455.
\]
Then by another condition $\int_{P\cap\{y_{3}>-b\}}a(y_{3}+b)dy=1$,
we know $a$ is equal to 
\[
72\left(140+(2650850-9150\sqrt{51})^{1/3}+5^{2/3}(106034+366\sqrt{51})^{1/3}\right)^{-1}\approx0.1728.
\]
\end{example}

\subsection{The moment-weight equalities}
\begin{thm}
\label{thm: M-W equality}Let $X$ be a smooth toric Fano variety
associated to Delzant polytope $P$. Then (\ref{eq: toric M-W ineq})
and (\ref{eq: relative toric M-W ineq}) take equality. Moreover,
we have the moment-weight equality: 
\begin{equation}
\inf_{\omega\in c_{1}(X)}\left\Vert 1-e^{h_{\omega}}\right\Vert _{L^{2}(\omega^{n})}=\sup_{(\mathcal{X},\mathcal{L})}\frac{-D^{\mathrm{NA}}(\mathcal{X},\mathcal{L})}{\left\Vert \mathcal{X},\mathcal{L}\right\Vert _{2}},\label{eq: M-W equality}
\end{equation}
and both sides are equal to $\left(V_{P}^{-1}\int_{P}(V_{P}\mathfrak{b}-1)^{2}\right)^{1/2}=\left(V_{P}\mathfrak{b}(0)-1\right)^{1/2}$. 
\end{thm}

\begin{proof}
We will use norm $\left\Vert \cdot\right\Vert _{2}$, and replace
it by $\left\Vert \cdot\right\Vert _{2,d\bar{y}}$ in the end. By
Corollary \ref{cor: inf of W}, the supremum in (\ref{eq: toric M-W ineq})
is equal to $\left\Vert V_{P}^{-1}-\mathfrak{b}\right\Vert _{2}$
and attained by $\mathfrak{b}$. For (\ref{eq: relative toric M-W ineq}),
by Lemma \ref{lem: inf of Wl}, the supremum is equal to $\left\Vert \mathfrak{b}-\ell_{P}\right\Vert _{2}$
and also attained by $\mathfrak{b}$. We take a sequence $\{b_{j}\}\subset\mathfrak{B}_{P}$
such that $b_{j}$ are positive and smooth on a neighborhood of $P$
and $b_{j}\rightarrow\mathfrak{b}$ in $L^{2}(P)$. By Theorem 1.1
in \cite{Berman Berndtsson} (or the continuation method in \cite{Wang-Zhu}),
there exists $\phi_{j}\in\mathcal{H}(P)$ such that $e^{h(\phi_{j})}(d\phi_{j}^{*})=V_{P}\cdot b_{j}$.
Then for both inequalities, the energy of $\phi_{j}$ on the LHS converges
to the supremum on the RHS, thus they take equality. 

For (\ref{eq: M-W equality}), since we already have the inequality
for (\ref{eq: M-W equality}), hence 
\[
\inf_{\textrm{inv.}\ \omega}Din(\omega)^{\frac{1}{2}}\geq\inf_{\omega\in c_{1}(X)}Din(\omega)^{\frac{1}{2}}\geq\sup_{(\mathcal{X},\mathcal{L})}\frac{-D^{\textrm{NA}}(\mathcal{X},\mathcal{L})}{\left\Vert \mathcal{X},\mathcal{L}\right\Vert _{2}}\geq\sup_{\textrm{convex}\ f}\frac{-D^{\textrm{NA}}(f)}{\left\Vert f-\bar{f}\right\Vert _{2,d\bar{y}}},
\]
the first infimum is taking over all $S_{N}$-invariant metric in
$c_{1}(X)$. The last $\geq$ is since the two quotients are the same
when $(\mathcal{X},\mathcal{L})$ is a toric test-configuration induced
by $f$. Since (\ref{eq: toric M-W ineq}) takes equality, it forces
(\ref{eq: M-W equality}) holds. 
\end{proof}

\end{document}